\documentclass[10pt]{amsart}
\usepackage{amsmath}
\usepackage{amsfonts}
\usepackage{amsthm}
\usepackage{amssymb}
\usepackage{amscd}
\usepackage[all]{xy}
\usepackage{enumerate}


\theoremstyle{plain}
\newtheorem{thm}{Theorem}[subsection]

\newtheorem{prop}[thm]{Proposition}

\newtheorem{cor}[thm]{Corollary}
\newtheorem{lem}[thm]{Lemma}
\newtheorem{cla}[thm]{Claim}

\theoremstyle{definition}
\newtheorem{defn}[thm]{Definition}

\newtheorem{prob}[thm]{Problem}

\newtheorem{nota}[thm]{Notation}

\newtheorem*{ackn}{Acknowledgment}

\theoremstyle{remark}
\newtheorem*{rem}{Remark}
\newcommand{\sB}{\mathcal{B}}
\newcommand{\sC}{\mathcal{C}}
\newcommand{\sD}{\mathcal{D}}
\newcommand{\sE}{\mathcal{E}}
\newcommand{\sF}{\mathcal{F}}
\newcommand{\sG}{\mathcal{G}}
\newcommand{\sH}{\mathcal{H}}
\newcommand{\sI}{\mathcal{I}}
\newcommand{\sJ}{\mathcal{J}}
\newcommand{\sL}{\mathcal{L}}
\newcommand{\sN}{\mathcal{N}}
\newcommand{\sM}{\mathcal{M}}
\newcommand{\sO}{\mathcal{O}}

\newcommand{\sQ}{\mathcal{Q}}

\newcommand{\sU}{\mathcal{U}}
\newcommand{\sV}{\mathcal{V}}

\newcommand{\sZ}{\mathcal{Z}}
\newcommand{\tA}{{\widetilde{A}}}

\newcommand{\mC}{\mathbb{C}}
\newcommand{\mF}{\mathbb{F}}

\newcommand{\mN}{\mathbb{N}}
\newcommand{\mP}{\mathbb{P}}

\newcommand{\mZ}{\mathbb{Z}}

\newcommand{\Ima}{\mathrm{Im}\,}

\newcommand{\ap}{\mathrm{ap}}
\newcommand{\Bs}{\mathrm{Bs}\,}

\newcommand{\CP}{\mathrm{CP}}

\newcommand{\expdim}{\mathrm{expdim}\,}

\newcommand{\Aut}{\mathrm{Aut}\,}
\newcommand{\Hilb}{\mathrm{Hilb}}
\newcommand{\Hom}{\mathrm{Hom}}
\newcommand{\sHom}{\mathcal{H}{\!}om\,}

\newcommand{\Pic}{\mathrm{Pic}\,}

\newcommand{\OG}{\mathrm{OG}}

\newcommand{\Supp}{\mathrm{Supp}\,}

\newcommand{\VSP}{\mathrm{VSP}\,}

\numberwithin{equation}{section}

\title[Scorza quartics of trigonal spin curve 
and their varieties of power sums 
]{Scorza quartics of trigonal spin curves and their varieties of power sums}

\author{Hiromichi Takagi}

\address{Graduate School of Mathematical Sciences \\
the University of Tokyo\\
Tokyo, 153-8914, Japan\\
\texttt{takagi@ms.u-tokyo.ac.jp}}

\author{Francesco Zucconi}

\address{D.I.M.I. \\
the University of Udine\\
Udine, 33100 Italy\\
\texttt{Francesco.Zucconi@dimi.uniud.it}}

\date{1.11, 2007}

\begin{document}

\maketitle

\markboth{Takagi and Zucconi}{Varieties of power sums}

\begin{abstract}
Our fundamental result is the
construction of new subvarieties in the varieties of power sums for
the Scorza quartic of any general
pairs of trigonal curves and non-effective theta characteristics.
This is a generalization of Mukai's description of 
smooth prime Fano threefolds of genus twelve as
the varieties of power sums for plane quartics. 
Among other applications, we give 
an affirmative answer to the conjecture of Dolgachev and Kanev on the existence of the Scorza quartic for any general pairs of
curves and non-effective theta characteristics. 
\end{abstract}

\tableofcontents

\section{Introduction}

\subsection{Varieties of power sums and the Waring problem}
\label{subsection:temptative}~

Throughout the paper, we work over $\mC$, the complex number field.

The problem of representing a homogeneous form as a sum of powers of
linear forms has been studied since the last decades of the
$19^{\rm{th}}$ century.  This is called the {\em{Waring problem}} for
a homogeneous form.  We are interested in the study of the global
structure of a suitable compactification of the variety parameterizing
all such representations of a homogeneous form.  To give a precise
definition of such a compactification consider a $(v+1)$-dimensional
vector space $V$.  Let $F\in S^m \check{V}$ be a homogeneous forms of
degree $m$ on $V$, where $\check{V}$ is the dual vector space of $V$.
Let $\mP_*\check{V}$ be the projective space parameterizing
one-dimensional vector subspaces in $\check{V}$, which is sometime
denoted by $\check{\mP}^v$.

\begin{defn}
    The {\em{varieties of power sums}} of $F$ is the following set
    with reduced structure:
    
    \[
    \VSP(F,n):=
    \overline{
    \{([H_1],\dots, [H_n])\mid H_1^m+\cdots+H_n^m=F\}}
    \subset \Hilb^n (\mP_*\check{V}).
    \]
    We call the {\em{Waring rank of $F$}}
    the minimum of $n$ such that $\VSP(F,n)\not =\emptyset$.
\end{defn}
There are other compactifications, for example, the one
in the $n$-th symmetric product of $\mP_*\check{V}$, but for our
treatment we need the one in the Hilbert scheme.

As far as we know,
the first global descriptions of positive dimensional
varieties of power sums for some homogeneous forms
were given by S. Mukai.

The most intensively studied cases of
varieties of power sums, including
Mukai's case, are where
$F$ is a general $(v+1)$-nary homogeneous form of degree $m$
for some $m, v\in \mN$, and 
$n$ is the Waring rank of $F$,
which we denote by $n(m,v)$.

By a standard parameter count,
we can easily compute the expected dimension
of $\VSP(F,n)$ for a general homogeneous form $F$. 
Since 
the dimension of the vector space of 
$(v+1)$-nary homogeneous forms of degree $m$ is
$\binom{m+v}{m}$,
the expected dimension is 
\[
\expdim \VSP(F,n):=
(v+1)n-\binom{m+v}{m}.
\]

Thus

it is expected that
\[
n(m,v)=\lceil \frac{1}{v+1}{\binom{m+v}{m}}\rceil.
\]
It is known, however,
that there are exceptions to $n(m,v)$ by the following result of
J. Alexander and A. Hirschowitz \cite{AH}:
\begin{center}
\begin{tabular}{|c|c|c|}
\hline
$m$ & $v$ & $n(m,v)$ \\
\hline\hline
$2$ & arbitrary & $v+1$ \\
\hline
$3$ &  $4$ & $8$ \\
\hline
$4$ &  $2$ & $6$ \\
\hline
$4$ &  $3$ & $10$ \\
\hline
$4$ &  $4$ & $15$ \\
\hline
\end{tabular}
\end{center}

Here is the table of the known descriptions of $\VSP(F,n(m,v))$.
\begin{flushleft}
{
\begin{tabular}{|c|c|c|c|c|}
\hline
$m$ & $v$ & $n(m,v)$ & $\VSP(F,n(m,v))$ & Ref. \\
\hline\hline
$2a-1$ &  $1$ & $a$ & $1$ point & Sylvester \\
\hline
$2$ &  $2$ & $3$ & quintic del Pezzo threefold & Mukai \cite{Mu2} \\
\hline
$3$ & $2$  & $4$  & $\mP^2$ & Dolgachev and Kanev \cite{DK} \\
\hline
$4$ & $2$  & $6$  & prime Fano threefold of genus twelve & \cite{Mu2} \\
\hline
$5$ & $2$  & $7$  & $1$ point & Hilbert, Richmond, Palatini \\
\hline
$6$ & $2$  & $10$  & polarized $K3$ surface of genus $20$ & \cite{Mu2} \\
\hline
$7$ & $2$  & $12$  & $5$ points & Dixon and Stuart \\
\hline
$8$ & $2$  & $15$  & $16$ points & \cite{Mu2} \\
\hline
$2$ & $3$  & $4$  & $G(2,5)$ & Ranestad and Schreier \cite{VSP} \\
\hline
$3$ & $3$  & $5$  & $1$ point & Sylvester's Pentahedral Theorem \\
\hline
$3$ & $4$  & $8$  & $W$ & \cite{VSP} \\
\hline
$3$ & $5$  & $10$  & $S$ & Iliev and Ranestad \cite{VSP2} \\
\hline
\end{tabular}}
\end{flushleft}
In the table, 
\begin{itemize}
\item $W$ is a fivefold 
and is the variety of lines in the
fivefold linear complete intersection
$\mP^{10}\cap \OG(5,10)\subset \mP^{15}$ of
the ten-dimensional orthogonal Grassmaniann $\OG(5,10)$,
\item
$S$ is a smooth symplectic fourfold obtained
as a deformation of the Hilbert square of a polarized $K3$ surface
of genus eight, and
\item
see the introduction of \cite{VSP} or \cite{doldual}
for the references of the results 
in the $19^{\rm{th}}$ and early $20^{\rm{th}}$ centuries.
\end{itemize}

As we can see in the table,
the study before Mukai's one were devoted only to the cases
where $\dim \VSP(F,n(m,v))=0$ and mostly the cases
where $F$ has a unique representation.
Recently, using the technique of birational geometry,
M. Mella proved in \cite{Waring}
that, if $m>v>1$,
then the uniqueness holds only in the
case where $(m,v)=(5,2)$.

In \cite{VSP3}, Iliev and Ranestad treat
some special $(v+1)$-nary cubics $F$ and prove that,
if $v\geq 8$, then 
the Waring rank of $F$ is less than that of a general cubic.

In \cite{IK} Iliev and Kanev study varieties of power sums
more systematically.

\subsection{Mukai's contribution}
\label{subsection:Mukai}~

Let $V_{22}$ be a smooth prime Fano threefold of genus twelve,
namely, a smooth projective threefold
such that $-K_{V_{22}}$ is ample, 
the class of $-K_{V_{22}}$ generates $\Pic V_{22}$, 
and the genus $g(V_{22}):=\frac{(-K_{V_{22}})^3}{2}+1$ is equal to twelve.
$V_{22}$ can be embedded into $\mP^{13}$ by the linear system $|-K_{V_{22}}|$. 
Mukai discovered the following remarkable result \cite[\S 6, Theorem 11]{Mu2}
(see also \cite{DK}, \cite{Schr}, and \cite[Theorem 3.12]{doldual}
for some details):
\begin{thm}\label{v22}
For a general ternary quartic form $F_4$, 
${\rm{VPS}}(F_4,6)\subset \Hilb^{6} \check{\mP}^{2}$
is a smooth prime Fano threefold of genus twelve,
where we use the dual notation for later convenience.
Moreover every general $V_{22}$ is of this form.
\end{thm}

To characterize a general $V_{22}$ he studied
the Hilbert scheme of lines on a general $V_{22}\subset \mP^{13}$
showing that it
is isomorphic to a smooth plane quartic curve $\sH_1\subset \mP^2$.
He thought how to recover $V_{22}$ by $\sH_1$.
For this, one more data was necessary.
Using the incidence relation on $\sH_1\times \sH_1$ defined 
by intersections of lines on $V_{22}$,
he found a non-effective theta characteristic $\theta$ on $\mathcal{H}_1$.
As explained in \cite[\S 6,7]{DK},
there is a beautiful result of G. Scorza  
which asserts that,
associated to the pair $(\sH_1,\theta)$,
there exists another plane quartic curve $\{F_4=0\}$
in the same ambient plane as $\sH_1$.
(By saluting Scorza,
$\{F_4=0\}$ is called the {\em{Scorza quartic}}.)
Then, finally, Mukai proved that
$V_{22}$ is recovered as $\VSP(F_4, 6)$.
Mukai observed that
conics on $V_{22}$ are parameterized by the plane $\sH_2$ and
$\sH_2$ is naturally considered as the plane $\check{\mP}^2$ 
dual to $\mP^2$.
Moreover, he showed,
for one representation of $F_4$ as a power sum
of linear forms $H_1,\dots,H_6$,
the six points $[H_1],\dots, [H_6]\in \check{\mP}^2$ 
correspond to six conics through one point of $V_{22}$.

Even if $F_4$ is taken as a special ternary quartic,
$\VSP(F_4,6)$ may be still a smooth prime Fano threefold of genus twelve.
Mukai \cite[\S 7]{Mu2} shows that,
if $F_4$ is the square of 
a non-degenerate quadratic form,
then $\VSP(F_4,6)$ is so called the
Mukai-Umemura threefold
discovered in \cite{MU}
as a smooth $\mathrm{SO}(3,\mC)$-equivariant compactification
of $\mathrm{SO}(3,\mC)/\mathrm{Icosa}$.
N. Manolache and F.-O. Schreyer
\cite{Abelian1} and F. Melliez and K. Ranestad \cite{Abelian2}
show that, if $F_4$ is the Klein quartic,
then $\VSP(F_4, 6)$ is a smooth compactification of
the moduli space of $(1,7)$-polarized abelian surfaces.

\subsection{Geometry of conics and lines and the main result}
\label{subsection:main}~

Our main result, given in the end of the section $2$, is
a generalization of Mukai's result Theorem \ref{v22};
we describe certain subvarieties of the varieties of power sum of
special quartic forms in any number $v+1$ of variables. The quartics
correspond to the ones of Theorem \ref{v22} if $v=2$.

For this we generalize Mukai's study of the geometries of lines and conics on 
$V_{22}$.
We recall Iskovskih's description of the so-called double projection of a 
$V_{22}$ from a general line as follows:
\begin{equation*}
\xymatrix{
& A' \ar[dl]_{f'}  &\dashrightarrow & 
 A \ar[dr]^{f} & \\
 V_{22}  &  &  & & B, }
\end{equation*}
where
\begin{itemize}
\item
$f'$ is the blow-up along a general line,
\item
$B$ is
the smooth quintic del Pezzo threefold,
namely, a smooth projective threefold
such that $-K_{B}=2H$, where $H$ is the ample generator
of $\Pic B$ and $H^3=5$, and
\item
$f$ is the blow-up along a smooth rational curve of degree five
(with respect to $H$).
\end{itemize}

Generalizing this situation we consider
a general smooth rational curve of degree $d$ on $B$, 
where $d$ is an arbitrary integer greater than or equal to $5$.
In \ref{subsection:CdB},
we establish the existence of such a $C$ and we study 
some of its properties, especially, the relations to
lines and conics on $B$ intersecting it. 
Let $f\colon A\to B$ be the blow up of $B$ along $C$. 
In \ref{subsubsection:lineA} and \ref{subsubsection:coniche},
we define lines and conics on $A$, which are appropriate
generalizations of lines and conics on $V_{22}$.
We say $l$ is a {\em{line}} on $A$ if $l$ is a reduced connected curve
with $-K_A\cdot l=1$, $E_{C}\cdot l=1$ and $p_a(l)=0$, 
where $E_{C}: =f^{-1}(C)$
is the exceptional
divisor of $f\colon A\to B$. 
We say $q$ is a {\em{conic}} on $A$ if $q$ is a reduced connected curve
with $-K_A\cdot q=2$, $E_{C}\cdot q=2$ and $p_{a}(q)=0$. 

We see that 
lines on $A$ are parameterized by a smooth trigonal canonical curve
$\mathcal{H}_1$ of genus $d-2$ (Corollary \ref{primaC}). 
Conics on $A$ turn out to be parameterized by a smooth surface $\sH_2$.
The study of $\sH_2$ is quite delicate.
For this purpose, we consider the intersection of lines and conics
and introduce the divisor
$D_l\subset \sH_2$ parameterizing conics 
which intersect a fixed line $l$. We show that $C$ has $\frac{(d-2)(d-3)}2$ 
bisecant lines and using this we can state the apparently simple
result:
\begin{thm}[see Theorem \ref{thm:H_2}]
\label{thm:main1} The surface $\sH_{2}$ which parameterizes conics on 
$A$ is smooth and 
it is obtained by the blow-up $\eta\colon\sH_{2}\rightarrow S^2 C\simeq \mP^2$ 
at the points $c_{i}$ where $c_{i}$ is the point 
of $S^2 C$ corresponding to the intersection of the bisecants $\beta_i$
and $C$,  $i=1,\ldots , \frac{(d-2)(d-3)}2$.
\end{thm}
Moreover, we show that
if $d\geq 6$, then
$|D_l|$ is very ample and embeds $\sH_2$ in $\check{\mP}^{d-3}$,
and if $d=5$, $|D_l|$ defines a birational morphism $\sH_2\to \check{\mP}^2$.
Here we use the dual notation for later convenience.
If $d\geq 6$, then $\sH_2$ is so called the {\em{White surface}}
(see \cite{White} and \cite{Gimi}).
It is interesting for us 
that the classical White surface naturally appears in this set up.

A deeper understanding of the geometry of conics requires the notion
of intersection of two conics and, more precisely, the divisor
$D_q\subset \sH_2$ parameterizing conics 
which intersect a fixed conic $q$.
It is easy to see that $D_q\sim 2D_l$.

Now assuming $d\geq 6$ we consider $\sH_2\subset \check{\mP}^{d-3}=\mP_*\check{V}$. 
By the double projection of $B$ from a general point $b$,
we see that 
there are $n:=\frac{(d-1)(d-2)}{2}$ conics (counted with multiplicities)
through $b$.
It is crucial that the number $n$ is equal to the dimension
of the quadratic forms on $\check{\mP}^{d-3}$. Nevertheless 
infinitely many conics on $A$
pass through a point on the strict transform of a bi-secant line of $C$.
Hence to have a finiteness result we have to consider the blow-up
$\rho\colon \widetilde{A}\to A$ along the strict transforms 
of bi-secant lines of $C$ on $B$.
Then by a careful analysis on mutually intersecting conics on $A$ we
construct a morphism $\Phi\colon \widetilde{A}\to 
\Hilb^{n}\check{\mP}^{d-3}$ obtained by an attaching process which
associates $n$
conics on $A$ to each point ${\widetilde{a}}$ of ${\widetilde{A}}$; 
see Definition \ref{attached} for the precise definition of attached conics.
To produce the quartic we are looking for, we show that the proper locus 
$\{[q]\in \sH_2 \mid [q]\in D_q\}$ on $\sH_2$ is cut out by a
quartic, 
whose equation is denoted by $\check{F}_4$. 
Moreover we show that $\check{F}_4$ 
is non-degenerate, this means
that the polar map induced by $\check{F}_4$ from 
$S^2 \check{V}$ to $S^2 V$ 
is an isomorphism. 
Then the required quartic $F_4$ is the dual quartic to $\check{F}_4$,
namely, the quartic form in $S^4 \check{V}$
such that its induced polar map from 
$S^2 {V}$ to $S^2 \check{V}$ is the inverse of
that of $\check{F}_4$. 

For the precise statement of our main result, 
we need the following definition:

\begin{defn}
For a subvariety $S$ of $\mP_*\check{V}$,
we set
\[
\VSP(F,n;S):=
\overline{
\{([H_1],\dots, [H_n])\mid [H_i] \in S, H_1^m+\cdots+H_n^m=F\}}
\subset \VSP(F,n)
\]
and we call it the {\em{varieties of power sums of $F$
confined in $S$}}.
\end{defn}
As far as we know, $\VSP(F,n;S)$ is essentially a new object to study.

Our main theoretical result is the following:
\begin{thm}[=Theorem \ref{diretto}]
\label{thm:main2}
There is an injection $\Phi\colon \widetilde{A}\to 
\Hilb^{n}\check{\mP}^{d-3}$ mapping a point $a$ of $\widetilde{A}$ to
the point representing the $n$ conics on $A$ attached to $a$.
Moreover the image is an irreducible component of $\VSP(F_4, n;\sH_2)$.
\end{thm}

In the sequel \ref{subsection:bypro}, we explain a more significant
geometrical meaning of 
the special quartic $F_4$.

Based on Mukai's result we can state the following conjecture: 
$\Phi$ is an embedding and
$\Ima \Phi=\VSP(F_4, n;\sH_2)$.

We remark that,
for $d\leq 8$,
the number $n$ is equal to 
the Waring rank of a general $(d-2)$-nary quartic,
and especially, the cases where $d=5,6,7$ cover 
exceptional cases of Alexander and Hirschowitz.

Even if $d=5$, we have a similar result, 
which is an elaboration on Theorem \ref{v22}.
The explanation is technical: see \ref{subsubsection:v22}.

\subsection{Applications}
\label{subsection:bypro}~

In the section $3$, we give some applications of our study of $A$ 
for a pair of a canonical curve of any genus and 
a non-effective theta characteristic, a {\em{spin curve}} for short.
 
Dolgachev and Kanev \cite[\S 9]{DK} give a modern account of Scorza's
beautiful construction of a certain quartic hypersurface,
so called the {\em{Scorza quartic}}, associated to every spin curve.  
It is expected that the Scorza quartic is useful for the study of a spin curve
but no deeper properties of the Scorza quartic were unknown.
Firstly, its construction is not so explicit. Secondly, 
Scorza's construction itself depends on three assumptions 
on spin curves (see \cite[(9.1) (A1)--(A3)]{DK}) and 
it were unknown whether these conditions are
fulfilled for a general spin curve of genus $>3$.
Thus the existence of the
Scorza quartic was conditional except for the genus $3$ case,
where Scorza himself solved the problem. 
We give contributions for these two subjects.

In \ref{subsection:theta},
using the incidence correspondence on $\sH_1\times \sH_1$ defined 
by intersections of lines on $A$,
we define a non-effective theta characteristic $\theta$ on the
trigonal curve $\mathcal{H}_1$.
This is a generalization of Mukai's result explained as in 
\ref{subsection:Mukai}.

In \ref{subsection:duality}, we observe 
that there is a natural duality between $\sH_1$ and $\sH_2$,
which induces the natural identification
$\mP^{d-3}= \mP^*H^0(\sH_1,K_{\sH_1})$,
where for clarity reasons we denote by
$\mP^{d-3}$ the projective space
dual to the ambient projective space $\check{\mP}^{d-3}$
of $\sH_2$,
and
by $\mP^*H^0(\sH_1,K_{\sH_1})$ the ambient projective space
of the canonical embedding of $\sH_1$.

In \ref{subsection:discr1}, we recall the definition of the discriminant loci 
and we compute it explicitly for $(\sH_1,\theta)$.
In \ref{subsection:Scorza}, we recall the precise definition of
the Scorza quartic for a spin curve.

By virtue of our explicit computation of the discriminant,
we prove in \ref{subsection:conj}
that the pair $(\mathcal{H}_1,\theta)$ satisfies the conditions 
\cite[(9.1) (A1)--(A3)]{DK}, which guarantee the existence of 
the Scorza quartic
for the pair $(\mathcal{H}_1,\theta)$.
Then, by a standard deformation theoretic argument,
we can then verify that the conditions (A1)--(A3) hold
also for a general spin curve,
hence we answer affirmatively to the Dolgachev-Kanev Conjecture:

\begin{thm}[=Theorem \ref{thm:DK}]
\label{thm:main3}
The Scorza quartic exists for a general spin curve.
\end{thm}

Moreover we can find explicitly the Scorza quartic for $(\sH_1,\theta)$. 
In fact, by definition, the Scorza quartic $\{F'_4=0\}$ for $(\sH_1,\theta)$
lives in $\mP^*H^0(\sH_1,K_{\sH_1})$ but, as we remark above, we can consider
$\{F'_4=0\}\subset \mP^{d-3}$.
In \ref{subsection:coincide}, we prove that 
the special quartic $\{F_4=0\}\subset \mP^{d-3}$ in Theorem \ref{thm:main2} coincides with
the Scorza quartic $\{F'_4=0\}$.

We recommend the readers who are interested only in the subsections 3.1--3.5
to skip the subsection 2.5.  

Finally, in \ref{subsection:moduli},
we show that $A$ is recovered from the pair $(\sH_1, \theta)$.
This implies that $(\sH_1,\theta)$'s fill up an open subset of
the moduli of trigonal spin curves. In particular, 
$\sH_1$ is a general trigonal curve for
a general $C$.

\subsection{Final remarks}
\label{subsection:temptativetwo}~

In this paper, we only consider
a general rational curve on $B$
but there are interesting special cases.
In the forthcoming paper, 
applying the method of this paper,
we will study the blow-ups of $B$ along 
special smooth rational curves of degree six
and pairs of canonical curves of genus four and
even theta characteristics. 
   
\begin{ackn}
We are thankful to Professor S. Mukai for 
valuable discussions and constant interest on this paper.
We received various useful comments from 
K. Takeuchi, A. Ohbuchi, S. Kondo,
to whom we are grateful.
The first author worked on this paper
partially while he was working at the
Research Institute for Mathematical Sciences, Kyoto University
until March, 2004, 
and when he was staying at the Johns Hopkins University under the
program of Japan-U.S. Mathematics Institute (JAMI)
in November 2005 and at the Max-Planck-Institut f\"ur Mathematik
from April, 2007 until March, 2008. 
Besides he gave talks on this paper in the seminars and the conferences 
at the Waseda University, the Tohoku University, 
the Johns Hopkins University, the Princeton University, and
Mathematisches Forschungsinstitut Oberwolfach. 
The second author gave a talk on this paper in the seminar at 
the S.I.S.S.A. (Trieste) on May 2007 and at the Workshop on Hd. MMP.
December 2007 (Warwick).
The authors 
worked jointly during the first author's stay at the Universit\`a di Udine
on August 2005, and the Levico Terme conference on Algebraic Geometry in
Higher dimensions on June 2007.
The authors are thankful to all institutes, and 
organizers of the seminars
for the warm hospitality they received. Both authors would like to 
dedicate this paper to Miles Reid with full gratitude. 
\end{ackn}
%
\section{Blowing up the quintic del Pezzo threefold $B$ 
along a $\mP^1$ of degree $d$}\label{part:genus8}~

In this section, we study the geometries of the
blow-up $A$ of the quintic del Pezzo threefold $B$ 
along a smooth rational curve of degree $d$, which is 
nothing but the special threefold we mention in the abstract.
In \ref{subsection:rev} we review
the description of lines and conics on $B$, and
$2$-ray games originating from $B$. Based on this, 
we construct in \ref{subsection:CdB} 
smooth rational curves of degree $d$, 
where $d$ is an arbitrary positive integer, 
having nice intersection properties with respect 
to lines and conics. 
The results in \ref{subsection:CdB} are delicate but their proof is
more or less based on standard parameter count.
In \ref{subsection:H1} and \ref{subsection:H2},
we study the families of curves on $A$ of degree one or two with respect to
the anti-canonical sheaf of ${A}$ (we call them {\em{lines}} and {\em{conics}} on $A$ respectively).
The curve $\sH_1$ parameterizing lines on $A$ and 
the surface $\sH_2$ parameterizing conics on $A$ are two of the main characters
in this paper. See Corollary \ref{primaC} and Theorem \ref{thm:H_2} for a quick
view of their properties.
Finally in \ref{subsection:VSPA}, we prove the main theorem
(Theorem \ref{diretto}). See \ref{subsubsection:v22} for the relationship of our result with
Mukai's one we mentioned in the introduction. 

\subsection{Review on geometries of $B$}
\label{subsection:rev}
\hspace{0.10in}

Let $V$ be a vector space with ${\rm{dim}}_{\mC}V=5$. The
Grassmannian $G(2,V)$ embeds into $\mP^{9}$ and we denote the image 
by $G\subset\mP^{9}$. It is well-known that the quintic del Pezzo $3$-fold,
i.e.,
the Fano $3$-fold $B$ of index $2$ and of degree $5$ 
can be realized as $B=G\cap \mP^{6}$, where $\mP^{6}\subset \mP^9$ 
is transversal to $G$ (see \cite{Fu2}, 
\cite[Thm 4.2 (iii), the proof p.511-p.514]{I1}).

Let $\sH_{{1}}^{B}$ and $\sH_{{2}}^{B}$ be the Hilbert scheme, 
respectively, of lines and of conics on $B$.
We collect basic known facts on lines and conics on $B$ almost without proof.

\subsubsection{Lines on $B$}
\label{subsubsection:FN}
Let $\pi\colon \mP\to \sH_{1}^{B}$ be the universal
family of lines on $B$ and 
$\varphi\colon \mP\to B$ the natural projection. 
By \cite[Lemma 2.3 and Theorem I]{FuNa}, $\mathcal{H}_{1}^{B}$ is isomorphic to $\mP^2$
and $\varphi$ is a finite morphism of degree three. 
In particular the number of lines passing through a point 
is three counted with multiplicities.  
We recall some basic facts about $\pi$ and $\varphi$ which we use in
the sequel.

Before that, we fix some notation.
\begin{nota}
For an irreducible curve $C$ on $B$, 
denote by $M(C)$ the locus $\subset \mP^2$ of lines intersecting $C$,
namely, $M(C):=\pi(\varphi^{-1}(C))$ with reduced structure. 
Since $\varphi$ is flat, $\varphi^{-1}(C)$ is purely one-dimensional.
If $\deg C\geq 2$, 
then $\varphi^{-1}(C)$ does not contain a fiber of 
$\pi$, thus $M(C)$ is a curve.
See Proposition \ref{prop:FN} for the description of $M(C)$
in case $C$ is a line.  
\end{nota}

\begin{defn}
A line $l$ on $B$ is called a {\em{special line}} if
$\sN_{l/B}\simeq 
\sO_{\mP^1}(-1)\oplus \sO_{\mP^1}(1).$
\end{defn}

\begin{prop} 
\label{prop:FN}
It holds$:$ 
\begin{enumerate}[$(1)$]
 \item
 for the branched locus $B_{\varphi}$ of 
 $\varphi\colon \mP\to B$ we have:
\begin{enumerate}[$({1}\text{-}1)$]
\item $B_{\varphi}\in |-K_{B}|$,
\item $\varphi^*B_{\varphi}=R_1+2R_2$,
\item  $R_1\simeq R_2\simeq \mP^1\times \mP^1$, and
\item $\varphi\colon R_1\to
B_{\varphi}$ and 
$\varphi\colon R_2\to
B_{\varphi}$ are injective,
\end{enumerate}
\item $R_{2}$ is contracted to a conic $Q_2$ by 
$\pi\colon \mP\to \sH_{1}^{B}$. Moreover $Q_{2}$ is 
the branched locus of $\pi_{|R_1}\colon R_{1}\rightarrow
\mathcal{H}^{B}_{1}$,
\item $Q_2$ parameterizes special lines.
If $l$ is not a special line on $B$,
then 
$\sN_{l/B}=\sO_{l}\oplus\sO_{l}$,
\item
if $l$ is a special line,
then $M(l)$ is a line, and
$M(l)$ is tangent to $Q_2$ at $[l]$.
If $l$ is not a special line,
then $M(l)$ is the disjoint union of a line and the point $[l]$.
{\em{By abuse of notation}}, we denote by $M(l)$ the one-dimensional part of
$M(l)$ for any line $l$. Vice-versa,
any line in $\mathcal{H}_{1}^{B}$ is of the form $M(l)$ for some line $l$,
\item 
the locus swept by lines intersecting $l$ is a hyperplane section $T_{l}$
  of $B$ whose singular locus is $l$. For every point $b$ of $T_l\setminus l$,
  there exists exactly one line which belongs to $M(l)$ 
  and passes through $b$.
Moreover, if $l$ is not special, then the normalization of $T_l$ is $\mF_1$
and the inverse image of the singular locus is the negative section of $\mF_1$,
or, if $l$ is special, then 
the normalization of $T_l$ is $\mF_3$
and the inverse image of the singular locus is the union of 
the negative section and a fiber,
and
\item
if $l$ is not a special line, then
$\varphi^{-1}(l)$ is the disjoint union of 
the fiber of $\pi$ corresponding to $l$,
and the smooth rational curve dominating $M(l)$.
\end{enumerate}
\end{prop}
\begin{proof} See \cite[\S 2]{FuNa} and \cite[\S 1]{IlievB5}. 
\end{proof}

By the proof of \cite{FuNa} we see that $B$ is stratified
according to the ramification of $\varphi\colon \mP\to B$ as follows:
\[
B=(B\setminus B_{\varphi})\cup 
(B_{\varphi}\setminus C_{\varphi})\cup C_{\varphi},
\] where 
$C_{\varphi}$ is a smooth rational normal sextic and 
if $b\in B\setminus B_{\varphi}$ exactly three
distinct lines pass through it, if $b\in
(B_{\varphi}\setminus 
C_{\varphi})$ exactly two distinct lines pass through it,
one of them is special,
and finally $C_{\varphi}$ is the loci of $b\in B$ 
through which it passes only one line.

\subsubsection{Conics on $B$}
\begin{prop}
\label{prop:conics}
The Hilbert scheme of conics on $B$ is isomorphic to 
$\mP^4=\mP_*\check{V}$. The support of a double line is a special line and
the double lines are parameterized by a
rational normal quartic curve $\Gamma\subset \mP_*\check{V}$ and the
secant variety of $\Gamma$ is a singular cubic hypersurface which
is the closure of the loci parameterizing reducible conics.
\end{prop}
\begin{proof} See \cite[Proposition 1.2.2]{IlievB5}.

\end{proof}

The identification is given by the map 
$sp\colon \sH_{2}^{B}\rightarrow \mP_*\check{V}$ with $[c]
\mapsto\langle Gr(c)\rangle=\mP^{3}_{c}\subset\mP_*V$, where for a
general conic $c\subset B$ we set
\[
Gr(c):=\cup\{r\in \mP_*V\mid [r]\in c\}\simeq \mP^{1}\times\mP^{1}.
\]

\subsubsection{Two-ray games based on $B$}

We are interested in the geometry of the $3$-fold $A$ obtained by
blowing-up $B$ along a curve $C:=C_{d}$ as constructed in 
Proposition \ref{prop:Cd}.
To understand this geometry we need to describe 
some two-ray games originating from $B$.

\begin{defn}
Let $b$ be a point of $B$.
We call the rational map from $B$
defined by the linear system of hyperplane sections
singular at $b$ {\em{the double projection from $b$}}.
\end{defn}

\begin{prop}
\label{prop:proj}
\begin{enumerate}[$(1)$]
\item
Let $b$ be a point of $B$. Then 
the target of the double projection from $b$ is $\mP^2$,
and the double projection from $b$ and 
the projection $B\dashrightarrow \overline{B}_b$ from $b$ 
fit into the following diagram:
\begin{equation}
\label{eq:pt1}
\xymatrix{ & B_b \ar[dl]_{\pi_{1 b}} \ar[dr] & & 
 B'_b \ar[dl] \ar[dr]^{\pi_{2 b}} & \\
 B &  & \overline{B}_b & & \mP^2, }
\end{equation}
where $\pi_{1 b}$ is the blow-up of $B$ at $b$,
$B_b\dashrightarrow B'_b$
is the flop of the strict transforms of lines through $b$, 
and $\pi_{2 b}\colon B'_b\rightarrow \mP^{2}$ 
is a $(\text{unique})$ $\mP^1$-bundle structure.
We denote by
$E_b$ the $\pi_{1b}$-exceptional divisor and
by $E'_b$ the strict transform of $E_b$ on $B'_b$.
Moreover we have the following descriptions$:$
\begin{enumerate}[$({1}$-$1)$]
\item 
\label{eq:pt2}
\[
-K_{B'_b}=H+L,
\]
where 
$H$ is the strict transform 
of a general hyperplane section of $B$,
and $L$ is the pull back of a line on $\mP^2$,
\item

if $b\not \in B_{\varphi}$,
then the strict transforms $l'_i$ of three lines $l_i$ through $b$ on $B_b$
have the normal bundle $\sO_{\mP^1}(-1)\oplus \sO_{\mP^1}(-1)$.
The flop $B_b\dashrightarrow B'_b$ is the Atiyah flop.
In particular, $E'_b\to E_b$ is the blow-up at
the three points $E_b\cap l'_i$.

If $b \in B_{\varphi}\setminus C_{\varphi}$, then 
$E_b\dashrightarrow E'_b$ can be described as follows$:$
let $l$ and $m$ be two lines through $b$, where
$l$ is special, and $m$ is not special.
Let $l'$ and $m'$ be the strict transforms of $l$ and $m$ on $B_b$.
First blow up $E_b$ at two points $t_1:=E_b\cap l'$ and $t_2:=E_b\cap m'$
and then blow up at a point $t_3$ on the exceptional curve $e$ over $t_1$. 
Finally, contract the
strict transform of $e$ to a point. Then we obtain $E'_b$
$($this is a degeneration of the case $(\mathrm{a})$$)$.

See $\textup{\cite{FuNa2}}$ in case of $b \in C_{\varphi}$, and 
\item
a fiber of $\pi_{2 b}$ not contained in $E'_b$
is the strict transform of a conic through $b$, 
or
the strict transform of a line $\not \ni b$ intersecting a line through $b$.

The description of the fibers of $\pi_{2 b}$ contained in $E'_b$ is 
as follows$:$

if $b\not \in B_{\varphi}$, then
$\pi_{2 b|E'_b}\colon E'_b \to \mP^2$ is the blow-down
of the strict transforms of three lines
connecting two of $E_b\cap l'_i$, namely,
$E_b\dashrightarrow \mP^2$ is the Cremona transformation.

Assume that $b \in B_{\varphi}\setminus C_{\varphi}$. 
Then $\pi_{2 b|E'_b}\colon E'_b \to \mP^2$ is the blow-down
of the strict transforms of two lines, one is the line 
connecting $t_1$ and $t_2$,
the other is the line whose strict transform passes through $t_3$. 
$E_b\dashrightarrow \mP^2$ is a degenerate Cremona transformation. 
See $\textup{\cite{FuNa2}}$ in case of $b \in C_{\varphi}$. 

\end{enumerate}
\item
Let $l$ be a line on $B$. Then the 
projection of $B$ from $l$ is decomposed as follow$:$
\begin{equation}
\label{eq:line}
\xymatrix{ & B_l \ar[dl]_{\pi_{1 l}} \ar[dr]^{\pi_{2 l}}\\
 B &  & Q, }
\end{equation}
where $\pi_{1l}$ is the blow-up along $l$ and 
$B \dashrightarrow Q$ is the projection from $l$ and 
$\pi_{2l}$ contracts onto a rational normal curve of degree $3$ the
strict transform of the loci swept by the lines of $B$
touching $l$. Moreover 
\begin{equation}
\label{eq:line2}
-K_{B_l}=H+H_Q,
\end{equation}
where $H$ and $H_Q$ are the pull backs of general hyperplane sections of 
$B$ and $Q$ respectively.
We denote by $E_l$ the $\pi_{1l}$-exceptional divisor.
\item
Let $q$ be a smooth conic on $B$. Then the 
projection of $B$ from $q$ behaves as follow$:$
\begin{equation}
\label{eq:conic}
\xymatrix{ & {B_{q}} \ar[dl]_{\pi_{1 q}} \ar[dr]^{\pi_{2 q}}\\
 B &  & \mP^{3} , }
\end{equation}
\noindent 
where $\pi_{1 q}$ is the blow-up of $B$ along $q$ and 
$\pi_{2 q}\colon {B_{q}}\rightarrow \mP^{3}$ is the
divisorial contraction of the strict transform $T_{q}$ of the loci
swept by
the lines touching $q$. 
Moreover
\begin{equation}
\label{eq:conic2}
-K_{{B_{q}}}=H+H_{\mP},
\end{equation}
where $H$ and $H_{\mP}$ are the pull backs of general hyperplane sections of 
$B$ and $\mP^3$ respectively.
\end{enumerate}
\end{prop}

\begin{proof} These results come from explicit computations and are
more or less known. 
Especially, for (2), refer \cite{Fu2}, and 
for (3) (and (2)), refer \cite{MM1},
No. 22 for (3) (No. 26 for (2)).
See also \cite{MM3}, p.533 (7.7) for a discussion.

(1) is less known. We have only found the paper \cite{FuNa2},
in which they deal with the most difficult case (c).
Here we only sketch 
the construction of the flop in the middle case (b) to intend the reader
to get a feeling of birational maps from $B$.

Let $b$ be a point of $B_{\varphi}\setminus C_{\varphi}$.
We use the notation of the statement of (1-2). 
The flop of $m'$ is the Atiyah flop. We describe the flop of $l'$.
By $\sN_{l/B}\simeq \sO_{\mP^1}(1)\oplus \sO_{\mP^1}(-1)$,
it holds that $\sN_{l'/B_b}\simeq \sO_{\mP^1}\oplus \sO_{\mP^1}(-2)$.
Hence the flop of $l'$ is a special case of Reid's one \cite[Part II]{Pagoda}.
We show that the width is two in Reid's sense.  
Let $T_1$ be the normalization of $T_l$. By Proposition \ref{prop:FN} (5), 
$T_1\simeq \mF_3$ and
the inverse image of the singular locus of $T_l$ is the union of
the negative section $C_0$ and a fiber $r$.
Let $\mu\colon \widetilde{B}_b\to B_b$ be the blow-up along $l'$ and 
$F$ the exceptional divisor. Let $T_2$ be the strict transform of $T_l$ on 
$\widetilde{B}_b$.
Then $T_2$ is the blow-up of $T_1$ at two points $s_1\in C_0$ and $s_2\in r$. 
Denote by $C'_0$ and $r'$ the strict transforms of $C_0$ and $r$.
We prove that $\sN_{r'/\widetilde{B}_b}\simeq \sO_{\mP^1}(-1)^{\oplus 2}$.
Note that $F\cap T_2=C'_0\cup r'$. The curves $C'_0$ and $r'$ are two sections
on $F$. Let $T'_1$ be the image of $T_2$ on $B_b$.
By $\sN_{l'/B_b}\simeq \sO_{\mP^1}\oplus \sO_{\mP^1}(-2)$
and $T_2=\mu^* T'_1-2F$,
it holds $F\simeq \mF_2$, and $T_{2|F}\sim 2G_0+3\gamma$,
where $G_0$ is the negative section of $F$ and 
$\gamma$ is a fiber of $F\to l'$.
Note that $F\cdot C'_0=(F_{|T_2}\cdot C'_0)_{T_2}=-3$ and 
$F\cdot r'=(F_{|T_2}\cdot r')_{T_2}=0$,
and 
$F\cdot G_0=0$ and $F\cdot (G_0+3\gamma)=-3$.
Thus we have $C'_0\sim G_0+3\gamma$ and $r'=G_0$ on $F$. 
Now we see that $-K_{\widetilde{B}_b}\cdot r'=
(\mu^*(-K_{B_b})-F)\cdot r'=0$.
Therefore, by $(r')^2=-1$ on $T_2$, it holds that
$\sN_{r'/\widetilde{B}_b}\simeq \sO_{\mP^1}(-1)^{\oplus 2}$.

It is easy to see that we can flop $r'$.
Let $\widetilde{B}_b\dashrightarrow \widetilde{B}'_b$ be the
flop of $r'$ (now we consider locally around $r'$).
Let $F'$ be the strict transform of $F$ on $\widetilde{B}'_b$.
By \cite{Pagoda}, $F'\simeq F$ and 
there is a blow-down $\widetilde{B}'_b\to \widetilde{B}''_b$ of $F'$ 
such that $\widetilde{B}''_b$ is smooth.
$\widetilde{B}_b\dashrightarrow \widetilde{B}''_b$ is the flop of $l'$.

By this description of the flop, 
we can easily obtain (1-2).
\end{proof}
As a first application of the above operations, we have the following result,
which we often use: 

\begin{cor}
\label{cor:uniqueconic}
Let $b_1$ and $b_2$ be two $($possibly infinitely near$)$ points on
$B$
such that there exists no line on $B$ through them.
Then there exists a unique conic on $B$ through $b_1$ and $b_2$.
\end{cor}

\begin{proof}
We project $B$ from $b_1$ as in (\ref{eq:pt1}).
Then the assertion follows by the description of fibers
of ${\pi_{2 b_{1}}}$ as in Proposition \ref{prop:proj} (1-3). 
\end{proof}

\subsection{Construction of smooth rational curves $C_d$ of 
degree $d$ on $B$}
\label{subsection:CdB}
~

We construct smooth rational curves of degree $d$ on $B$ 
with certain properties.

\begin{prop}
\label{prop:Cd0}
There exists a smooth rational curve $C_d$ of degree $d$ on $B$ such that
\renewcommand{\labelenumi}{\textup{(\alph{enumi})}}
\begin{enumerate}
\item
a general line on $B$ intersecting $C_d$ is uni-secant,
\item 
$C_d$ is obtained as a smoothing of 
the union of a smooth rational curve $C_{d-1}$ of degree $d-1$
on $B$ and a general uni-secant line of it on $B$, 
\item
$\sN_{C_{d}/B}\simeq\sO_{\mP^{1}}(d-1)\oplus\sO_{\mP^{1}}(d-1)$. 
In particular 
$h^1(\sN_{C_{d}/B})=0$ and $h^0(\sN_{C_{d}/B})=2d$, and
\item
if $d=5$, 
then $C_5$ is a normal rational curve and 
is contained in a unique hyperplane section $S$, 
which is smooth.
If $d\geq 6$, then
$C_d$ is not contained in a hyperplane section.
\end{enumerate}
\end{prop}
\begin{proof}
    We argue by induction on $d$.

If $d=1$, we have the assertion since 
$\sN_{C_1/Q}\simeq\sO_{\mP^{1}}\oplus\sO_{\mP^{1}}$
for a general line $C_1$.

Now assume that $C_{d-1}$ is a general smooth rational curve of degree $d-1$
on $B$.
By induction, a general secant line $l$ of $C$ on $Q$ is uni-secant.
Set $Z:=C_{d-1}\cup l$ and $\sN_{Z/Q}: =\sHom_{\sO_{Q}}(\sI_{Z},\sO_{Q})$.
By induction, the normal bundle of $C_{d-1}$ satisfies (c). 

Thus, by $\sN_{l/Q}\simeq\sO_{\mP^{1}}\oplus\sO_{\mP^{1}}$
and \cite[Theorem 4.1 and its proof]{HH},
it holds $h^1(\sN_{Z/Q})=0$, and moreover
$Z:=C_{d-1}\cup l$ is strongly smoothable,
namely, we can find a smoothing $C_d$ 
of $Z$ with the smooth total space.
By the upper semi-continuity theorem,
we have 
$h^1(\sN_{C_d/Q})=0$.

Thus the Hilbert scheme of $Q$ is smooth at $[C_d]$ and is of dimension 
$2d$, which is also the dimension of the component of the Hilbert scheme containing
$[C_d]$. On the other hand, varying $C_{d-1}$ and $l$, we have
a family of reducible curves of dimension $2(d-1)+1=2d-1$.
Thus the smoothing constructed as above is general in 
the component of the Hilbert scheme whose generic point parameterizes
smooth rational curves of degree $d$.

It is easy to see that a general line $m$ intersecting $C_{d-1}$
does not intersect $l$, thus $m$ is a uni-secant line of $C_{d-1}\cup l$.
This implies (a) for $C_d$ by a deformation theoretic argument.

To check the form of the normal bundle, simply assume by induction that
$\sN_{C_{d-1}/Q}=\sO_{\mP^{1}}(d-2)\oplus\sO_{\mP^{1}}(d-2)$. Set 
$\sN_{C_{d}/Q}:=\sO_{\mP^{1}}(a_{d})\oplus\sO_{\mP^{1}}(b_{d})$ ($a_d\geq b_d$)
for the smoothing $C_d$ of $Z$. 
We show that $a_d=b_d=d-1$.

It suffices to prove $h^0(\sN_{Z/Q}(-d))=0$.
In fact, then, by the upper semi-continuous theorem,
we have $h^0(\sN_{C_d/Q}(-d))=0$ and $a_d, b_d\leq d-1$.
Thus, 
by $a_d+b_d=2d-2$, it holds $a_d=b_d=d-1$.
 
The equality $h^0(\sN_{Z/Q}(-d))=0$ easily 
follows from the following three exact sequences, 
where $t:=C_{d-1}\cap l$:

\begin{equation*}
    0\rightarrow \sN_{Z/Q}\rightarrow
    \sN_{Z/Q\vert C_{d-1} }\oplus\sN_{Z/Q\vert l}\rightarrow
    \sN_{Z/Q}\otimes_{\sO_{Q}}\sO_t\rightarrow 0.
 \end{equation*}
 
 \begin{equation*}
    0\rightarrow \sN_{C_{d-1}/Q}\rightarrow
    \sN_{Z/Q\vert C_{d-1}}\rightarrow T^{1}(t)\rightarrow 0.
 \end{equation*}

 \begin{equation*}
    0\rightarrow \sN_{l/Q}\rightarrow
    \sN_{Z/Q\vert l}\rightarrow T^{1}(t)\rightarrow 0.
 \end{equation*}

Finally we prove (d). 
In case $d=5$,  
we have only to notice that a general hyperplane section of $B_{5}$ is
a del Pezzo surface of degree $5$ which contains a smooth
$C_5$. For $d\geq 6$, the assertion follows by induction.
\end{proof}

We denote by $\sH_d^{B}$ the union of components
of the Hilbert scheme of $B$ whose general points
parameterize smooth rational curves of degree $d$
obtained inductively as in Proposition \ref{prop:Cd0}.

The following proposition describe relations between 
$C_d$ and lines and conics on $B$.
\begin{prop}
\label{prop:Cd1}
A general $C_d$ as in Proposition $\ref{prop:Cd0}$
satisfies the following conditions$:$
\begin{enumerate}[$(1)$]
\item
there exist no $k$-secant lines of $C_d$ on $B$ with $k\geq 3$,
\item there exist at most finitely many bi-secant lines of $C_{d}$ on
$B$, and any of them intersects $C_{d}$ simply,
\item bi-secant lines of $C_{d}$ on
$B$ are mutually disjoint,
\item neither a bi-secant line nor
a line through the intersection point between a bi-secant line and $C_d$
is a special line,
\item
there exist at most finitely many points $b$ outside $C_d$
such that all the lines through $b$ intersect $C_d$,
and such points exist outside bi-secant lines of $C_d$,
\item 
there exist no $k$-secant conics of $C_{d}$ with $k\geq 5$,
\item
there exist at most finitely many quadri-secant conics of $C_d$ on
$B$, and no quadri-secant conic is tangent to $C_d$, and
\item $q_{|C_d}$ has no point of multiplicity greater than two
for any multi-secant conic $q$.
\end{enumerate}
\end{prop}

\begin{proof}
We can prove the assertions 
by simple dimension count based upon Proposition \ref{prop:Cd0}.
We assume that $d\geq 4$ since otherwise
we can verify the assertion easily.

$(1)$. Let $\sD$ be the closure of the set
\begin{align*}
\{([C_d], [l])\mid  C_d\cap l\ \text{consists of
$3$ points}
\}  \subset \sH_d^{B}\times \sH_1^{B}.
\end{align*}
Let $\pi_{d}\colon \sD\rightarrow\sH_d^{B}$ and  
$\pi_{1}\colon \sD\rightarrow\sH_1^{B}$ be the natural
morphisms induced by
the projections. The claim follows if we show that
${\rm{dim}}_{\mC}\sD\leq 2d-1$
since $\dim \sH_d^{B}=2d$.

Thus we estimate $\dim_{\mC} \Hom^{2d} (\mP^1, B;p_{i}\mapsto s_{i}, i=1,2,3)$ 
at $[\pi]$, 
where $p_{i}$, $i=1,2,3$ are fixed  points of $\mP^1$, $[\pi]$ 
is a general point and the degree is measured by $-K_{B}$. 
By $d\geq 4$ and Proposition \ref{prop:Cd0} (c), it holds that 
$h^{0}(\mP^{1},\pi^*T_{B}(-p_{1}-p_{2}-p_{3}))=2d-6$ and
$h^{1}(\mP^{1},\pi^*T_{B}(-p_{1}-p_{2}-p_{3}))=0$. 
Then 
\[
\dim_{\mC} \Hom^{2d} (\mP^1, B, p_{i}\mapsto s_{i}, i=1,2,3)_{[\pi]}
=h^0(\pi^*T_{B} (-p_{1}-p_{2}-p_{3}))
=2d-6.
\]
This implies that $\dim_{\mC}\pi_{1}^{-1}([l])\leq 2d-6 +3=2d-3$ since
the three points can be chosen arbitrarily. Then ${\rm{dim}}_{\mC}\sD\leq 2d-1$
since ${\rm{dim}}_{\mC}\sH_1^{B}=2$.
    
$(2)$. Now let 
$\sD$ be the closure of the set
\begin{align*}
\{([C_d], [l])\mid C_d\cap l\ \text{consists of
$2$ points} \}  
\subset \sH_d^{B}\times \sH_1^{B}.
\end{align*}
As before, let $\pi_{d}\colon \sD\rightarrow\sH_d^{B}$ and  
$\pi_{1}\colon \sD\rightarrow\sH_1^{B}$ be the natural
morphisms induced by
the projections.
By $d\geq 4$ and Proposition \ref{prop:Cd0} (c), it holds that
$h^{0}(\mP^{1},\pi^*T_{B}(-p_{1}-p_{2}))=2d-3$ and 
$h^{1}(\mP^{1},\pi^*T_{B}(-p_{1}-p_{2}))=0$. Then
\[
\dim_{\mC} \Hom^{2d} (\mP^1, B, p_{i}\mapsto s_{i}, i=1,2)_{[\pi]}=
h^0(\pi^*T_{B} (-p_{1}-p_{2}))
=2d-3.
\]
Since 
$\dim_{\mC}\Aut (\mP^1, p_{1},p_{2})=1$ it holds that
$\dim_{\mC}\pi_{1}^{-1}([l])\leq 2d-3+2-1=2d-2$. Hence 
${\rm{dim}}_{\mC}\sD=2d$. 
Hence $C_{d}$ has only a
finite number of bisecant lines. 

We now show that the loci where $C_{d}$
has a tangent bisecant is a codimension one loci inside $\sH_d^{B}$.
Let $B_{t}$ be the blow-up of $B$ in a point $t\in C_{d}$ and let $l$ 
be a bi-secant which is tangent to $C_{d}$ at $t$ (if it exists). Let $E$ be 
the exceptional divisor, and $C'$ and $l'$ the strict transforms of $C$ and 
$l$ respectively. By hypothesis there exists a unique point $s\in
E\cap C'\cap l'$. We estimate $\dim_{\mC} \Hom^{d-2} (\mP^1, B_{t}, 
p\mapsto s)_{[\pi]}$, where
$p$ is a fixed  point of $\mP^1$, $[\pi]$ is a general point, and 
the degree is measured by $-K_{B_t}$. 
In this case $h^0(\pi^*T_{B_{t}} (-p))= 2d-2$
hence $\dim_{\mC}\pi_{1}^{-1}([l])\leq 2d-2+ 1 -2=2d-3$. This implies 
the claim.

The cases $(3)$, $(4)$ and $(5)$  are similar. 
Thus we only give few comments for (5).
Set $\sD$ be the closure of the set
\begin{align*}
\{([C_d], [l_1],[l_2],[l_3])\mid 
C_d\cap l_i\not =\emptyset\, (i=1,2,3),\\l_1\cap l_2\cap l_3\not =\emptyset,
l_1\cap l_2\cap l_3\not \in C_d,
l_i \ \text{are distinct}\} \\
\subset \sH_d^{B}\times \sH_1^{B}\times \sH_1^{B}\times \sH_1^{B}.
\end{align*}
For the former half of (5), we have only to prove that $\dim \sD\leq 2d$.
This can be carried out by a similar dimension count as above.
For the latter half of (5), we use the inductive construction 
of $C_d$ besides dimension count. We can omit 
the proof of (6)--(8) since are definitely similar to those of
(1)--(3).



 \end{proof}

\begin{nota}
Denote by $\beta_i$ ($i=1,\cdots, s$) bi-secant lines of $C_d$.
\end{nota}

In the following proposition,
we describe some more relations between $C_d$ and lines on $B$
by using $M(C_d)\subset \sH^B_1$. 
\begin{prop}
\label{prop:Cd}
A general $C_d$ as in Proposition $\ref{prop:Cd0}$
satisfies the following conditions$:$

\begin{enumerate}[$(1)$]
\item
$C_d$ intersects $B_{\varphi}$ simply,

\item 
$M_d:=M(C_d)$ intersects $Q_2$ simply,
\item
$M_d$ is an irreducible curve of degree $d$
with only simple nodes
$($recall that in Proposition $\ref{prop:FN}$,
we abuse the notation by denoting the one-dimensional part of
$\pi(\varphi^{-1}(C_1))$ by
$M(C_1)$$)$,
\item
for a general line $l$ intersecting $C_d$,
$M_d\cup M(l)$ has only simple nodes as its singularities, and
\item
$M_d\cup M(\beta_i)$ has only simple nodes as its singularities.

\end{enumerate}
\end{prop}

\begin{proof}
We show the assertion inductively using
the smoothing construction of $C_d$ from
the union of $C_{d-1}$ and a general uni-secant line $l$
of $C_{d-1}$.

In case of $d=1$, by letting $C_1$ be a general line,
the assertion follows from Proposition \ref{prop:FN}. 
By induction on $d$ assume that we have a smooth $C_{d-1}$ ($d\geq
2$) satisfying (1)--(5).
We verify $C_{d-1}\cup l$ satisfies the following (1)'--(5)',
which are suitable modifications of (1)--(5):\\
(1)' $C_{d-1}\cup l$ intersects $B_{\varphi}$ simply by (1) for $C_{d-1}$
and generality of $l$.\\
(2)'
$M_{d-1}\cup M(l)$ intersects $Q_2$ simply by (2) for $C_{d-1}$
and generality of $l$.\\
(3)'
$M_{d-1} \cup M(l)$ is not irreducible but is of degree $d$
and has only simple nodes by (4) for $C_{d-1}$.\\
(4)'
$M_{d-1}\cup M(l)\cup M(m)$ has only simple nodes as its singularities
for a general line $m$ intersecting $C_{d-1}$.

Indeed, since $m$ is also general, 
$M_{d-1} \cup M(m)$ has only simple nodes by (4) for $C_{d-1}$.
Thus we have only to prove that
$M_{d-1}\cap M(l)\cap M(m)=\emptyset$, namely,
there is no secant line of $C_{d-1}$ intersecting both $l$ and $m$.
Fix a general $l$ and move $m$.
If there are secant lines $r_m$ of $C_{d-1}$ intersecting both $l$ and $m$
for general $m$'s,
then $r_m$ moves whence we have $M(l)\subset M_{d-1}$, a contradiction.
\\
(5)'
For a bi-secant line $\beta$ of $C_{d-1}\cup l$ except the lines
through $C_{d-1}\cap l$,
the curve $M_{d-1}\cup M(l)\cup M(\beta)$ has only simple nodes as its singularities.

Indeed, if $\beta$ is a bi-secant line of $C_{d-1}$, then the assertion follows from
(5) for $C_{d-1}$ by a similar way to the proof of (4)'.
Suppose that $\beta$ is a uni-secant line of $C_{d-1}$ intersecting $l$.
We have only to prove that
there is no secant line of $C$ intersecting both $l$ and $\beta$.
If there is such a line $r$, then
$l$, $\beta$ and $r$ pass through one point.
This does not occur for general $l$ and $\beta$ by
Proposition \ref{prop:Cd1} (5).

Thus, by a deformation theoretic argument, 
we see that $C_d$ satisfies (1)--(5). 
\end{proof}

\begin{nota}
\label{nota:Cb}
Consider the double projection from $b$, see proposition \ref{prop:proj}
[$(1)$].
Throughout the paper,
we denote by $C'_b$, $C''_b$ and $C_b$ 
the strict transforms of $C:=C_d$ on $B_b$, $B'_b$ and $\mP^2$
respectively. 
\end{nota}

The following result is one of the key properties of the component
$\sH^{B}_{d}$. Its importance and difficulty
lies in the actual fact that it holds for
every $b\in B$.

\begin{prop}
\label{prop:Cd2}
Let $C_d$ be a general smooth rational curve of 
degree $d$ on $B$ constructed as in Proposition $\ref{prop:Cd0}$.
Then,
for any point $b\in B$,
the restriction of $\pi_b$ to $C_d$ is birational if $d \geq 5$.
\end{prop}
\begin{proof}
We prove the assertion by induction
based on the construction of $C_d$ from 
$C_{d-1}\cup l$, where
$l$ is a general uni-secant line of $C_{d-1}$ on $B$.

Assume that $d=5$ and
$\pi_{b|C_5}$ is not birational for a $b$.
Then $C_b$ is a line or conic in $\mP^2$.
Let $S$ be the pull-back of $C_b$ by $\pi_{2b}$.
If $C_b$ is a line, then $C_5$ is contained in a
singular hyperplane section, which is the strict transform of
$S$ on $B$ 
(recall that $B\dashrightarrow \mP^2$ is the double projection from $b$). 
This contradicts Proposition \ref{prop:Cd0} (d).
Assume that $C_b$ is a conic.

The only possibility is that
$L\cdot C''_b=4$ and $C''_b\to C_b$ is a double cover
since $\deg C_b \cdot \deg (C''_b\to C_b)\leq 5$. 
By Proposition \ref{prop:proj} (1-1),
it holds $H\cdot C''_b=6$.
Then by $L=H-2E'_b$, we have $E'_b\cdot C''_b=1$.
Note that ${E'_b}^2 S=2$.
On $S$, we can write 
$E'_{b|S}\sim C_0+pl$ and $C''_b\sim 2C_0+ql$ $(p,q\geq 0)$,
where $C_0$ is the negative section of $S$ and $l$ is a fiber of $S\to C''_b$.
Set $e:=-C^2_0$.
By $E'_b\cdot C''_b=1$ and ${E'_b}^2 S=2$,
we have $q+2p-2e=1$ and $2p-e=2$.
Thus $e=2p-2$ and $q=2p-3$.
Since $C''_b$ is irreducible,
$q\geq 2e$, whence $2p-3\geq 2(2p-2)$, i.e., $p=0$ and $q=-3$, a contradiction.

Assume that $d\geq 6$.
Let $\sC\to \Delta$ be the one-parameter smoothing of $C_{d-1}\cup l$
such that $\sC$ is smooth.
We consider the trivial family of the double projections
$B\times \Delta\dashrightarrow \mP^2\times \Delta$
from $b\times \Delta$.
Denote by $\sC'_b, \sC''_b$ and $\sC_b$ the strict transforms of
$\sC$ on $B'_b\times \Delta$, $B''_b\times \Delta$
and $\mP^2\times \Delta$ respectively.
We also denote by $C'_{d-1,b}$, $C''_{d-1,b}$, and $C_{d-1,b}$
the strict transforms of
$C_{d-1}$ on $B'_b$, $B''_b$ and $\mP^2$ respectively.
It suffices to
prove that there exists at least one point on 
$C_{d-1}$ where $\sC\dashrightarrow \sC_b$ is birational.

Indeed, set 
\[
\sN:=
\{(b, t)\in B\times \Delta\mid \sC\dashrightarrow \sC_b 
\ \text{is not birational at any point of } \ \sC_t \}
\]
and let $\Delta'\subset \Delta$ be the image of $\sN$ 
by the projection to $\Delta$. 
$\sN$ is a closed subset, and so is $\Delta'$ since $B\times \Delta\to \Delta$
is proper. Thus $\Delta'$ consists of finitely many points since the origin is not contained in $\Delta'$.
For a point $t\in \Delta$ sufficiently near the origin,
$\sC_t\dashrightarrow \sC_{t,b}$ is birational for any $b$.

By induction, we may assume that $C_{d-1}\dashrightarrow C_{d-1,b}$
is birational.
Note that $C_{d-1,b}$ is not a line since
otherwise $C_{d-1}$ is contained in 
a singular hyperplane section as we see above in the case of $C_5$,
a contradiction. 
As for $l$, if $b\not \in l$, then
the image of $l$ is a line or a point on $\mP^2$.
If $b\in l$, then the strict transform of $l$ on $B_b$ is a flopping
curve. Thus $\sC_b$ contains the line corresponding to $l$.
We investigate the other possible irreducible components of 
the central fiber $\sC_{b,0}$ of $\sC_b\to \Delta$.  
If $b\not \in C_{d-1}\cup l$, then
the only possibility is that 
$\sC_{b,0}$ contains the image of a flopped curve,
which is a line on $\mP^2$.
Thus 
$\sC\dashrightarrow \sC_b$ is birational at a point of $C_{d-1}$.
If $b\in C_{d-1}\cup l$, then
$\sC_{b,0}$ contains the image $m_b$ of 
the strict transform $m''_b$ of a line $m'_b$
in $E_b$ through $E_b\cap (C'_{d-1,b}\cup l'_b)$, 
where $l'_b$ is the strict transform of $l$ on $B_b$. 
The line $m'_b$ is nothing but the exceptional curve for $\sC'_b\to \sC$
(recall that $\sC$ is a smooth surface).
Moreover, if $b\in l$, then
by the description of $E_b\dashrightarrow \mP^2$,
$m_b$ is a line since $l_b$ is a flopping curve.
Thus 
$\sC\dashrightarrow \sC_b$ is birational at a point of $C_{d-1}$.
Suppose that $b\in C_{d-1}\setminus l$.
If $m'_b$ intersects a flopping curve,
$m_b$ is a line or a point, thus 
$\sC\dashrightarrow \sC_b$ is birational at a point of $C_{d-1}$.
In the other case, $m_b$ is a conic.
If $b\not \in \cup_i \beta_i$, then
$\deg C_{d-1,b}=d-3$ by Proposition \ref{prop:proj} (1-1).
By $d\geq 6$, $C_{d-1,b}$ is not a conic, thus 
$\sC\dashrightarrow \sC_b$ is birational at a point of $C_{d-1}$.
Assume $b\in \beta_i$.
Then 
$\deg C_{d-1,b}=d-4$. We have only to show that
if $d=6$, then $C_{d-1,b}\neq m_b$.
By Proposition \ref{prop:Cd1} (4),
the flop is of type (a) in Proposition \ref{prop:proj} (1-2).
The line $m'_b$ intersects three lines which are the strict transforms
of three fibers of $\pi_b$ contained in $E'_b$.
On the other hand, by $E'_b\cdot C''_{d-1,b}=2$,
the curve $C''_{d-1,b}$ intersects
at most two fibers of $\pi$ contained in $E'_b$.
Thus $C_{d-1,b}\neq m_b$, and 
$\sC\dashrightarrow \sC_b$ is birational at a point of $C_{d-1}$.
\end{proof}

We restate the proposition in terms of the relation 
between $C_d$ and multi-secant conics of $C_d$ on $B$ as follows:
\begin{cor}
\label{cor:finite}
Let $b$ be a point of $B$ not in any bi-secant line of $C_d$ on $B$.
If $d\geq 5$, then
there exist finitely many $k$-secant conics of $C_d$ on $B$ 
through $b$ with $k\geq 2$ if $b\not \in C_d$
$(\text{resp. with $k\geq 3$ if $b \in C_d$})$.
\end{cor}

\begin{proof}
For a point $b\in B$ outside bi-secant lines of $C_d$ on $B$,
there exist a finite number of 
singular multi-secant conics of $C_d$ through $b$
since the number of lines through $b$ is finite, and
the number of lines intersecting both a line through $b$ and $C_d$
is also finite by Proposition \ref{prop:Cd} (3).
Therefore
we have only to consider smooth multi-secant conics $q$ of $C_d$ through $b$.
By Proposition \ref{prop:proj} (1-3),
the strict transform $q'$ of 
such a conic $q$ on $B'_b$ is a fiber of $\pi_{2b}$.
If $b\not \in C_d$,
then $q'$ intersects $C'_b$ twice or more counted with multiplicities, 
thus by Proposition \ref{prop:Cd2}, the finiteness of such a $q$ follows.
We can prove the assertion in case of $b\in C_d$ similarly, 
thus we omit the proof. 
\end{proof}

\begin{rem}
We refine this statement in Lemmas \ref{fuoripiatto} and
\ref{finitezzafinale}.
\end{rem}

\subsection{Curve $\sH_1$ parameterizing marked lines}
\label{subsection:H1}~

We fix a general $C:=C_{d}$ as in \ref{subsection:CdB}.
Let $f\colon A\to B$ be the blow-up along $C$.
We start the study of the geometry of $A$. The first step consists of 
finding the curves, if any, which replace the lines of ordinary
geometry.

\subsubsection{Construction of $\sH_1$ and marked lines}

Set $\sH_1:=\varphi^{-1} C\subset \mP$ and $M:=M_d$.
We begin with a few corollaries of Proposition \ref{prop:Cd}:
   
\begin{cor}\label{primaC}
If $d\geq 2$, then 
 $\sH_1$ is a smooth curve of genus $d-2$
with the triple cover $\sH_1\to C$.
In particular, if $d\geq 3$, then 
$\sH_1$ is a smooth non-hyperelliptic trigonal curve of genus $d-2$.
\end{cor}

\begin{proof}  
By Propositions \ref{prop:FN} and \ref{prop:Cd} (1),
it holds that $\sH_1$ is smooth and
the ramification for $\sH_1\to C$ is simple 
by Proposition \ref{prop:Cd} (1).
Since $B_{\varphi}\in |-K_{B}|$ and $d={\rm{deg\,{C}}}$,
we can compute $g(\sH_1)$ by the Hurwitz formula:
\[
  \text{$2g(\sH_1)-2=3\times (-2)+d\times 2$, equivalently,  
  $g(\sH_1)=d-2$}.      
\]
\end{proof}

\begin{cor}
    \label{orabisecanti}
The number $s$ of nodes of $M$ is $\frac{(d-2)(d-3)}{2}$, whence
$C$ has $\frac{(d-2)(d-3)}{2}$ bi-secant lines on $B$.
\end{cor} 

\begin{proof}    
By the inductive construction of $C$ 
we see that $\pi_{|\sH_1}\colon\sH_1\rightarrow M$ is birational. By \ref{prop:Cd}
$(3)$ $p_a(M)=\frac{(d-1)(d-2)}{2}$. Then by \ref{prop:Cd}
$(3)$ we know the number of nodes of $M$ since $g(\sH_1)=d-2$. 
The latter half follows since a bi-secant line of 
$C$ corresponds to a node of $M$. 
\end{proof}


Now we select some lines on $B$ which we use in the sequel. Note that
\[
\sH_1=\{([l],t)\mid [l]\in M, t \in C\cap l\} \subset M\times C,
\]
and the elements of $\sH_1$ deserve a name:

\begin{defn}
The pair of a secant line $l$ of $C$ on $B$ and 
a point $t\in C\cap  l$ is called a {\it{marked line}}.

Let $(l,t)$ be a marked line.
If $C\cap l$ is one point, then 
$\{t\}=C\cap l$ is uniquely determined. 
For a bi-secant line $\beta_i$ of $C$, 
there are two choices of $t$.  
Thus $\sH_1$ parameterizes marked lines.
\end{defn}



\subsubsection{Lines on the blow-up $A$ of $B$ along $C_d$}
~
\label{subsubsection:lineA}

We prove that each marked line corresponds to a curve 
of anticanonical degree $1$ on the blow-up
$A$ of $B$ along $C$. This gives us a suitable notion
of line on $A$.

\begin{nota}
~
\begin{enumerate}[$(1)$] 
\item \parbox[t]{0.80\textwidth}{
Let $f\colon A\to B$ be the blowing up along $C$ and 
$E_{C}$ the $f$-exceptional divisor,} 

\item
 $\{p_{i1}, p_{i2}\}=C\cap \beta_{i}\subset B$,
\item
$\zeta_{ij}=f^{-1}(p_{ij})\subset E_{C}\subset A$, and
\item
$\beta'_{i}\cap \zeta_{ij}= p'_{ij}\in E_{C}\subset A$,
\end{enumerate}
where $i=1,\ldots, s$ and $j=1,2$.
\end{nota}

\begin{defn}
    \label{linea}
    We say that a connected curve $l\subset {A}$ 
    is a {\em{line}} on ${A}$ if
\begin{enumerate}[(i)]
\item
$-K_{{A}} \cdot l=1$, and
\item
${{E_{C}}}\cdot l=1$.
\end{enumerate}
\end{defn}

We point out that since $-K_A=f^*(-K_{B})-E_{C}$ and $E_{C}\cdot
l=1$ then
$f(l)$ is a line on $B$ {\it{intersecting}} $C$.
More precisely:
\begin{prop}\label{lineeA}
A line $l$ on $A$ is one of the following curves on $A:$
\renewcommand{\labelenumi}{\textup{(\roman{enumi})}}
\begin{enumerate}
\item
the strict transform of a uni-secant line of $C$ on $B$,
or
\item
the union $l_{ij}=\beta'_{i}\cup \zeta_{ij}$, where $i=1,\ldots, s$
and $j=1,2$.
\end{enumerate}
In particular $l$ is reduced and $p_a(l)=0$.
\end{prop}

\begin{nota}
For a line $l$ on $A$, we usually denote by $\overline{l}$ its image on
$B$.
\end{nota}

\begin{cor}\label{parameter} 
    The curve $\sH_1\subset\mP$ is the Hilbert scheme of
     the lines of $A$.
\end{cor}

\begin{proof} 
Let $\sH'_1$ be the Hilbert scheme of lines on $A$,
which is a locally closed subset of the Hilbert scheme
of $A$.
By the obstruction calculation of the normal bundles of the
components of
lines on $A$, it is easy to see that $\sH'_1$ is a smooth curve.
Denote by $\sU_1\to \sH'_1$ the universal family of the lines on 
$A$ and let $\overline \sU_1$ be the image of $\sU_1$
on $B\times \sH'_1$ (with induced reduced structure).

\begin{cla} 
\label{cla:basic}
$\overline {\sU}\to \sH'_{1}$ is a $\mP^1$-bundle.
\end{cla}
\begin{proof}[Proof of the claim]
Let $\sL$ be the pull-back of 
the ample generator of $\Pic B$
by 
\[
{\sU_{1}}\hookrightarrow 
A \times \sH'_{1}\to 
B\times \sH'_{1}\to B.
\]
Since $\varrho\colon \sU_{1}\to \sH'_{1}$ is flat
and $h^0(l, \sL_{|l})=2$ for a line $l$ on ${B}$,
$\sE:=\varrho_*\sL$ is a locally free sheaf of rank two.
$\mP(\sE)$ is nothing but the $\mP^1$-bundle
contained in $B\times \sH'_{1}$ whose fiber
is the image of a line on ${A}$.
This implies that $\mP(\sE)=\overline \sU$ as schemes and 
$\overline \sU$ is a $\mP^1$-bundle.
\end{proof}
By the claim same we have a natural morphism
$\sH'_1\to \mP^2$, 
whose image is $M$.
By Proposition \ref{lineeA} $\sH'_1 \to M$ is birational and surjective. 
Since $\sH'_1$ and $\sH_1$ are smooth, they are both normalizations of
$M$, then $\sH'_1\simeq \sH_1$.
\end{proof}

\begin{rem}
For a bi-secant line $\beta_i$,
we have two choices of marking, $p_{i1}$ or $p_{i2}$.
We describe which line on $A$ corresponds to $(\beta_i,p_{ij})$. 
Denote by $\sU_1\to \sH_1$ the universal family of the lines on 
$A$ and consider the following diagram:
\begin{equation*}
\xymatrix{{\sU_1}\ar[d] \subset & 
A\times \sH_1 \ar[d]\\
{\overline{\sU}_1}  \subset &
B \times \sH_1.}
\end{equation*}
Then
$\sU_1\to \overline{\sU}_1$ is the blow-up
along $(C\times \sH_1)\cap \overline{\sU}_1$,
which is the union of a section of $\overline{\sU}_1\to \sH_1$
consisting markings
and finite set of points
$(p_{i,3-j}, [\beta_i, p_{ij}])$.
Thus 
the marked line $(\beta_i, p_{ij})$ corresponds to
the line $l_{i,3-j}$.
\end{rem}

\subsection{Surface $\sH_2$ parameterizing marked conics}
\label{subsection:H2}
~

Now we define a notion of {\it{conic}} on $A$. We proceed as in the
case of lines, first defining the notion of {\it{marked conic}}.

\subsubsection{Construction of $\sH_2$ and marked conics}
\label{subsubsection:H_2}

\begin{defn}\label{conichemarked}
The pair of a $k$-secant conic $q$ on $B$ with $k\geq 2$
and a zero-dimensional subscheme $\eta \subset C$
of length two contained in $q_{|C}$ 
is called a {\it{marked conic}}.
\end{defn}

From now on, we assume that $d\geq 3$.

Marked conics are parameterized by
\[
\sH'_2:=\{([q],[\eta])\mid [q]\in \overline{\sH}'_2, \eta\subset C\cap
q\}
\subset \overline{\sH}'_2\times S^2 C
\]
with reduced structure,
where
$\overline{\sH}'_2\subset \mP^4$ is the locus of
multi-secant conics of $C$ on $B$.

By Corollary \ref{cor:uniqueconic} and $d\not =1$,
the natural projection of $\sH'_2\to S^2 C$ 
is one to one outside $[\beta_{i|C}]$ and the diagonal of $S^2 C$, thus 
by the Zariski main theorem, it is an isomorphism
outside $[\beta_{i|C}]$ and the diagonal of $S^2 C$.

We denote by 
$e'_i$ the fiber of  
$\sH'_2\to S^2 C$ over a $[\beta_{i|C}]$.
Since $B$ is the intersection of quadrics,
any conic cannot intersect a line twice properly.
Thus any conic $\supset \beta_{i|C}$ contains $\beta_i$.
This implies that $e'_i\simeq \mP^1$, and $e'_i$
parameterizes marked conics of the form
\[
\{([\beta_i\cup \alpha], [\beta_{i|C}])\mid 
\text{$\alpha$ is a line such that $\alpha\cap \beta_i\not
=\emptyset$}\}.
\]

Over the diagonal of $S^2 C$,
$\sH'_2\to S^2 C$ is finite since 
for $t\in C$,
there exist a finite number of reducible conics with $t$
as a singular point or conics tangent to $C$ at $t$.

Hence $\sH'_2$ is the union of 
the unique two-dimensional component, which dominates $S^2 C$,
and possibly lower dimensional components mapped into
the diagonal of $S^2 C$ or $e'_i$.
Note that $\sH'_2\to \overline{\sH}'_2$ is finite since
choices of markings of a multi-secant conic of $C$ is finitely many
by $d\geq 3$. 
\begin{cla}
$e'_i$ is contained 
in the unique two-dimensional component of $\sH'_2$.
\end{cla}

\begin{proof}
We have only to prove that $\overline{\sH}'_2$ is two-dimensional
near the generic point of the image of $e'_i$ since $\sH'_2\to \overline{\sH}'_2$ is one to one near the generic point of the image of $e'_i$.
Let $\sV_2\to \sH_2^{B}\simeq \mP^4$ be the universal family of
conics on $B$ and $\overline{\sH}''_2$
the inverse image of $C\times C$ 
by $\sV_2\times _{\mP^4} \sV_2 \to B\times B$.
Since the morphism $\sV_2\times _{\mP^4} \sV_2 \to \sV_2\to \mP^4$ is flat,
$\sV_2\times _{\mP^4} \sV_2$ is purely six-dimensional.
Thus any component of $\overline{\sH}''_2$ 
has dimension greater than or equal to two.
Though the inverse image of
the diagonal of $C\times C$ is three-dimensional,
any other component of $\overline{\sH}''_2$ is at most two-dimensional 
by a similar investigation to $\sH'_2$.
Thus $\overline{\sH}'_2$ is two-dimensional
near the generic point of the image of $e'_i$ since
$\overline{\sH}'_2$ is the image of the two-dimensional part of 
$\overline{\sH}''_2$ by $\sV_2\times _{\mP^4} \sV_2\to \mP^4$
near the generic point of the image of $e'_i$.
\end{proof}

\begin{nota}
Let $\sH_2$ be the normalization of the unique two-dimensional component of
$\sH'_2$ and $\overline{\sH}_2\subset \overline{\sH}'_2$ 
the image of $\sH_2$. 
Denote by $\eta$ the natural morphism $\sH_2\to S^2 C$.
Set 
\[
c_i:=[\beta_{i|C}] \in S^2 C\simeq \mP^2,
\]
and
\[ 
e_i:={\eta}^{-1}(c_i),
\] 
where $i=1,\ldots, s$.
\end{nota}

By the above consideration,
$\eta\colon \sH_2\to S^2 C$ is isomorphic outside $[\beta_{i|C}]$
and $\sH_2\to \overline{\sH}_2$ is the normalization.
Thus we see that $\sH_2$ parameterizes marked conics
outside the inverse image of $c_i$. We need to understand the inverse
image by $\eta$ of the diagonal.

\begin{cla}
\label{cla:reduced}
Assume that $([q],[2b])\in \sH_2$ for $b\in C$ and a conic $q$. Then
\begin{enumerate}[$(1)$]
\item
$q$ is reduced,
\item
if $q$ is smooth at $b$, then $q$ is tangent to $C$ at $b$, and
\item
if $q$ is singular at $b$, then
the strict transform of $q$ is connected on $A$.
Moreover, $b\not \in \beta_i$ nor $B_{\varphi}$.
\end{enumerate}
\end{cla}

\begin{proof}
By Proposition \ref{prop:proj} (1-3) and a degeneration argument, $q$ corresponds to the fiber of $\pi_{2b}$ through 
the point $t'$ in $C''_b\cap E'_b$ coming from $t:=C'_b\cap E_b$.

(1) Assume by contradiction that $q$ is non-reduced.
By Proposition \ref{prop:conics}, $q$ is a multiple of 
a special line $l$. By Proposition \ref{prop:Cd1} (4),
$l$ is a uni-secant line of $C$.
Let $m$ be the other line through $b$ (by generality of $C$, we have $l\not =m$). Let $l'$ and $m'$ be the strict transforms of $l$ and $m$ on $B_b$ respectively. By Proposition \ref{prop:proj} (1-3),
the fiber of $\pi_{2b}$ through $t'$ is the strict transform
of the line in $E_b$ joining $l'\cap E_b$ and $m'\cap E_b$.
Hence by the assumption,
$l'\cap E_b$, $m'\cap E_b$ and $C'_b\cap E_b$ are collinear.
By dimension count similar to 
the proof of Proposition \ref{prop:Cd1}, 
we can prove that a general $C$ does not satisfy this condition.

(2) This follows from the previous discussion.

(3) Set $q=l_1\cup l_2$, where
$l_1$ and $l_2$ are the irreducible components of $q$,
and let $l'_i$ be the strict transform of $l_i$ on $B_b$.
By (1), it holds $l_1\neq l_2$.
Then the fiber of $\pi_{2b}$ corresponding to $q$ 
is the strict transform of the line on $E_b$ through $E_b\cap l'_1$
and $E_b\cap l'_2$.
Note that $A$ is obtained from $B_b$ by blow- up $B_b$ along $C'_b$
and then contracting the strict transform of $E_b$.
Thus the former half of the assertion follows.
The latter half follows again by simple dimension count.
\end{proof}


\subsubsection{Conics on $A$}
\label{subsubsection:coniche}

\begin{defn}
    \label{conica}
    We say that a curve $q\subset {A}$ 
    is a {\em{conic}} on ${A}$ if
    \\[1ex]

\begin{math}
\begin{aligned}
    \text{(i) } & \text{$q$ is connected and reduced, \quad}  & 
    \text{(ii) } & -K_{{A}} \cdot q = 2,   \\
    \text{(iii) } & {{E_{C}}}\cdot q=2, \ \text{and}
    &
    \text{(iv) } & p_a(q) =0.
\end{aligned}
\end{math}
\end{defn}

Using this definition, we can classify conics on $A$
similarly to Proposition \ref{lineeA}:

\begin{prop}\label{classificazioneconicheinA}
Let $q$ be a conic on ${A}$.
Then $\overline q:=f(q)\subset B$ is a $k$-secant conic of 
$C$ with $k\geq 2$.
Moreover one of the following holds$:$
\renewcommand{\labelenumi}{\textup{(\alph{enumi})}}
\begin{enumerate}
\item
$\overline q$ is smooth at $\overline q\cap C$. 
$q$ is the union of the strict transform $q'$ 
of $\overline q$
and $k-2$ distinct fibers $\zeta_1,\dots, \zeta_{k-2}$ of $E_C$ 
such that $\zeta_i\cap q' \not =\emptyset$,
\item
$\overline{q}$ is the union of two uni-secant lines $\overline l$ and
$\overline m$ such that 
$C\cap \overline l\cap \overline m \not =\emptyset$.
$q$ is the union of the strict transforms $l$ and $m$ 
of $\overline l$ and
$\overline m$ respectively $($we assume that $l\cap m\not =\emptyset$$)$, or 
\item
$\overline{q}$ is the union of $\beta_i$ and a line $\overline{r}$ through a
$p_{ij}$. $q$ is the union of the fiber $\zeta_{ij}$ over $p_{ij}$
and the strict transforms $\beta'_i$ and $r'$ 
of $\beta_i$ and $\overline{r}$ respectively.
\end{enumerate}  
\end{prop}



\begin{nota}
We usually denote by $\overline{q}\subset B$ the image of a conic $q$ on $A$.
\end{nota}
 
Let $\sH^A_2$ be the normalization of the two-dimensional part of the Hilbert scheme of conics on $A$,
which is a locally closed subset of the Hilbert scheme
of $A$.
Let
$\mu \colon {\sU_2}\to \sH^A_2$ be the pull-back of
the universal family of conics on ${A}$.

\begin{lem}\label{useful?}
    Let
$\overline{\sU}_2$ be the image of ${\sU_2}$
on $B\times \sH^A_2$ $($with induced reduced structure$)$ then
$\overline \sU_2\to \sH^A_2$ is a conic bundle. 
\end{lem}

\begin{proof}
The proof is similar to that of Claim \ref{cla:basic}.

Let $\sL$ be the pull-back of 
the ample generator of $\Pic B$
by 
\[
{\sU_{2}}\hookrightarrow 
A \times \sH^A_2\to 
B\times \sH^A_2\to B.
\]

Since $\mu\colon {\sU_2}\to \sH^A_2$ is flat
and $h^0(q, \sL_{|q})=3$ for a conic
$q$ on ${A}$ (recall that $q$ is reduced), then
$\sE:=\mu_*\sL$ is a locally free sheaf of rank $3$.
Letting
$\mP^{6}=\langle B \rangle$, $\mP(\sE)$ is the $\mP^2$-bundle
contained in $\mP^{6}\times\sH^A_{2}$ whose fiber
is the plane spanned by the image of a conic on ${A}$.
Let $\mathcal{Q}:=(B\times \sH^A_2) \cap \mP(\sE)$,
where the intersection is taken in $\mP^{6}\times \sH^A_2$.  
A scheme theoretic fiber of $\sQ \to \sH^A_2$ 
is the image of a conic of ${A}$ since
$B$ is the intersection of quadrics.
Then $\sQ=\overline{\sU}_2$ as schemes and 
$\overline{\sU}_2$ is a conic bundle.
\end{proof}

\begin{prop}\label{equality}
The two surfaces $\sH_2^{A}$ and $\sH_2$ are isomorphic.
\end{prop}

\begin{proof} 
By Lemma \ref{useful?},
there exists a natural morphism $\overline{\nu}\colon \sH^A_2\to \overline{\sH}'_2$.
By Proposition \ref{classificazioneconicheinA},
$\overline{\nu}$ is finite and birational, 
hence $\overline{\nu}$ lifts to the morphism $\nu\colon \sH^A_2 \to \sH_2$
since $\sH_2\to \overline{\sH}_2$ is the normalization.
By the Zariski main theorem, $\nu$ is an inclusion.
By Claim \ref{cla:reduced} (1) and (3), 
and Proposition \ref{classificazioneconicheinA}, 
$\nu$ is also surjective.
\end{proof}
 
 By Proposition \ref{equality} we can pass freely from conics on $A$,
 that is elements of $\sH_2^{A}$ to marked conics and vice-versa
 according to the kind of argument we will need. In particular we can
 speak of the universal family $\mu\colon \sU_{2}\rightarrow \sH_{2}$ of marked 
 conics meaning $\sU_{2}:=\sU^{A}_{2}$ and $\sH_2^{A}$ identified to
 $\sH_{2}$ via $\nu$. 
 
\begin{cor}
The Hilbert scheme of conics on $A$ is an irreducible surface
$($and $\sH_2$ is the normalization$)$.
The normalization is injective, 
namely, $\sH_2$ parameterizes conics on $A$ in one to one way.  
\end{cor}

\begin{proof}
By Proposition \ref{classificazioneconicheinA},
the image of $\sH_2$ in the Hilbert scheme parameterizes
all the conics, thus the first part follows.

For the second part,
we have already seen that $\sH_2$
parameterizes conics on $A$ in one to one way outside $\cup_i e_i$.
Let $\alpha$ be a general line intersecting $\beta_i$,
and $\alpha'$ the strict transform of $\alpha$ on $A$.
By easy obstruction calculation, we see that
the Hilbert scheme of conics on $A$ is smooth at $[\beta'_i\cup \alpha']$.
Thus general points of $e_i$ also parameterizes conics on $A$.
Then, however, since $e'_i\simeq \mP^1$, where
$e'_i$ is the inverse image of $[\beta_{i|C}]$ by $\sH'_2\to S^2 C$,
it holds that $e_i\simeq e'_i \simeq \mP^1$
($\sH_2\to S^2 C$ has only connected fibers).
This implies the assertion.
\end{proof}
\medskip


\subsubsection{Description of  $\sH_2$}
~
 
 We want to investigate further
the morphism $\eta\colon \sH_2\to S^2 C\simeq \mP^2$. 

\begin{nota}
For a point $b\in C$, set
\[
L_b:=\overline{\{[q]\in \sH_2 \mid \exists, b'\not =b,
f(q)\cap C=\{b, b' \} \}}.
\]
By Corollary \ref{cor:uniqueconic}, 
$\eta(L_b)$ is a line in $S^2 C \simeq \mP^2$.
\end{nota}

To understand better $\eta\colon \sH_2\to\mP^2$ 
we need to find special loci inside $\sH_{2}$. A natural step
is to study the locus of conics which intersect a fixed line.

Let ${\sU}'_1 \subset {\sU_2}\times \sH_1$ be the
pull-back of ${\sU_1}$ via the following diagram:
\begin{equation}
\label{eq:familyB1}
\xymatrix{\sU'_1\subset {\sU_2}\times \sH_1 \ar[d] 
\ar[r] & A \times \sH_1 \supset \sU_1 \ar[d]\\
\widehat{\sD}_1 \subset \sH_2\times \sH_1 \ar[r] & \sH_1,} 
\end{equation}
where $\widehat{\sD}_1$ is the image of $\sU'_1$ on 
$\sH_2\times \sH_1$.

By definition
\[
\widehat{\sD}_1=\{([q], [l])\mid q \cap l \not =\emptyset\}
\subset \sH_2\times \sH_1.
\]

First we need to know which component of $\widehat{\sD}_1$
is divisorial or dominates $\sH_1$.

Let $\psi\colon\sU_{2}\to A$ be the morphism obtained
    via the universal family $\mu\colon\sU_{2}\to\sH_{2}$. Next lemma 
    is necessary to prove the finitess of $\psi$ outside $\cup_{i=1}^{s}
    \beta'_i\subset A$.

\begin{lem}
\label{lem:Cb}
Let $\overline{l}$ be a general uni-secant line of $C$ and
$l_b\subset \mP^2$ the image of $\overline l$ by the double projection from a point $b$.
For a general point $b\not \in C$,
$\deg C_b=d$ and $C_b\cup l_b$ has only simple nodes.
For a general point $b$ of $C$,
$\deg C_b=d-2$ and $C_b\cup l_b$ has only simple nodes.
\end{lem}

\begin{proof}
The claims for $\deg C_b$ follows from Propositions 
\ref{prop:proj} (1-1) and \ref{prop:Cd2}.
As for the singularity of $C_b\cup l_b$, the claim follows from
simple dimension count. For simplicity, we only prove that
for a general point $b$, $C_b$ has only simple nodes.
By Proposition \ref{prop:Cd1}, we may assume that 
any multi-secant conic through $b$ is smooth, bi-secant and
intersects $C$ simply.
Let $\overline{q}$ be a smooth bi-secant conic through $b$.
We may assume that $\sN_{\overline{q}/B}\simeq \sO_{\mP^1}(1)^{\oplus 2}$.
Let $q'$ be the strict transform of $\overline{q}$ on $B'_b$. 
Let $\widetilde{B}'\to B'_b$ be the blow-up along $q'$,
$E_{q'}$ the exceptional divisor and $\widetilde{C}''$ the strict transform of
$C''_b$. 
Note that $E_{q'}\simeq \mP^1\times \mP^1$
since $\sN_{q'/B'_b}\simeq \sO_{\mP^1}^{\oplus 2}$.
Then $C_b$ has simple nodes at the image of $q'$ if and only if
the two points in $E_{q'}\cap \widetilde{C}''$ does not belong to the same
ruling with the opposite direction to a fiber of $E_{q'}\to q'$. 
Let $\widetilde{B}_{\overline{q}}\to B$ be the blow-up along $\overline{q}$,
$E_{\overline{q}}$ the exceptional divisor and $\widetilde{C}$ the
strict transform of $C$.
It is easy to see that 
a ruling of $E_{\overline{q}}$ 
with the opposite direction to a fiber of $E_{\overline{q}}\to \overline{q}$
corresponds to that 
of $E_{q'}$ with the opposite direction to a fiber of $E_{q'}\to q'$.
Thus
$C_b$ has simple nodes at the image of $q'$ if and only if
the two points in $E_{\overline{q}}\cap \widetilde{C}$ does not belong to the same
ruling with the opposite direction to a fiber of $E_{\overline{q}}\to \overline{q}$. 
We can show that this is the case for a general $b$ by simple dimension count.
\end{proof} 

From now on, we assume $d\geq 5$ throughout the paper
since we need Proposition \ref{prop:Cd2}.

We do not have the finiteness of $\psi$ all over $A$.
To obtain a finite morphism,
we blow-up $A$ more in \ref{subsubsection:special}. Till now we can
prove:

\begin{prop}\label{fuoripiatto}
 $\psi$ is finite of degree $n:=\frac{(d-1)(d-2)}{2}$ 
 and flat outside $\cup_{i=1}^{s} \beta'_i$.     
 \end{prop}

\begin{proof}  
    Let
    $a\in A\setminus\cup_{i=1}^{s}\beta'_{i}$ and set $b:=f(a)$.
    If $b\not \in C$, then the finiteness of $\psi$ over $a$
    follows from Corollary \ref{cor:finite}. 
    Moreover, by Lemma \ref{lem:Cb},
    the number of conics through a general $a$ is $n$.
    Thus $\deg \psi=n$. We will prove that $\psi$ is finite over 
    $a\in E_C\setminus \cup_{i=1}^{s} \beta'_i$.
    Once we prove this,
    the assertion follows. Indeed,
    $\sU_2$ is Cohen-Macaulay since 
    $\sH_{2}$ is smooth
    and any fiber of $\sU_2\to\sH_2$ is reduced, thus  
    $\psi$ is flat.
    
    Let $a\in E_C\setminus \cup_{i=1}^{s} \beta'_i$.
    The assertion is equivalent to that
    only finitely many conics belonging to $L_b$ pass through $a$.
    If $b\not \in \cup_{i=1}^{s} \beta_i$, then $L_b$ is irreducible.
    If $b\in \cup_{i=1}^{s} \beta_i$, then $L_b=L'_b+e_i$,
    where $L'_b$ is the strict transform of $\eta(L_b)$ whence is irreducible.
    Note that almost all the conics belonging to $e_i$
    does not pass through $a \not \in \cup_{i=1}^{s} \beta'_i$.
    Let $S_{b}\subset {A}$ be the locus swept by the 
    conics of the family $L_{b}$ if $b\not \in \cup_{i=1}^{s} \beta_i$, 
    or the locus swept by the 
    conics of the family $L'_{b}$ if $b\in \cup_{i=1}^{s} \beta_i$.
    $S_b$ is irreducible.
    Let $\overline{S}_b:=f(S_b)$,
    $\overline{S}'_b$ and $\overline{S}''_b$ 
    the strict transforms of $\overline{S}_b$ on $B_b$ and 
    $B'_b$ respectively.
    Then $\overline{S}''_b=\pi_{2b}^* C_b$. 
    Let $d_b:=\deg C_b$. By Proposition \ref{prop:proj} (1-1),
    $d_b=d-2$ if $b\not \in \cup_{i=1}^{s} \beta_i$, or
    $d-3$ if $b\in \cup_{i=1}^{s} \beta_i$.
    Since $\overline{S}''_b\sim d_b L$ and $L=H-2E'_b$, we have
    $\overline{S}'_{b|{E_b}}$ is a curve of degree $2d_b$ in $E_b\simeq \mP^2$.
    
    Since $A$ is obtained from $B_b$ by blowing up $C'_b$ and then
    contracting the strict transform of $E_b$,
    a point $a$ over $b$ corresponds to a line $l_a$ in $E_b$ through 
    $t:=E_b\cap C'_b$. 
    If $C''_b$ does not intersect fibers of $\pi_{2b}$ contained in  $E'_b$, 
    then
    $\overline{S}'_{b|{E_b}}$ is irreducible. Thus no $l_a$ is 
    contained in $\overline{S}'_{b|{E_b}}$ and we are done.
    Assume that  
    $C''_b$ intersects a fiber $l'$ of $\pi_{2b}$ contained in $E'_b$. 
    By Claim \ref{cla:reduced} (3),
    $b\not \in B_{\varphi}$ nor $b \not \in \cup_{i=1}^{s} \beta_i$
    for a general $C$.
    Since $L_b$ is irreducible,
    it suffices to prove the finiteness and nonemptyness of the set of conics
    through a general point $a$ over $b$. 
    Equivalently, we have only to show that 
    a general $l_a$ intersects 
    $\overline{S}'_{b|E_b}$ outside $t$.    
    Since $l'$ intersects $C''_b$ simply at one point,
    $C_b$ is smooth at the image $t'$ of $l'$ on $\mP^2$.     
    Thus   
    $\overline{S}'_{b|{E_b}}=C'''_b+l$, 
    where $C'''_b$ and $l$ are the strict transforms
    of $C_b$ and $l'$. Note that $C'''_b$ is smooth at $t$ and
    $\deg C'''_b=2d_b-1=2d-5\geq 5$ by $d\geq 5$.
    Thus a general $l_a$ intersect $C'''_b$ outside $t$.
\end{proof}

\begin{rem}
Though we do not need it later,
we describe the fiber of $\psi$ 
over a general point $a \in E_C\setminus \cup_{i=1}^{s} \beta'_i$
for reader's convenience.

Set $b:=f(a)$. 
As in the proof of Proposition \ref{fuoripiatto}, 
a point $a$ over $b$ corresponds to a line $l_a$ in $E_b$
passing through $E_b\cap C'_b$.
By Lemma \ref{lem:Cb},
it holds that $\deg C_b=d-2$ and $C_b$ has $\frac{(d-3)(d-4)}{2}$ simple nodes
for a general $b\in C$. This means that  
$\frac{(d-3)(d-4)}{2}$ tri-secant conics pass through $b$. 
By Proposition \ref{classificazioneconicheinA}, corresponding to a tri-secant conic
$\overline q$,
there is a unique conic $q$ on $A$ containing the fiber of $E_C$
over $b$ and such a conic on $A$ contains $a$. 
Thus we obtain $\frac{(d-3)(d-4)}{2}$ conics through $a$.
By definition of $L_b$, these conics do not belong to $L_b$.

We need more
$n-\frac{(d-3)(d-4)}{2}=2d-5$ conics through $a$.
We show that there exist $2(d-2)-1$ conics through $a$ on $A$
coming from the family parameterized by $L_b$.
We use the notation of the proof of Proposition \ref{fuoripiatto}.
For a general $b\in C$,
$C''_b$ does not intersect fibers of $\pi_{2b}$ contained in $E'_b$. 
Thus $\overline{S}'_{b|E_b}$ is an irreducible curve
of degree $2(d-2)$ on $E_b$.
Thus there are $2(d-2)$ intersection points of
$\overline{S}'_{b|E_b}$ and $l_a$. Among those, the intersection
point $C'_{b}\cap E_b$ does not correspond to a conic on $A$
through $a$ since it comes from the tangent of $C$. 
Thus we have $2(d-2)-1$ conics as desired.
\end{rem}


We need to study mutual intersection of a conic and a line in special cases. 
Let $\sF\subset\sH_{2}\times\sH_{1}$ be
    the image in $\sH_2\times \sH_1$ of 
the inverse image of $((\cup \beta'_i)\times \sH_1) \cap \sU_1$; that 
is
\[
\sF:=\{([q],[l])\mid q \cap \beta'_i\cap l\not =\emptyset\}. 
\]
A point $([q],[l])\in\sF$ iff
i) $l=l_{ij}$ and 
$q\cap \beta'_i\not =\emptyset$ or ii)
$l\not =l_{ij}$, $l\cap\beta'_i\neq\emptyset$ and 
$q \cap \beta'_i\cap l\not =\emptyset$. 
For every $i=1,\ldots ,s$, $j=1,2$ the family of those $([q],[l])$
which satisfies i) or ii) has dimension one and clearly does not dominate
$\sH_{1}$. 

\begin{cor}\label{divisoriali}
Any component of $\widehat{\sD}_1$ which is not contained in
$\sF$ dominates $\sH_1$.
Moreover, any non-divisorial component of $\widehat{\sD}_1$
outside $\sF$
$(\text{if it exists})$ is
a one-dimensional component whose generic point
parameterizes reducible conics, namely,
a one-dimensional component of 
\[
\{([q],[l])\mid l\subset q\}.
\]
\end{cor}

\begin{rem}
Here we leave the possibility that
a one-dimensional component whose generic point
parameterizes reducible conics is contained in 
a divisorial component of $\widehat{\sD}_1$.
We, however, prove that 
this is not the case in Corollary \ref{noncontenere}.
Hence, finally, the fiber of $\widehat{\sD}_1\to \sH_1$ over
a general $[l]\in \sH_{1}$ parameterizes 
conics which properly intersect $l$. 
\end{rem}

\begin{proof} 
    By Proposition \ref{fuoripiatto}, $\sU_2\to A$ is finite and flat outside $\cup
    \beta'_i$. Then $\sU_{2}\times\sH_{1}\to A\times\sH_{1}$ is flat
    outside $(\cup \beta'_i)\times \sH_1$. By base change,
    $\sU'_{1}\rightarrow \sU_{1}$ is flat and finite outside 
    $((\cup \beta'_i)\times \sH_1) \cap \sU_1$. Then every 
    irreducible component of $\sU'_1$ which is not mapped to 
$((\cup \beta'_i)\times \sH_1) \cap \sU_1$ is two-dimensional, and
dominates $\sU_1$, hence dominates $\sH_1$.
Therefore
any component of $\widehat{\sD}_1$ which is not contained in
$\sF$ dominates $\sH_1$.

We find a possible non-divisorial component of
$\widehat{\sD}_1$ outside $\sF$.
Let $\gamma\subset \sU'_1$ be a curve mapped to a point, say, 
$([q],[l])$ 
on $\sH_2\times \sH_1$. The image of $\gamma$
on $A$ is an irreducible component of $q$, say, $q_1$. 
The image of $\gamma$ on $\sU_1$ 
is $q_1 \times [l]$, thus $q_1$ is also an irreducible
component of $l$.
We have the following three possibilities:

\begin{enumerate}[$(1)$]
\item
$l$ is irreducible, hence $q_1=l$ and
$q=l\cup m$, where $m$ is another line. 
Such $([q],[l])$ form the one-dimensional family of reducible conics,
\item
$l=l_{ij}$ and $\beta'_i\subset q$.
Namely $[q]\in e_i$, or $q=\beta'_i\cup \alpha\cup \zeta_{ik}$,
where $\alpha$ is the strict transform of a line on $B$
intersecting $\beta_i$ and $C$ outside $\beta_i\cap C$, or
\item
$l=l_{ij}$ and $\zeta_{ij}\subset q$ and
$f(q)$ is a tri- or quadri-secant conic of $C$
such that $p_{ij}\in f(q)$.
\end{enumerate}
Thus we have the second assertion.
\end{proof}

\begin{nota}
Let ${\sD}_1 \subset \sH_2\times \sH_1$ be 
the divisorial part of $\widehat{\sD}_1$. 
Since $\sH_1$ is a smooth curve ${\sD}_1 \to \sH_1$ is flat.  
Let ${D}_l$ be the fiber of 
${\sD}_1\to \sH_1$ over $[l]\in \sH_1$. Clearly we can write
$D_{l}\hookrightarrow \sH_{2}$.
\end{nota}

Next two lemmas are basic to understand the geometry of $\sH_{2}$.

\begin{lem}
\label{lem:twolines}
Let $\overline{l}_1$ and $\overline{l}_2$ be 
two general secant lines of $C$
such that ${\overline{l}_1}\cap {\overline{l}_2}=\emptyset$.
Let $B\dashrightarrow Q \dashrightarrow \mP^2$
be the successive linear projections 
from $\overline{l}_1$ and then
the strict transform of $\overline{l}_2$ on $Q$.
Let $\overline{l}$ be another general secant line of $C$,
and $C'$ and $\overline{l}'\subset\mP^{2}$  be the images of $C$
and $\overline{l}$ respectively. 
Then $C\cup \overline{l}\dashrightarrow C'\cup \overline{l}'$ 
is birational and 
$C'\cup \overline{l}'$ has only simple nodes as its singularities,
where, by birational, we means that $\deg C'\cup \overline{l}'=d-1$.
In particular $($since $\deg C'=d-2$ and $C'$ is rational$)$
$C'$ has $\frac{(d-3)(d-4)}{2}$ simple nodes, equivalently,
there exist 
$\frac{(d-3)(d-4)}{2}$ bi-secant conics of $C$ intersecting both
$\overline{l}_1$ and $\overline{l}_2$.
\end{lem}

\begin{proof}
We show the assertion using the inductive construction of $C=C_d$.
The assertion follows for $d=3$ directly.
Consider a smoothing from $C_{d-1}\cup \overline{m}$ to $C_d$.
Let $\overline{m}_1$ and $\overline{m}_2$ 
two general secant lines of $C_{d-1}$
such that ${\overline{m}_1}\cap {\overline{m}_2}=\emptyset$.
Let $B\dashrightarrow Q \dashrightarrow \mP^2$
be the successive linear projections 
from $\overline{m}_1$ and then from
the strict transform of $\overline{m}_2$ on $Q$.
Let $\overline{r}$ be another general secant line of $C_{d-1}$,
and $C'_{d-1}, \overline{m}'$ and $\overline{r}' \subset \mP^2$
be the images of $C_{d-1}$, $\overline{m}$ and $\overline{r}$ respectively. 
Then we have only to show that
$C_{d-1}\cup \overline{m}\cup \overline{r} 
\dashrightarrow C'_{d-1}\cup \overline{m}'\cup \overline{r}'$ 
is birational and 
$C'_{d-1}\cup \overline{m}'\cup \overline{r}'$ has 
only simple nodes as its singularities
assuming $C_{d-1}\cup \overline{r}\dashrightarrow C'_{d-1}\cup \overline{r}'$ is birational and $C'_{d-1}\cup \overline{r}'$ has only simple nodes as its singularities.

Since $\overline{m}$ is also general,
$C_{d-1}\cup \overline{m}\dashrightarrow C'_{d-1}\cup \overline{m}'$ is birational and $C'_{d-1}\cup \overline{m}'$ has only simple nodes as its singularities.
Thus $C_{d-1}\cup \overline{m}\cup \overline{r} 
\dashrightarrow C'_{d-1}\cup \overline{m}'\cup \overline{r}'$ 
is clearly birational.
To show $C'_{d-1}\cup \overline{m}'\cup \overline{r}'$ has 
only simple nodes as its singularities,
it suffices to prove that
there are no secant conics of $C_{d-1}$
intersecting all the
$\overline{m}_1$, $\overline{m}_2$, $\overline{m}$ and $\overline{r}$.
This follows from the fact that a secant conic $\overline{q}$ of $C_{d-1}$
intersects finitely many secant lines of $C_{d-1}$
by $M(\overline{q})\not \subset M(C_{d-1})$.
\end{proof}

\begin{lem}
\label{lem:twolines'}
Let $\overline{l}_0$ be a general secant line of $C$.
Let $B\dashrightarrow Q \dashrightarrow \mP^2$
be the successive linear projections 
from $\overline{l}_0$ and then
the strict transform of a $\beta_i$ on $Q$.
Let $\overline{l}$ be another general secant line of $C$,
and $C'$ and $\overline{l}'\subset\mP^{2}$  be the images of $C$
and $\overline{l}$ respectively. 
Then $C\cup \overline{l}\dashrightarrow C'\cup \overline{l}'$ 
is birational and 
$C'\cup \overline{l}'$ has only simple nodes as its singularities. 
In particular $($since $\deg C'=d-3$ and $C'$ is rational$)$
$C'$ has $\frac{(d-4)(d-5)}{2}$ simple nodes,
equivalently,
there exist $\frac{(d-4)(d-5)}{2}$ bi-secant conics of $C$ intersecting
$\beta_i$ and $\overline{l}_0$ except conics containing $\beta_i$.

\end{lem}

\begin{proof}
Similarly to the previous lemma, 
we show the assertion using the inductive construction of $C=C_d$.
The assertion follows for $d=4$ directly.
Consider a smoothing from $C_{d-1}\cup \overline{m}$ to $C_d$.
Let $\overline{m}_0$ be a general secant line of $C_{d-1}$,
and $\beta$ a bi-secant line of $C_{d-1}\cup \overline{m}$
different from two lines through $C_{d-1}\cap \overline{m}$. 
Let $B\dashrightarrow Q \dashrightarrow \mP^2$
be the successive linear projections 
from $\overline{m}_0$ and then
the strict transform of $\beta$ on $Q$.
Let $\overline{r}$ be another general secant line of $C_{d-1}$,
and $C'_{d-1}, \overline{m}'$ and $\overline{r}' \subset \mP^2$
be the images of $C_{d-1}$, $\overline{m}$ and $\overline{r}$ respectively. 

First we suppose that $\beta$ is a bi-secant line of $C_{d-1}$.
Then we have only to show that
$C_{d-1}\cup \overline{m}\cup \overline{r} 
\dashrightarrow C'_{d-1}\cup \overline{m}'\cup \overline{r}'$ 
is birational and 
$C'_{d-1}\cup \overline{m}'\cup \overline{r}'$ has 
only simple nodes as its singularities
assuming $C_{d-1}\cup \overline{r}\dashrightarrow C'_{d-1}\cup \overline{r}'$ is birational and 
$C'_{d-1}\cup \overline{r}'$ has only simple nodes as its singularities.
The proof is the same as that of Lemma \ref{lem:twolines}, so we omit it.

Next suppose that $\beta$ is a uni-secant line of $C_{d-1}$ 
intersecting $\overline{m}$ outside $C_{d-1}\cap \overline{m}$.
Note that, by the projection $B\dashrightarrow \mP^2$,
$\overline{m}$ is contracted to a point.
Moreover,
$\beta$ is a general uni-secant line since so is $\overline{m}$.
Thus, by Lemma \ref{lem:twolines}, 
$C_{d-1}\cup \overline{m}\cup \overline{r} 
\dashrightarrow C'_{d-1}\cup \overline{m}'\cup \overline{r}'$ 
is birational and 
$C'_{d-1}\cup \overline{m}'\cup \overline{r}'$ has 
only simple nodes as its singularities.
\end{proof}

Now we reach the precise description of $\sH_2$.

\begin{thm}
\label{thm:H_2}
$\eta\colon \sH_2\to \mP^2$
is the blow-up at $c_1,\dots, c_s$
and $e_i$ are $\eta$-exceptional curves. Moreover
$\sH_2$ has the following properties:
\begin{enumerate}[$(1)$]
\item
\[
D_l\sim (d-3)h-\sum_{i=1}^{s} e_i,
\]
where $h$ is the strict transform of a general line on $\mP^2$.
\item
\[
h^1(\sH_{2},\sO_{\sH_{2}}((d-4)h-\sum_{i=1}^{s} e_i))=0.
\]
\item
$|D_l|$ is base point free.
In case of $d=5$, the image of $\Phi_{|D_l|}$ is $\check{\mP}^2$.
In case of $d\geq 6$, $D_l$ is very ample and
$|D_l|$ embeds $\sH_2$ into $\check{\mP}^{d-3}$.
\item
If $d\geq 6$, then 
$\sH_2\subset \check{\mP}^{d-3}$ is 
projectively Cohen-Macaulay and is the intersection of cubics . 
\end{enumerate}
\end{thm}

\begin{rem}
\begin{enumerate}[(1)]
\item
If $d\geq 6$, then $\sH_2\subset \check{\mP}^{d-3}$ is so called the 
{\em{White surface}} (see \cite{White} and \cite{Gimi}).
In \cite{Man},
the white surface attains the maximal degree among 
projectively Cohen Macaulay rational surfaces in
a fixed projective space.
\item
We use the dual notation $\check{\mP}^{d-3}$
for later convenience.
\end{enumerate} 
\end{rem}

\begin{proof}
(1) 
Let $\pi_C\colon C\times C\to S^2 C$ be the natural projection
and $L'_b$ a ruling of $C\times C$ in one fixed direction 
such that $\pi_C(L'_b)=\eta(L_b)$.
By applying the Bertini theorem to $|L'_b|$, 
we see that $\pi_C^*\eta(D_l)$ and $L'_b$ intersect simply for a general
$b\in C$ whence $\eta(D_l)$ intersects $\eta(L_{b})$ simply since 
$\pi_C$ is \'etale at $\pi_C^*\eta(D_l)\cap L'_b$. Then $D_l$ intersects $L_b$ 
simply since $\eta$ is isomorphic at $D_{l}\cap L_{b}$.
For a general $b\in C$, we consider the double projection
$\pi_b\colon B\dashrightarrow \mP^2$ from $b$ 
as in Proposition \ref{prop:proj} (1).
Let $\overline{l}$ be a general line on $B$ intersecting $C$
and $l_{b}$ the image of $\overline l$ by $\pi_{b}$. 
We can assume that $b\not =c:=C\cap {\overline l}$.
Obviously $l_{b}$ is a line. By Lemma \ref{lem:Cb}, 
$\deg C_b=d-2$ and the curve $C_b\cup l_b$ has only simple nodes.
Hence the number of points in $C_b\cap l_b$ is
$d-2$, one of them counts for the unique conic through $b$ and $c$.
This last conic gives a conic on $A$ which does not intersects $l$.
The other $d-3$ points count for elements of $D_{l}$.
Then $\eta(D_l)$ is a curve of degree $d-3$.

Let $\overline{l}_1$ and $\overline{l}_2$ be 
two general secant lines of $C$
such that ${\overline{l}_1}\cap {\overline{l}_2}=\emptyset$.
By Lemma \ref{lem:twolines},
$\# (D_{l_1}\cap D_{l_2})=\frac{(d-3)(d-4)}{2}$.
This immediately gives for the intersection product
$D_{l_1}\cdot D_{l_2}\geq \frac{(d-3)(d-4)}{2}$.
On the other hand, $D_l\cap e_i\not =\emptyset$ for a general $l$
since $D_l\cap e_i$ contains the point corresponding to
a marked conic
$(\beta_i\cup \alpha, \beta_{i|C})$, where $\alpha$ is the unique line
intersecting $\beta_i$ and $l$.
Moreover, for two general $l_1$ and $l_2$, 
$D_{l_1}\cap D_{l_2}\cap e_i=\emptyset$. Now the curves $e_{i}$ have
negative self intersection then $D_{l_1}\cdot D_{l_2}\leq
(d-3)^2-s=\frac{(d-3)(d-4)}{2}$.
Therefore
$D_{l_1}\cdot D_{l_2}=\frac{(d-3)(d-4)}{2}$. 
Moreover $e_{i}^{2}=-1$ and since 
$e_{i}\cap e_{j}=\emptyset$ we obtain that $\eta\colon \sH_2\to \mP^2$
is the blow-up at $c_1,\dots, c_s$.
Consequently,
$D_l\sim (d-3)h-\sum_{i=1}^{s} e_i$ for a general $[l] \in \sH_1$, and,
by the flatness of $\sD_1 \to \sH_1$,
that holds for any $[l] \in \sH_1$. 

(2) Let $L'_{p_{ij}}=L_{p_{ij}}-e_i$ (note that $e_i\subset L_{p_{ij}}$).
We see that
$L'_{p_{ij}}\subset D_{l_{ij}}$ and 
$D_{l_{i1}}-L'_{p_{i1}}=D_{l_{i2}}-L'_{p_{i2}}$,
which we denote by $D_{\beta_i}$. 
Note that
\[
D_{\beta_i}\sim (d-4)h-\sum_{k\not =i} e_k.
\]
It is easy to see that $D_{\beta_i}$ have the following properties:
\begin{eqnarray}
D_{\beta_i}\cap e_i=\emptyset.
\label{eqnarray:1}\\
D_{\beta_i}\cap D_{\beta_j}\cap D_{\beta_k}=\emptyset.
\label{eqnarray:2}
\end{eqnarray}
We only prove (\ref{eqnarray:1}).
Since $D_{\beta_i}\cap e_i\not =\emptyset$ would imply that $e_i$ is a component of
$D_{\beta_i}$, it suffices to prove that, for a general $l$,
$D_{\beta_i}\cap D_l$ does not contain
a point of $e_i$.
By Lemma \ref{lem:twolines'},
$D_{\beta_i}\cap D_l$ contains
$\frac{(d-4)(d-5)}{2}$ points
corresponding to bi-secant conics intersecting
$\beta_i$ and $l$ except conics containing $\beta_i$.
On the other hand,
we have $D_l\cdot D_{\beta_i}=\frac{(d-4)(d-5)}{2}$,
thus the conics we count in Lemma \ref{lem:twolines'}
correspond to all the intersection of $D_{\beta_i}\cap D_l$.
Consequently, $D_{\beta_i}\cap D_l$ does not contain
a point of $e_i$.

By (\ref{eqnarray:1}) and the trivial equality
\[
(d-4)h-\sum_{i\geq k+1} e_i=D_{\beta_k}+ e_1+\cdots+ e_{k-1},
\]
\noindent 
we obtain $e_k\not \subset \Bs |(d-4)h-\sum_{i\geq k+1} e_i|$.

Since $\sO_{\sH_{2}}((d-4)h-\sum_{i\geq k+1}
e_i)\otimes_{\sO_{\sH_{2}}}\sO_{e_{k}}\simeq \sO_{e_{k}}$ we have that
\[
H^0(\sH_{2},\sO_{\sH_{2}}((d-4)h-\sum_{i\geq k+1} e_i)) \to  
H^0(\sH_{2}, \sO_{e_k})
\]
is surjective.
Hence by the exact sequence
\[
0\to \sO_{\sH_{2}}((d-4)h-\sum_{i\geq k} e_i)\to 
\sO_{\sH_{2}}((d-4)h-\sum_{i\geq k+1} e_i)\to \sO_{e_k}\to 0,
\]
we have
$H^1(\sH_{2}, \sO_{\sH_{2}}((d-4)h-\sum_{i=1}^{s} e_i))
\simeq H^1(\sH_{2}, \sO_{\sH_{2}}(d-4)h)$. Since 
it is easy to see that $h^1(\sH_{2},\sO_{\sH_{2}}(d-4)h)=0$,
we have (2).

(3) Since no conic on $A$ intersects all the line on $A$,
$|D_l|$ has no base point.
In case $d=5$, the image of $\Phi_{|D_l|}$ is $\mP^2$ by $(D_l)^2=1$.

Assuming $d\geq 6$, we prove that $D_l$ is very ample.
By (2) and \cite[Theorem 3.1]{DG},
it suffices to prove
that 
\[
h^0(\sH_{2}, \sO_{\sH_{2}}(h-\sum_{j=1}^{d-3} e_{i_j}))=0
\] 
for any set of $d-3$ exceptional curves
$e_{i_1},\dots, e_{i_{d-3}}$.
Assume by contradiction that
there exists an effective divisor 
$L\in |h-\sum_{j=1}^{d-3} e_{i_j}|$ 
for a set of $d-3$ exceptional curves $e_{i_1},\dots, e_{i_{d-3}}$.  
By $\frac{(d-2)(d-3)}{2}-(d-3)\geq 3$, we find at least three $e_i$
such that
$i\not \in \{j_1,\dots,j_{d-3} \}$.
For an $i\not \in \{j_1,\dots,j_{d-3}\}$,
noting $D_l\sim D_{\beta_i}+h-e_{i}е$, 
$D_l\cdot L=0$, and $L\cdot (h-e_i)>0$, we have
$L\subset D_{\beta_i}$.
This contradicts (\ref{eqnarray:2}) since the number of $i$ 
such that $i\not \in \{j_1,\dots,j_{d-3} \}$ is at least $3$.

We show that
$h^0(\sH_{2}, \sO_{\sH_{2}}(D_l))=d-2$.
By the Riemann-Roch theorem,
$\chi(\sO_{\sH_{2}}(D_l))=d-2$.
Since $h^2(\sH_{2}, \sO_{\sH_{2}}(D_l))=
h^0(\sH_{2}, \sO_{\sH_{2}}(-D_l+K_{\sH_2}))=0$,
we see that
$h^0(\sH_{2}, \sO_{\sH_{2}}(D_l))=d-2$ is equivalent to
$h^1(\sH_{2}, \sO_{\sH_{2}}(D_l))=0$.
Since
$|D_l|$ has no base point, so is $|(d-3)h-\sum_{i\geq k+1} e_i|$.
Thus
the proof that $h^1(\sH_{2}, \sO_{\sH_{2}}(D_l))=0$ 
is almost the same as the above one showing (2)
and we omit it.

(4) follows from
\cite[Proposition 1.1]{Gimi}.
\end{proof}

\begin{rem}
In case of $d=5$,
the morphism defined by $|D_l|$ contracts
three curves $D_{e_i}$ ($i=1, 2, 3$),
which are nothing but the strict transforms of
three lines passing through two of $c_j$.
Namely, the composite $S^2 C\leftarrow \sH_2 \rightarrow \check{\mP}^2$
is the Cremona transformation. 
\end{rem}

The following corollary contains the nontrivial 
result that for a general $[l]\in\sH_{1}$, $D_{l}$ parameterizes
conics which properly intersect $l$.

\begin{cor}\label{noncontenere} 
    For a general $[l]\in\sH_{1}$, $D_{l}$
    does not contain any point corresponding to the line pairs $l\cup m$
    with $[m]\in\sH_{1}$.
   
\end{cor}
\begin{proof} Fix $[m]\in\sH_{1}$ such that $l\cup m$ is a line pair. 
    If $({\overline m},b)$ is the marked line given by $m$
    then we have $d-2$ line pairs $l\cup m$, $l_{1}\cup m$,\ldots
    ,$l_{d-3}\cup m$. Since $L_{b}\sim h$ then $h\cdot D_{l}= d-3$ and
    definitely $[l_{1}\cup m],\ldots
    ,[l_{d-3}\cup m]\in D_{l}$. Thus $[l\cup m]\not\in D_{l}$.
\end{proof}


\subsection{Varieties of power sums for special non-degenerate
quartics $F_4$}
\label{subsection:VSPA}
~

In Proposition \ref{fuoripiatto} we have seen that $\psi\colon \sU_{2}\rightarrow A$ is
finite and flat outside $\cup_{i=1}^{n}\beta'_{i}$. We can modify 
the morphism $\psi\colon \sU_{2}\rightarrow A$ to obtain a finite one.
See Proposition \ref{finitezzafinale}, which is the goal of 
\ref{subsubsection:special}.
This and our understanding of the geometry of $\sH_{2}$ 
give an important morphism whose target is $\VSP(F_4,n;\sH_2)$: see
Theorem \ref{diretto}.

\subsubsection{Special blow-up $\widetilde A$ of $A$}
\label{subsubsection:special}
Similarly to (\ref{eq:familyB1}),
we consider the following diagram:
\begin{equation}
\label{eq:family2}
\xymatrix{\sU'_2\subset {\sU}_2\times \sH_2 \ar[d] 
\ar[r]^{(\psi, \mathrm{id})} & A \times \sH_2 \supset \sU_2 \ar[d]\\
\widehat{\sD}_2 \subset \sH_2\times \sH_2 \ar[r] & \sH_2.} 
\end{equation}

Let $\sU'_2\subset \sU_2\times \sH_2$ be the pull-back of $\sU_2$
and $\widehat{\sD}_2$ the image of $\sU'_2$ on $\sH_{2}\times \sH_{2}$.
Similarly to the investigation of the diagram (\ref{eq:familyB1}),
we see that the image $\sF'$ in $\sH_2\times \sH_2$ of the inverse image of 
$\cup_{i=1}^{n}\beta'_i\times \sH_2$ is not divisorial nor does not dominate $\sH_2$.
Moreover, 
any component of $\widehat{\sD}_2$ outside $\sF'$ dominates $\sH_2$, and
is divisorial or possibly the diagonal of $\sH_2\times \sH_2$.
Note that dislike the diagram (\ref{eq:familyB1}),
there is no other non-divisorial component in this case.
Compare the proof of Corollary \ref{divisoriali}.
Here we leave the possibility that 
the diagonal of $\sH_2\times \sH_2$ is contained in
the divisorial component of $\widehat{\sD}_2$.
We, however, prove this is not the case in 
Lemma \ref{noselfintersection}. 

Let $\sD_2 \subset \sH_2\times \sH_2$ be the union of the
divisorial components of $\widehat{\sD}_2$ with reduced structure.
$\sD_2$ is Cartier since $\sH_2\times \sH_2$ is smooth.
$\sD_2 \to \sH_2$ is flat 
since $\sD_2$ is Cohen-Macaulay, $\sH_2$ is smooth and 
$\sD_2 \to \sH_2$ is equi-dimensional.  
Let $D_q$ be the fiber of 
$\sD_2\to \sH_2$ over $[q]\in \sH_2$ via the projection to
the second factor.

We are almost ready to define the 
modification of $\psi\colon \sU_{2}\rightarrow A$ we are looking for.
To find the range 
we consider the blow-up of $A$ along 
$\cup_{i=1}^{n}\beta'_{i}$ and we denote it by 
$\rho\colon{\widetilde A}\rightarrow A$.

\begin{lem}\label{sullebisecanti}
    $$
    \sN_{\beta'_{i}/A}=\sO_{\beta'_{i}}(-1)\oplus\sO_{\beta'_{i}}(-1).
    $$
\end{lem}
\begin{proof}
We prove the assertion by using the inductive construction of $C_d$.
The assertion is clear for $d=1$ since $C_1$ has no bi-secant line.
Suppose the assertion holds for $C_{d-1}$. 
Choose a general uni-secant line $\overline{l}\subset B$ of $C_{d-1}$.
Let $\overline{m}_1,\dots, \overline{m}_{d-2}$ be the lines on $B$ 
intersecting both $C_{d-1}$ and $\overline{l}$ outside 
$C_{d-1}\cap \overline{l}$.
Let $A'\to B$ be the blow-up along $C_{d-1}\cup \overline{l}$.
Note that the smoothing $C_{d-1}\cup \overline{l}$ to $C_d$ induces 
that of $A'$ to $A$.
Let $\widetilde{m}_i$ be
the strict transform of $\overline{m}_i$ on $A'$.
By the smoothing construction of $C_d$ from $C_{d-1}\cup \overline{l}$
and the assumption on induction, 
we have only to prove $\sN_{\widetilde{m}_i/A'}=
\sO_{\mP^1}(-1)\oplus\sO_{\mP^1}(-1)$.
Let $A'_1\to B$ be the blow-up along $\overline{l}$ and 
$A'_2\to A'_1$ the blow-up along the strict transform of $C_{d-1}$.
Denote by $m'_i$ and $m''_i$ the strict transform of $\overline{m}_i$ on 
$A'_1$ and $A'_2$ respectively.
Then $\sN_{\widetilde{m}_i/ A'}=\sN_{m''_i/ A'_2}$.
Since $m'_i$ is a fiber of $A'_1\to Q$ (cf. Proposition \ref{prop:proj} (2)), 
we have $\sN_{m'_i/A'_1}=\sO_{\mP^1}\oplus\sO_{\mP^1}(-1)$.
Let $F$ be the exceptional divisor of $A'_1\to Q$ and $F'$ the strict transform
of $F$ on $A'_2$. We may suppose $F$ and $C'_{d-1}$ intersect transversely,
thus $F'\to F$ is the blow-up at $d-2$ points $m'_i\cap C'_{d-1}$ 
($i=1,\dots, d-2$). Thus $F'\cdot m''_i=-1$ and $\sN_{m''_i/F'}=
\sO_{\mP^1}(-1)$, and this implies the assertion.
\end{proof}

We add the following piece of notation:
\begin{nota}
\begin{enumerate}[$(1)$] 
\item
$E_i:=\rho^{-1}(\beta'_{i})$. 
By Lemma \ref{sullebisecanti},
$E_{i}\simeq \mP^{1}\times\mP^{1}$,
\item
$f_{i}:=$ a general fiber of $\rho_{\vert E_{i}}\colon
 E_{i}\rightarrow\beta'_{i}$, 
\item $\gamma_{i}:=$ a general fiber of the other projection
$E_{i}\rightarrow\mP^{1}$, 
\item $\widetilde{E}_{C}:=$ the strict transform of $E_{C}$, and
\item 
$\widetilde{\zeta}_{ij}:=$ the strict transform of 
the fiber ${\zeta}_{ij}$ of ${E}_C$ over $p_{ij}\in C\cap\beta_{i}$,
\end{enumerate}
where $i=1,\ldots, s$ and $j=1,2$.
\end{nota}

The domain of the finite morphism is ${\widetilde \sU_{2}}:=
\sU_{2}\times_A{\widetilde A}$; in other words,
$\widetilde \sU_{2}$ is the blow-up of $\sU_2$ along 
$\sU_2\cap (\cup_{i=1}^{s}\beta'_{i}\times \sH_2$). We obtain that the
natural morphism
${\widetilde\sU_{2}}\rightarrow {\widetilde A}$ is finite after an
analysis of the morphism ${\sU_{2}}\rightarrow { A}$ in the
neighborhood of some special conics and via the suitable notion of
conic on ${ \widetilde A}$.
 
Note that, by Proposition \ref{prop:Cd} (5), 
there are $d-4$ lines $\alpha_1,\dots,\alpha_{d-4}$ distinct 
from $\beta_i$ and
intersecting both $C$ and $\beta_i$ outside $C\cap \beta_i$.
Set $t_k:=\alpha_k\cap C$.
Corresponding to $\alpha_k$,
there are two marked conics $(\alpha_k\cup \beta_i; p_{i1}, t_k)$
and $(\alpha_k\cup \beta_i; p_{i2}, t_k)$, which does not
belong to $e_i$ (by the choice of marking).
We denote by $\xi_{ijk}$ the conics on $A$ corresponding to
$(\alpha_k\cup \beta_i; p_{ij}, t_k)$, 
where $i=1,\dots, s$, $j=1,2$, and $k=1,\dots, d-4$.

\begin{lem}
\label{lem:0}
$\xi_{ijk}$ does not belong to $D_{\beta_i}$.
\end{lem}

\begin{proof}
By the projection from $\beta_i$, the image $\overline{q}$ 
of a general conic $q$
belonging to $D_{\beta_i}$ maps to a bi-secant line 
of the image $C'\subset Q$ of $C$, and $\alpha_k$ maps to a point $p_{ijk}$.
Let $p'_{ij}$ be the point of $C'$ corresponding to $p_{ij}$.
Let $F$ be the exceptional divisor over $\beta_i$, and
$F'$ the image of $F$ on $Q$. 
We say a ruling of $F'\simeq \mP^1\times \mP^1$ is
horizontal if it does not come from a fiber of $F\to \beta_i$. 
If $[\xi_{ijk}]\in D_{\beta_i}$,
then $\xi_{ijk}$ corresponds to a bi-secant line of $C'$,
which must be the horizontal ruling of $F'$
through $p'_{ij}$ and $p_{ijk}$.
By inductive construction of $C$, we can prove that
$p'_{ij}$ and $p_{ijk}$ do not lie on a horizontal ruling. 
Thus we have the claim. 
\end{proof}

\begin{defn}
    \label{conicabis}
    We say that a curve $\widetilde{q}\subset \widetilde{A}$ 
    is a {\em{conic}} on $\widetilde{A}$ if
\begin{enumerate}[(i)]
\item
$\widetilde{q}$ is connected and reduced,
\item
$-K_{\widetilde{A}} \cdot \widetilde{q}=2$,
\item
${{\widetilde E_{C}}}\cdot \widetilde{q}=2$, 
\item
$E_i\cdot \widetilde{q}=0$,
and
\item
$p_a(\widetilde{q})=0$.
\end{enumerate}
\end{defn}

Similarly to the case of conics on $A$, we know
there exists a unique two-dimensional component of 
the Hilbert scheme of $\widetilde{A}$
parameterizing conics on $\widetilde{A}$.
Let $\sH^{\widetilde{A}}_2$ be the normalization of the two-dimensional component. Similarly to the proof of Proposition \ref{equality},
we can show $\sH^{\widetilde{A}}_2$ is proper and 
there is a natural birational morphism 
$\sH^{\widetilde{A}}_2\to \overline{\sH}_2$.
Since $\sH_2\to \overline{\sH}_2$ is the normalization,
we have also a natural morphism $\sH^{\widetilde{A}}_2\to \sH_2$.
We do not need a full classification of conics on $\widetilde{A}$ but
only the following:

\begin{lem}
\label{lem:1}
\begin{enumerate}[$(1)$]
\item
There is a unique conic $\widetilde{q}$ on $\widetilde{A}$
corresponding to a conic $q$ on $A$ belonging to $e_i$,
and moreover, $\widetilde{q}$ is isomorphic to $q$ over
the component $\beta'_i$.
In particular,
$\sH^{\widetilde{A}}_2\to \sH_2$ is isomorphic near $e_i$.
\item
A conic belonging to $D_{\beta_i}$ is smooth near $\beta'_i$.
There is a unique conic $\widetilde{q}$ on $\widetilde{A}$
corresponding to a conic $q$ on $A$ belonging to $D_{\beta_i}$,
and, over $\beta'_i$,
$\widetilde{q}$ is isomorphic to the union of $q$ 
and the fiber of $E_i$ over $q\cap \beta'_i$.
In particular,
$\sH^{\widetilde{A}}_2\to \sH_2$ is isomorphic near $D_{\beta_i}$.
\end{enumerate}
\end{lem}
 
\begin{proof}
This follows from an explicit calculation as in the proof
of Proposition \ref{lineeA}.
For the first statement of (2), 
we use Proposition \ref{prop:Cd1} (5)
and Lemma \ref{lem:0}.
\end{proof}

Let $\Gamma:=\sU_2\cap (\cup_{i=1}^{s}\beta'_{i}\times \sH_2)$.
Outside $\cup_i \beta'_i\times e_i$,
$\Gamma$ is set-theoretically 
the disjoint union of
\[
\Gamma_i:=\{(x, [q])\mid [q]\in D_{\beta_i}, x\in q\cap \beta'_i\}\,
(i=1,\dots, s),
\]
which is a section of $\mu$ over $D_{\beta_i}$,
and 
\[
\Gamma_{ijk}:=\{(x,[\xi_{ijk}])\mid x\in \beta'_i\}\, (k=1,\dots,d-4,\, j=1,2).
\]

\begin{lem}
\label{lem:2}
Along $\Gamma_{ijk}$,
$\sU_2$ is smooth and $\Gamma$ is reduced. 
\end{lem}

\begin{proof}
To show that 
$\sU_2$ is smooth near $\Gamma_{ijk}$,
we have only to see that
the conic $\xi_{ijk}$ is strongly smoothable.
Note that $\sN_{\beta'_i/A}\simeq \sO_{\mP_1}(-1)^{\oplus 2}$,
$\sN_{\alpha'_k/A}\simeq \sO_{\mP_1}\oplus \sO_{\mP_1}(-1)$
and
$\sN_{\zeta_{i,3-j/A}}\simeq \sO_{\mP_1}\oplus \sO_{\mP_1}(-1)$.
We apply \cite[Theorem 4.1]{HH} by setting
$C=\beta'_i$, $D=\alpha'_k\cup \zeta_{i,3-j}$ and 
$S=(\alpha'_k\cap \beta'_i)\cup (\zeta_{i,3-j}\cap \beta'_i)$.
We check the conditions a) and b) of [ibid.].
The condition a) clearly holds.
The condition b) follows from the following two facts:
\begin{enumerate}[(1)]
\item
let $F$ be the exceptional divisor of 
the blow up of $B$ along $\alpha_k$.
Note that $F\simeq \mP^1\times \mP^1$.
We say a fiber of $F\to \mP^1$ in the other direction to
$F\to \alpha_k$ a horizontal fiber.
Then the intersection points of 
the strict transform of $C$ and $F$, and 
the strict transform of $\beta_i$ and $F$
do not lie on a common horizontal fiber.

This can be proved by the inductive construction of $C=C_d$
in a similar fashion to the proof of Lemma \ref{sullebisecanti}, and
\item
let $G$ be the exceptional divisor of 
the blow up of $A$ along $\zeta_{i,3-j}$.
Note that $G\simeq \mF_1$.
Then the intersection points of 
the strict transform of $\beta'_i$ and $G$
does not lie on the negative section of $G$.

This can be easily proved by noting
$\zeta_{i,3-j}$ is a fiber of $E$.
\end{enumerate}
Thus, by \cite[Theorem 4.1]{HH}, $\xi_{ijk}$ is strongly smoothable.

Second, we prove that $\Gamma$ is reduced along $\Gamma_{ijk}$.
We have only to prove that
${\sU}_2\to A$ is unramified along $\Gamma_{ijk}$
since then $\Gamma$ is the \'etale pull-back of $\beta'_i$ near
$\Gamma_{ijk}$, hence is reduced.

By the inductive construction of $C=C_d$ and
the following exact sequence:
\[
0\to \sN_{\beta'_i/A}\to \sN_{\xi_{ijk}/A|\beta'_i}\to T^1_S\to 0,
\]
we can prove that $\sN_{\xi_{ijk}/A|\beta'_i}\simeq \sO_{\mP^1}^{\oplus 2}$.
Thus $H^0(\sN_{\xi_{ijk}/A})\otimes \sO_{\xi_{ijk}}\to \sN_{\xi_{ijk}/A}$
is surjective at a point of $\Gamma_{ijk}$ since
it factor through the surjection
$H^0(\sN_{\xi_{ijk}/A|\beta'_i})\otimes \sO_{\beta'_i}\to 
\sN_{\xi_{ijk}/A|\beta'_i}$.
Thus 
${\sU}_2\to A$ is unramified along $\Gamma_{ijk}$.
\end{proof}

The next proposition contains the finitess result we need. 

\begin{prop}\label{finitezzafinale}
    $\widetilde \sU_2$ is Cohen-Macaulay and the natural morphism 
    ${\widetilde\psi}\colon{\widetilde
    \sU_{2}}\rightarrow {\widetilde A}$ is finite $($of degree
    $n:=\frac{(d-1)(d-2)}{2}$$)$.
In particular,
${\widetilde\psi}$ is flat.
\end{prop}

\begin{proof}
Lemma \ref{lem:1} shows that
$\widetilde{\sU}_2\to \sH_2$ is isomorphic to
the universal family of conics on $\widetilde{A}$
over $e_i$ and $D_{\beta_i}$.
Thus $\widetilde{\sU}_2$ is Cohen-Macaulay over
$e_i$ and $D_{\beta_i}$ since 
so are the fibers.  
Note that $\widetilde{\sU}_2\to \sU_2$ is the blow-up along $\Gamma_i$ near $\Gamma_i$
and is an isomorphism near $\beta'_i\times e_i$. 

Lemma \ref{lem:2} shows that
$\widetilde{\sU}_2\to \sU_2$ is the blow-up along $\Gamma_{ijk}$
near $\xi_{ijk}\times [\xi_{ijk}]$, and
$\widetilde{\sU}_2$ is smooth over $\Gamma_{ijk}$.
Thus $\widetilde{\sU}_2$ is Cohen-Macaulay.
To see $\widetilde \psi$ is finite,
we have only to note that
the inverse images of $\beta'_i\times e_i$ on $\widetilde{\sU}_2$ and
the exceptional divisor of $\widetilde \sU_2\to \sU_2$
are not contracted by $\widetilde \psi$.
\end{proof}

From now on in the section 3,
we assume that $d\geq 6$ and we consider $\sH_2\subset \check{\mP}^{d-3}$.

Consider  the following
diagram:

\begin{equation}
\label{iniezionehilbert}
\xymatrix{ & \widetilde{\sU_{2}}е \ar[dl]_{{\widetilde{\mu}}}
\ar[dr]^{{\widetilde\psi}}\\
 \sH_{2}е &  & {\widetilde A}. }
\end{equation}

\begin{defn}\label{attached}
    Let $\widetilde a$ be a point of $\widetilde A$.
    We say that $[\widetilde \psi^{-1}(\widetilde a)]\in 
\Hilb^{n} \check{\mP}^{d-3}$ is {\it{the cluster of conics attached to}} $\widetilde a$ and denote it by $[\sZ_{\widetilde a}]$. A conic $q$ such that $[q]\in \Supp \sZ_{\widetilde a}$ is 
called {\it{a conic attached to $\widetilde a$.}}
 \end{defn}

\begin{rem}
Though we do not need it later,
we describe the fiber of $\widetilde \psi$ over a general point 
${\widetilde a}\in E_i$ for some $i$
for reader's convenience.
In other words, we exhibit $n$ conics 
attached to $\widetilde a$.

Set $a:=\rho({\widetilde a})\in A$ and
    $b:= f(a)\in \beta_i$. We use notations of Proposition \ref{fuoripiatto}. 
    Since $\deg C_b=d-2$,
the number of bi-secant conics through $b$ not belonging to
the family $e_i$ is given by the number of double points of $C_b$, 
which is $\frac{(d-3)(d-4)}{2}$. Moreover
$2(d-4)$ conics $\xi_{ijk}$ through $a$.

The number of remaining conics is $3=n-\frac{(d-3)(d-4)}{2}-2(d-4)$. 
Such conics will belong to $e_i$.
By Lemma \ref{lem:1},
$\widetilde{\sU}_2\to \sH_2$ is isomorphic to
the universal family of conics on $\widetilde{A}$
over $e_i$.
Thus a desired conic on $A$ is the image of a conic $\widetilde{q}$ on 
$\widetilde{A}$
such that $\widetilde{a}\in \widetilde{q}$ and $\rho(\widetilde{q})$
belongs to $e_i$.
We show there are three such conics.
Let $S_i$ be the strict transform on $\widetilde{A}$ of 
the locus of lines intersecting $\beta_i$. 
Then it is easy to see that $S_{i|E_i}$ does not contain 
any fiber $\gamma_{i}$ of the second projection 
$\sigma_{i}\colon E_i\rightarrow\mP^{1}$.
Moreover $S_{i|E_i}\sim2\gamma_i+3f_{i}$. 
Let $\gamma'_{i}$ be the fiber of $\sigma_{i}$ through $\widetilde{a}$.
Then $\gamma'_{i}е$ intersect $S_i$ at three points.
Corresponding to these three points, there are
three lines on $B$ intersecting $\beta_i$.
Denote by $l_1$, $l_2$ and $l_3\subset A$ the strict transforms of these three 
lines. Then $\beta'_i\cup l_j$ $(j=1,2,3)$ are the conics on $A$ what we want.
\end{rem}

By Proposition \ref{finitezzafinale} and 
the universal property of Hilbert schemes, we obtain
a naturally defined map $\Psi\colon
{\widetilde A}\rightarrow{\rm{Hilb}}^{n}\check{\mP}^{d-3}$.
This is clearly injective because 
$n$ conics attached to a point ${\widetilde a}\in
    {\widetilde A}$ uniquely determines $\widetilde a$.

The task is to understand the image of $\Psi$.



\subsubsection{Morphism from $\tA$ to VSP}
~

To understand the image of $\Psi\colon
{\widetilde A}\rightarrow{\rm{Hilb}}^{n}\check{\mP}^{d-3}$ we
construct explicitly a quartic polynomial which plays the role of the
plane quartic in the Mukai's interpretation of $V_{22}$.

\begin{lem}\label{noselfintersection}
$\sD_2$ does not contain the diagonal of $\sH_2\times \sH_2$.
In particular we have the following:

let $\widetilde a$  be a general point of $\widetilde{A}$ and
$q_{1}, q_{2},\ldots ,q_{n} \in\sH_{2}$
the conics attached to $\widetilde a$. Then 
$$D_{q_i}([q_i])\neq 0$$
\noindent
for $1\leq i\leq n$.
\end{lem}
\begin{proof}
Here we assume $d\geq 3$. 
It suffices to prove that
$D_q([q])\neq 0$
for a general $[q]\in \sH_2$.
This is equivalent to that
the image ${D}^{\flat}_q$ on $\overline{\sH}_2$ of $D_q$ does not contain
$[\overline{q}]$. Note that $D^{\flat}_q$ is the closure
of the locus of multi-secant conics of $C$ intersecting 
properly $\overline{q}$.
Now the assertion follows from the inductive construction of $C_d$ from
$C_{d-1}\cup \overline{l}$.
From now on, we denote ${D}^{\flat}_{q}$ for $C_d$ by ${D}^{\flat}_{q,d}$.
If $d=3$, then $D_q\sim 0$, thus 
the assertion trivially true.
If ${D}^{\flat}_{q',d-1}([\overline{q}'])\neq 0$ 
for a general multi-secant conic $\overline{q}'$ of $C_{d-1}$,
then ${D}^{\flat}_{q,d}([\overline{q}])\neq 0$ for 
a general multi-secant conic $\overline{q}$ of $C_{d}$. 
\end{proof}

The proof of the following lemma is almost identical to the one of 
Theorem \ref{thm:H_2}; then we omit it: 

\begin{lem}
\label{lem:Dq}
$D_q\sim
2(d-3)h-2\sum_{i=1}^{e} e_i$
for a conic $q$,
namely, $D_q$ is a quadric section of $\sH_2\subset \check{\mP}^{d-3}$.
\end{lem}

We proceed to construct the quartic polynomial. By the seesaw theorem,
it holds that $\sD_2\sim p_1^*D_q+p_2^*D_q$.
Consider 
the morphism $\sH_2\times \sH_2$ into 
$\check{\mP}^{d-2}\times \check{\mP}^{d-3}$
defined by $|p_1^*D_l+p_2^*D_l|$, which 
is an embedding since $d\geq 6$.
Since it is easy to see that 
\[
H^0(\sH_2\times \sH_2,\sD_2)
\simeq H^0(\check{\mP}^{d-3}\times \check{\mP}^{d-3},
\sO(2,2)),
\] 
it holds that $\sD_2$ is the restriction of a unique $(2,2)$-divisor on 
$\check{\mP}^{d-3}\times \check{\mP}^{d-3}$, 
which we denote by $\{\widetilde{\sD}_2=0\}$.
Since $\{\widetilde{\sD}_2=0\}$
is also symmetric, 
we may take the equation $\widetilde{\sD}_2$
so that it is the bi-homogenization of 
an equation $\check{F}_4$ of a quartic in 
$\check{\mP}^{d-3}$ (cf. \cite[\S 1]{DK}).
Moreover the fiber of $\{\widetilde{\sD}_2=0\}$ over a point 
$p\in \check{\mP}^{d-3}$ is defined by the polar $P_p(\check{F}_4)$, 
which we denote by $\widetilde{D}_p$.
For $[q]\in \sH_2$, we denote $\widetilde{D}_{[q]}$ simply by
$\widetilde{D}_{q}$. By construction,
$D_q=\{\widetilde{D}_{q}=0\}\cap (\sH_2\times \sH_2)$. 
We may choose the defining equation $H_q$ of 
the hyperplane of $\mP^{d-3}$ corresponding to $[q]$ such that
$P_{H^2_q}(\check{F}_4)=\widetilde{D}_q$.

From now on, we write $\mP^{d-3}=\mP_* V$, where $V$ is 
the $d-2$-dimensional vector space.
The crucial point in the following assertions is that 
the number of the conics attached to a point of $\widetilde A$
coincides with ${\rm{dim}}_{\mC}S^{2}V$.

   Let $\widetilde{a}$ be a general point of $\widetilde{A}$ and 
    $q_{1},\ldots, q_{n}$ are the conics attached to $\widetilde{a}$.
    By the definition of $\widetilde{D}_{q_i}$ and 
    generality of $\widetilde{a}$, 
we have the following (we use Lemma \ref{noselfintersection}):
\begin{equation}
\label{eq:1'}
\text{$\widetilde{D}_{q_j}([q_i])=0$ $(j\not =i)$
and $\widetilde{D}_{q_i}([q_i])\neq 0$}.
\end{equation}
 
    (\ref{eq:1'}) implies 
    $\widetilde{D}_{q_1},\dots, \widetilde{D}_{q_n}$ 
    are linearly independent. 
Thus by $P_{H_{q_i}^2}(\check{F}_4)=\widetilde{D}_{q_i}$, 
it holds that
the apolarity map 
\begin{eqnarray*}
\ap_{\check{F}_4}\colon S^2 \check{V} & \to & S^2 V\\
G & \mapsto & P_G(\check{F}_4)
\end{eqnarray*}
is an isomorphism.
Moreover, $H_{q_1}^2, \dots,H_{q_n}^2$ are linearly independent.
Thus $\check{F}_4$ is {\em{non-degenerate}} in the sense of Dolgachev.
By \cite[\S 2.3]{doldual},
there exists a unique quartic form $F_4$ 
such that 
$\ap_{F_4}=\ap^{-1}_{\check{F}_4}$.
In particular, it holds
\[
P_{\widetilde{D}_q}(F_4)=H_q^2.
\]
$F_4$ is called 
{\em{the quartic form dual to $\check{F}_4$}}.


To see the relation between the set of conics attached to a
general point of ${\widetilde A}$ and the representation of $F_4$ 
as a sum of powers of linear forms we need to find conditions
which force $n$ conics to be attached to ${\widetilde a}\in {\widetilde A}$.
Next lemma is sufficient for our purposes.
\begin{lem}
\label{lem:GeomA}
Let $q_1,\dots, q_n$ be $n$ distinct conics such that
\begin{enumerate}[$(1)$]
\item
$\widetilde{D}_{q_i}([q_j])=0$ for all $i\not =j$,
\item
all the $\overline{q}_i$ are smooth, 
\item
if three of $\overline{q}_i$ pass through a point $b$,  
then any other $\overline{q}_i$ does not intersect a line through $b$
outside $b$, and
\item
no two of $\overline{q}_i$ intersect at a point of $C\cup \cup_i \beta_i$.
\end{enumerate}
Then the
$q_i$'s are attached to a point of $\tA$.
\end{lem}
\begin{proof}

By the assumption (1),
$\overline{q}_1,\dots, \overline{q}_n$ are 
mutually intersecting multi-secant conics of $C$.
By the assumption (4), it suffices to prove
$\overline{q}_1,\dots, \overline{q}_n$ pass through one point of $B$. \\

\noindent
{\bf{Step 1.}}
Let $b\in B$ be a point such that
five of $\overline{q}_i$, 
say, $\overline{q}_1,\dots, \overline{q}_5$ pass through $b$.
Then all the $\overline{q}_i$ pass through $b$.

By the double projection from $b$,
$\overline{q}_1,\dots,\overline{q}_5$ are mapped to 
points $p_1,\dots,p_5$ on $\mP^2$.
Suppose by contradiction 
that a smooth conic $\overline{q}_j$ does not pass through $b$.
Let $q'_j$, $q''_j$ and $\widetilde{q}_j$ be the strict transforms
of $\overline{q}_j$ on $B_b$, $B'_b$ and $\mP^2$, and set 
$S:=\pi_{2b}^*\widetilde{q}_j$. 
By the assumption (3), $\overline{q}_j$ does not intersect
a line through $b$.
Thus $\widetilde{q}_j$ is a smooth conic 
through $p_1,\dots,p_5$. The conic $\widetilde{q}_j$ is unique
since a conic through five points is unique.
It holds that
$-K_{B'_b}\cdot q''_j=4$ and $S\cdot q''_j=4$, thus
$S\simeq \mF_2$ and $q''_j$ is the negative section.
This implies that $q_j$ is also unique.
By reordering, we may assume that $j=n$.
We have the configuration such that
all the conics pass through $b$ except $q_n$.
Denote by $p_i$ the image of $q_i$ ($i\neq n$).
Then $\widetilde{q}_n$ and $C_b$ intersect at $p_i$.
By $d\geq 6$, it holds $\deg C_b\geq 3$, thus
$\widetilde{q}_n\neq C_b$. 
By the assumption (4), $b\not \in C$.
Therefore $\widetilde{q}_n$ and $C_b$ intersect
at $n-1$ singular points of $C_b$.
Since $\deg C_b\leq d$, it holds
$2(n-1)\leq 2d$, a contradiction.\\
{\bf{Step 2.}}
If four conics $\overline{q}_1, \dots, \overline{q}_4$ 
pass through one point $b$, 
then all the conics pass through $b$.

By contradiction and Step 1,
we may assume that all the conics except 
$\overline{q}_1,\dots,\overline{q}_4$
do not pass through $b$.
Pick up two any conics, say, $\overline{q}_5$ and $\overline{q}_6$, 
not passing through $b$.
Considering the double projection from $b$ as in Step 1.
Denote by $\widetilde{q}_j$ $(j\geq 5)$ 
the image of $\overline{q}_j$ on $\mP^2$. 
By the assumption (3), $\overline{q}_5$ and $\overline{q}_6$
do not intersect a line through $b$,
thus $\widetilde{q}_5$ and $\widetilde{q}_6$ are conics on $\mP^2$.
Therefore $\overline{q}_5 \cap \overline{q}_6$ lies 
on one of $\overline{q}_1,\dots,\overline{q}_4$
since otherwise 
$\widetilde{q}_5$ and $\widetilde{q}_6$ would intersect at five points
and this is a contradiction as in Step 1.
Thus any two conics intersect on $\overline{q}_1,\dots, \overline{q}_4$.
Let $p_i$ be the intersection 
$\overline{q}_i\cap \overline{q}_5$ for $i=1,\dots, 4$.
Then $\overline{q}_j$ ($j\geq 5$) pass through one of $p_i$.
Thus one of $p_i$, say, $p_1$, there
pass through at least $\lceil \frac{(n-5)}{4} \rceil$ conics.
By Step 1,
$\lceil \frac{(n-5)}{4}\rceil \leq 2$ 
(already $\overline{q}_1$ and $\overline{q}_5$ pass through $p_1$).
This implies $d=6$. We exclude this case in Step 3.
Note that if $d=6$, then the four conics $\overline{q}_1$,
$\overline{q}_2$, $\overline{q}_5$, and $\overline{q}_6$
mutually intersect and the all the intersection points are different.
By reordering conics, we assume that 
$\overline{q}_i$ $(1\leq i \leq 4)$ satisfy this property.\\ 
{\bf{Step 3.}} 
We complete the proof.

Assume by contradiction that
$\overline{q}_1,\dots,\overline{q}_n$  do not pass through one point on $B$. 
If $d\geq 7$, then, by Steps 1 and 2,
\begin{equation}
\label{eq:three}
\text{
at most three of $\overline{q}_i$'s pass through any intersection point.}
\end{equation}
Let $m$ be the number of conics in a maximal tree $T$ of $\overline{q}_i$'s 
such that two conics in $T$ pass through any intersection point. 
Note that $T$ is connected since $\overline{q}_i$'s mutually 
intersect. 
The number of the intersection points of $\overline{q}_i$'s contained in 
$T$ is $\frac{m(m-1)}{2}$.

By the maximality of $T$,
a conic not belonging to $T$ passes through one of the
intersection points of conics in $T$.
By (\ref{eq:three}),
no two conics not belonging to $T$ 
pass through one of the  intersection point of conics in $T$.
Hence it holds $\frac{m(m-1)}{2} +m \geq n$.
This implies that $m\geq d-2$ by $n=\frac{(d-1)(d-2)}{2}$.
By reordering, we assume that 
$\overline{q}_1,\dots, \overline{q}_m$ belong to $T$.
If $d=6$, then we take $\overline{q}_1,\dots, \overline{q}_4$
as in the last part of Step 2.
Consider the projection $B \dashrightarrow \mP^3$ from $\overline{q}_1$.
Then $\overline{q}_2,\dots,\overline{q}_m$ are mapped 
to lines $l_2,\dots, l_m$ intersecting
mutually on $\mP^3$ and the intersection points are different.
Thus  $l_2,\dots, l_m$ span a plane, which in turn shows that
$\overline{q}_1,\dots, \overline{q}_m$ span a hyperplane section $H$ on $B$.
Since $C$ intersects $\overline{q}_i$ at two point or more, 
$C$ intersects $H$ at $2m$ points or more by the assumption (4).
But $2m\geq 2(d-2) >d$, $C$ must be contained in $H$,
a contradiction to Proposition \ref{prop:Cd0} (d). 
\end{proof}

We think the next theorem to be of theoretical relevance in itself
and as a first result to understand varieties of sum of powers confined
in a subvariety.

\begin{thm}\label{diretto}
$\Ima \Phi$ is an irreducible component of
$\VSP(F_4,n;\sH_2)$.
\end{thm}

\begin{proof}
Set 
\[
Z:=\{([H_{1}],\ldots, [H_{n}])
    \in \Hilb^{n} \check{\mP}^{d-3} \mid H^{4}_{1}+\ldots
 +H^{4}_{n}=F_4, [H_i]\in \sH_2\}.
\]

    For a general point $\widetilde{a}$ and conics  
    $q_{1},\ldots, q_{n}$ attached to $\widetilde{a}$, 
we have (\ref{eq:1'}).
Conversely, $n$ conics $q_i$ satisfying (\ref{eq:1'}) and
the assumptions (2)--(4) of Lemma \ref{lem:GeomA}
determine a point of $\widetilde{A}$.
Note that the assumptions (2)--(4) of Lemma \ref{lem:GeomA}
are open conditions.
Thus   
we have only to prove that (\ref{eq:1'}) is equivalent to
\begin{equation}
\label{eq:2'}
    \text{$\alpha_1 H^{4}_{q_1}+\ldots +\alpha_n H^{4}_{q_n}=F_4$
    with some nonzero $\alpha_i \in \mC$.} 
\end{equation}

We see that (\ref{eq:2'}) is equivalent to
\begin{equation}
\label{eq:3'}
\text{If $\{G=0\}\subset \check{\mP}^{d-3}$ is any quartic
through
    $[q_1], \cdots, [q_n]$, then $P_{F_4}(G)=0$.}
\end{equation}

Indeed, 
by the apolarity pairing,
$\langle G, H^4_{q_i} \rangle=0\Leftrightarrow  G([q_i])=0$,
thus, the assumption on $G$ is equivalent to
$G\subset \langle H^4_{q_1},\dots,H^4_{q_n}\rangle^{\perp}.$
Therefore
(\ref{eq:2'}) is equivalent to
$\langle H^4_{q_1},\dots,H^4_{q_n}\rangle^{\perp}\subset 
\langle F_4 \rangle^{\perp}.$
Since $F_4$ is non-degenerate,
this is equivalent to (\ref{eq:2'}).    

We show (\ref{eq:1'}) implies (\ref{eq:3'}).
If (\ref{eq:1'}) holds, then
    $\widetilde{D}_{q_i}$ ($i\not =1$)
    generate the space of linear forms passing through $[q_1]$,
    we may write 
$G=Q_2\widetilde{D}_{q_2}+\cdots+Q_n\widetilde{D}_{q_n}$,
    where $Q_i$ are quadratic forms on $\check{\mP}^{d-3}$.
    By $G([q_i])=0$ for $i\not =1$,
    we have $Q_i([q_i]) \widetilde{D}_{q_i}([q_i])=0$.
    $\widetilde{D}_{q_i}([q_i])\not=0$ implies that $Q_i([q_i])=0$. Thus
    $Q_i$ is a linear combination of
    $\widetilde{D}_{q_j}$ $(j\not =i)$.
    Consequently,
    $G$ is a linear combination of
    $\widetilde{D}_{q_i}\widetilde{D}_{q_j}$ $(1\leq i<j\leq n)$.
    Thus $P_{F_4}(G)=0$ follows from 
    that
   \[
   P_{F_4}(\widetilde{D}_{q_i}\widetilde{D}_{q_j})=
   P_{H_{q_i}}(\widetilde{D}_{q_j})=\widetilde{D}_{q_j}([q_i])=0.
   \]
Finally we show (\ref{eq:2'}) implies (\ref{eq:1'}).
By (\ref{eq:2'}), it holds
\[
H_{q_i}^2=P_{\widetilde{D}_{q_i}}(F_4)=
\sum \alpha_j \langle \widetilde{D}_{q_i}, H_{q_j}^4 \rangle H_{q_j}^2.
\]
Since $H_{q_j}^2$ are linearly independent, 
(\ref{eq:1'}) holds.
\end{proof}

\begin{defn}
We say $\Ima \Phi$ is the {\em{main component}} of $\VSP(n,F_4;\sH_2)$.
\end{defn}

The following lemma characterizes the main component of 
$\VSP(n,F_4;\sH_2)$, which will play a crucial role in \ref{subsection:moduli}:

\begin{lem}
\label{lem:main}
Let $(\sH_2^k)^o$ and 
$(\Hilb^k \check{\mP}^{d-3})^o$ $(k\in \mN)$ be the complements of
all the small diagonals of $\sH_2^k$ $(k$ times product of $\sH_2)$ and $\Hilb^k \check{\mP}^{d-3}$
respectively.
Set
\[
\VSP^o(F_4,n;\sH_2):=
\{([H_1],\dots, [H_n])\mid [H_i] \in \sH_2, H_1^m+\cdots+H_n^m=F_4\}.
\]
Let $V^o$ be the inverse image of
$\VSP^o(F_4,n;\sH_2)$ by the natural projection
$(\sH_2^n)^o\to (\Hilb^n \check{\mP}^{d-3})^o$.
Let $(\sH_2^n)^o\to (\sH_2^2)^o$ be the projection
to any of two factors.
Then a component of $V^o$ dominating $\sD_2$
dominates the main component of 
$\VSP(F_4,n;\sH_2)$.
\end{lem}

\begin{proof}
Let $([q_1],[q_2])\in \sD_2\cap (\sH_2^2)^o$ be a general point
and $\{q_i\}$ ($i=1,\dots,n$) any set of mutually conjugate
$n$ conics including $q_1$ and $q_2$.
Since $q_1$ and $q_2$ are general, we may assume that
all the $q_i$ are general. 
By Lemma \ref{lem:GeomA} and Theorem \ref{diretto},
it suffices to prove that
$q_1,\dots,q_n$ satisfies 
the conditions (2)--(4) of Lemma \ref{lem:GeomA}.

(2). Let $\overline{r}_1$ and $\overline{r}_2$ are mutually intersecting
smooth conics on $B$
and $\overline{r}_3$ a line pair on $B$ intersecting
both $\overline{r}_1$ and $\overline{r}_2$.
Since the Hilbert scheme of conics on $B$ is $4$-dimensional,
the pair of $\overline{r}_1$ and $\overline{r}_2$ depends on $7$ parameters.
If we fix $\overline{r}_1$ and $\overline{r}_2$,
then $\overline{r}_3$ depends on $1$ parameter. 
Thus the configuration $\overline{r}_1$, $\overline{r}_2$, $\overline{r}_3$
depends on $8$ parameters.
Fix $\overline{r}_1$, $\overline{r}_2$ and $\overline{r}_3$.
We count the number of parameters of $C_d$ such that
$C_d$ intersects each of $\overline{r}_i$ ($i=1,2,3$) twice.
The number of parameters is 
$h^0((\sO_{\mP^1}(d-1)\oplus \sO_{\mP^1}(d-1))\otimes \sO_{\mP^1}(-6))+6=2d-12+6=2d-6$,
where $+6$ means the sum of the numbers of parameters of two points on 
$\overline{r}_i$ ($i=1,2,3$).
By $2d-6+8=2d+2$,
a general $C_d$ has $2$-dimensional pairs of mutually intersecting bi-secant 
conics which intersect at least one bi-secant line pair of $C_d$.
Thus general pairs of mutually intersecting bi-secant conics of $C_d$, 
which form a $3$-dimensional family, do not intersect
a bi-secant line pair of $C_d$.
  
(3).
Assume by contradiction that  
$\overline{q}_i$, $\overline{q}_j$ and $\overline{q}_k$
pass through a point $b$, and
$\overline{q}_l$ does not pass through $b$ but intersects a line
through $b$. 
Then by the double projection from $b$,
$\overline{q}_l$ is mapped to a line through 
the three singular points of the image of $C_b$ corresponding to  
$\overline{q}_i$, $\overline{q}_j$ and $\overline{q}_k$.
Thus we have only to prove that for a general point of $b$ on $B$,
three double points of the image of $C_b$ do not lie on a line.

Fix a general point $b\in B$.
Let $\overline{r}_1$, $\overline{r}_2$, $\overline{r}_3$ be three conics on $B$ through $b$
such that by the double projection from $b$,
they are mapped to three colinear points on $\mP^2$.
The number of parameters of $C_d$'s intersecting each of $\overline{r}_i$ twice
is 
$h^0((\sO_{\mP^1}(d-1)\oplus \sO_{\mP^1}(d-1))\otimes \sO_{\mP^1}(-6))=2d-12$
since $h^1((\sO_{\mP^1}(d-1)\oplus \sO_{\mP^1}(d-1))\otimes \sO_{\mP^1}(-6))=0$
Note that the number of parameters of
$\overline{r}_1$, $\overline{r}_2$, $\overline{r}_3$
is $5$ since
that of lines in $\mP^2$ is $2$,
and that of three points on a line
is $3$.
Thus the number of parameters of $C_d$'s 
such that its image of the double projection from $b$
has three colinear double points is at most $2d-1$.
Hence a general $C_d$ does not satisfy this property.

(4).
Let $r_1$ and $r_2$ be a general pair of mutually conjugate conics on $A$
such that
$\overline{r}_1$ and $\overline{r}_2$ are smooth, and
$\overline{r}_1$ and $\overline{r}_2$ 
intersect at a point on $C\cup \cup_i \beta_i$.
Such general pairs of conics $r_1$ and $r_2$ 
form a two-dimensional family
since $\dim C\cup \cup_i \beta_i=1$ and
if one point $t$ of $C\cup \cup_i \beta_i$ is fixed,
then such pairs of conics such that
$t\in \overline{r}_1\cap \overline{r}_2$ 
form a one-dimensional family.
For a general pair of $r_1$ and $r_2$,
the number of the sets of $n$ mutually conjugate conics
including $r_1$ and $r_2$ is finite since
$D_{r_1}$ and $D_{r_2}$ has no common component.
Thus $\{q_i\}$ does not contain such a pair by generality
whence $\{q_i\}$ satisfies (4).   
\end{proof}

\subsubsection{Relation with Mukai's result}
\label{subsubsection:v22}
Here we sketch how the argument goes on if $d=5$ and
explain a relation of our result with
Theorem \ref{v22}.

Assume that $d=5$.
Associated to the birational morphism
$\sH_2\rightarrow\check{\mP}^{2}$,
there exists
a non finite birational morphism  
\[
\Phi\colon \widetilde{A}\rightarrow V_{22}:=\VSP(F_4,6)
\subset \Hilb^{6} \check{\mP}^{2},
\] 
which fits into the following diagram: 
\begin{equation*}
\xymatrix{ & & \widetilde{A} \ar[dl]_{\rho} \ar[dr]^{\rho'}
\ar @/^2pc/ [ddrr]^{\Phi} & &\\
& Aе \ar[dl]_{f}  &\dashrightarrow & 
 A' \ar[dr]^{f'} & \\
 B &  &  & & V_{22}, }
\end{equation*}

\noindent 
where 
\begin{itemize}
\item
$V_{22}$ is a smooth prime Fano threefold of genus twelve,
\item
$\rho'$ is the blow-down of the three $\rho$-exceptional divisors 
$E_i$ ($i=1,2,3$)
over the strict transform $\beta'_i$ in the other direction.
In other words, $A\dashrightarrow A'$ is the flops of $\beta'_1$, $\beta'_2$
and $\beta'_3$ (cf. Lemma \ref{sullebisecanti}), and
\item
the morphism $f'$ 
contracts the strict transform of 
the unique hyperplane section $S$ containing $C$
(see Proposition \ref{prop:Cd0} (d))
to a general line on $V_{22}$.
\end{itemize}
The rational map $V_{22}\dashrightarrow B$ is the famous double projection
of $V_{22}$ from a general line $m$ 
first discovered by Iskovskih (see \cite{I2}).

We explain how $f'$ and $\rho'$ are interpreted in our context. 
As we remarked after the proof of Theorem \ref{thm:H_2},
the morphism $\sH_2\to \check{\mP}^2$ defined by $|D_l|$ contracts
three curves $D_{e_i}$ which parameterize conics
intersecting $\beta'_i$.
By noting $S$ is covered by the images of such conics,
this corresponds to that 
the morphism $f'$ contracts the strict transform of $S$.

We can see that any conic on $A$ except one belonging to $D_{e_i}$
corresponds to that on $V_{22}$ in the usual sense, and
the component of Hilbert scheme of $V_{22}$ parameterizing
conics is naturally isomorphic to $\check{\mP}^2$.
The three conics on $V_{22}$ corresponding to the images of $D_{e_i}$
are $\beta''_i\cup m$, where $\beta''_i$ are the images of the flopped curve
corresponding to $\beta'_i$.

Let $a\in E_i$. Then six conics on $A$ attached to $a$
are $\xi_{ij1}$ ($j=1,2$), a conic $q_a$ from $D_{e_i}$ and
three conics from $e_i$ 
(see the remark at the end of \ref{subsubsection:special}).
Moreover, if $a$ moves in a fiber $\gamma$ 
of the other projection $E_i\to \mP^1$,
then only the conic $q_a$ from $D_{e_i}$ varies.
By the contraction $\sH_2\to \check{\mP}^2$, 
there is no difference among points on $\gamma$.
This is the meaning of the contraction $\rho'$ of $E_i$
in the other direction.

Finally we remark that $\sH_1$ is also naturally isomorphic to    
the component of Hilbert scheme of $V_{22}$ parameterizing lines.



\section{The existence of the Scorza quartic} 
In this section we will use the geometries of $\sH_{1}$ and $\sH_{2}$
to give an affirmative answer to the conjecture of Dolgachev and Kanev
 \cite[Introduction p. 218]{DK} (see Theorem \ref{thm:DK}).

\subsection{Theta-correspondence on $\mathcal{H}_1\times \sH_1$}
 \label{subsection:theta}
 ~

In this subsection, we regard $\sH_1$ as the component of the Hilbert scheme
of $A$ parameterizing lines on $A$.
   
We will define a non-effective theta characteristic on $\sH_1$ by
investigating the following set:
\[
  I:=
  \{([l_1], [l_2])\in \mathcal{H}_1\times \mathcal{H}_1 \mid
\text{$l_1$ and $l_2$ intersect}\}.
\]

We need a more precise and technical definition of $I$. First we
reconsider the desingularization morphism
$\pi_{|\sH_1}\colon\sH_{1}\rightarrow M\subset\mP^{2}$; see Corollary \ref{orabisecanti}. 

\begin{lem}
\label{lem:general}
$h^{0}(\sH_1,(\pi_{|\sH_1})^*\sO_M(1))=3$.
\end{lem}

\begin{proof}
Let $h\colon S \to \sH_1^{B}\simeq \mP^2$ be the blow-up of
$\sH_1^{B}$ at
the $s=\frac{(d-2)(d-3)}{2}$ nodes of $M$. 
Then $\sH_1\sim d \lambda -2\sum_{i=1}^{s}\varepsilon_i$,
where $\lambda$ is the pull-back of a general line and $\varepsilon_i$ 
are exceptional curves.
By the exact sequence
\[
0\to \sO_S(\lambda-\sH_1) \to \sO_S(\lambda)  \to 
\sO_{\sH_1}((\pi_{|\sH_1})^*\sO_M(1))\to 0
\]
together with
$h^0(\sO_S(\lambda))=3$ and 
$h^0(\sO_S(\lambda-\sH_1))=h^1(\sO_S(\lambda))=0$,
we see that
$h^{0}(\sH_1,(\pi_{|\sH_1})^*\sO_M(1))=3$ is equivalent to
$h^1(\sO_S(\lambda-\sH_1))=0$.
By the Riemann-Roch theorem, we have
$\chi(\sO_S(\lambda-\sH_1))=0$. 
Thus by $h^0(\sO_S(\lambda-\sH_1))=0$,
$h^1(\sO_S(\lambda-\sH_1))=0$ is equivalent to
$h^2(\sO_S(\lambda-\sH_1))=0$.
By the Serre duality,
$h^2(\sO_S(\lambda-\sH_1))=h^0(\sO_S((d-4)\lambda-\sum_{i=1}^{s}\varepsilon_i)$.
Thus
we have only to prove that
there exists no plane curve of degree $d-4$ through $s$ nodes of $M$.
We prove this fact by using the inductive construction of $C=C_d$.
In case $d=2$, the assertion is obvious.
From now on in the proof, we put the suffix $d$ to the object
depending on $d$.
For example, $s_d:=\frac{(d-2)(d-3)}{2}$.
Assuming  
$h^0(\sO_{S_d}((d-4)\lambda_d-\sum_{i=1}^{s_d}\varepsilon_{i,d})=0$,
we prove 
$h^0(\sO_{S_{d+1}}((d-3)\lambda_{d+1}-\sum_{i=1}^{s_{d+1}}
\varepsilon_{i,d+1})=0$.

Recall that we constructed $C_{d+1}$ by the smoothing of the union of $C_d$ 
and a general uni-secant line $\overline l$ of $C_d$.
By a standard degeneration argument,
we have only to prove that
there exists no plane curve of degree $d-3$ through $s_{d+1}$ nodes of 
$M_d\cup M(\overline{l})$, where $s_d$ of $s_{d+1}$ nodes are those of $C_d$ and
the remaining $s_{d+1}-s_d=d-2$ nodes are $M_d\cap M(\overline{l})$ except the two points
corresponding
to the two other lines $\overline l'$, 
$\overline l''$ through $C_d\cap \overline{l}$. 
Assume that 
there exists a plane curve $G$ of degree $d-3$ through $s_{d+1}$
nodes of $M_d\cup M(\overline{l})$. Then 
$G\cap M(\overline{l})$ contains at least $d-2$ points. 
Since $\deg G=d-3$, this implies $M(\overline{l})\subset G$. Thus 
there exists a plane curve of degree $d-4$ through $s_d$ nodes of
$M_d$, a contradiction.  
\end{proof}

We denote by $\delta$ the $g^1_3$ on $\sH_1$ which defines 
$\varphi_{|\sH_1}\colon \sH_1\rightarrow C$.
Let $l$, $l'$ and $l''$ be three lines on $A$ such that 
$[l]+[l']+[l'']\sim \delta$. Then $\overline{l}$, $\overline{l}'$
and $\overline{l}''$
are lines through one point of $C$.
Set 
\[
\theta:=(\pi_{|\sH_1})^*\sO_M(1)-\delta.
\]

 Let $l$ be any line on $A$ and
    $l', l''$ lines 
    such that $[l]+[l']+[l'']\sim \delta$.  
    By $\theta+[l]=\pi_{|\sH_1}^*\sO_M(1)-[l']-[l'']$ and
    Lemma \ref{lem:general}, we have 
    $h^{0}(\sH_1, \sO_{\sH_1}(\theta +[l]))=1$.
Let $p_i\colon \sH_1\times \sH_1 \to \sH_1$ ($i=1,2$) be the two
projections and 
    $\Delta$ the diagonal of $\sH_1\times \sH_1$.
    Set $\sL:=\sO_{\sH_1\times \sH_1}({p_2}^*\theta+\Delta)$. 
    By $h^{0}(\sH_1, \sO_{\sH_1}(\theta +[l]))=1$ for any $[l]\in
\sH_1$,
    we see that $p_{1*}\sL$ is an invertible sheaf.
    Define an ideal sheaf $\sI$ by ${p_1}^*p_{1*}\sL=\sL\otimes \sI$. 
    $\sI$ is an invertible sheaf and let $I$ be the divisor defined
by $\sI$. 
Then we can extract the following definition:
\begin{defn}\label{mukaiuno}
$I$ is called the {\it{theta-correspondence}}.
We will denote by $I([l])$ the fiber of $I\to \sH_1$ over $[l]$.
\end{defn}

The following result is a generalization of Mukai's result
\cite[\S 4, Theorem]{Mukai12} in our setting:

\begin{prop}\label{prop:sopratheta}
$\theta$ is a non-effective theta characteristic.
\end{prop}

\begin{proof}
    By invoking  \cite[Lemma 7.2.1]{DK} and the definition of $I$,
    it suffices to prove the following:
\begin{enumerate}[(a)]
\item
    $h^{0}(\sH_1, \sO_{\sH_1}(\theta +[l]))=1$ for any $[l]\in \sH_1$,
\item
    $I$ is reduced,
\item
$I$ is disjoint from the diagonal,
\item
$I$ is symmetric, and
\item
$I$ is a $(g(\sH_1), g(\sH_1))$-correspondence.
\end{enumerate}    

 Let $l$ be any line on $A$ and
    $l', l''$ lines 
    such that $[l]+[l']+[l'']\sim \delta$.  

We have proved (a) already.

Noting that 
 the line in $\mP^{2}$ joining $[\overline{l}']$ and $[\overline{l}'']$ parameterizes the lines
on $B$ intersecting $\overline{l}$, we see that 
the fiber of $I\to \sH_1$ over a general $[l]$ is reduced. 
Hence $I$ is reduced.

We prove (c).
It is equivalent to show that
the support of $I([l])$ does not contain $[l]$.
By definition $\theta +[l]=
\pi_{|\sH_1}^*\sO_M(1)-[l']-[l'']$.
If $\overline{l}$ is a uni-secant and is not special, then
$M(\overline{l})$ does not contain $[\overline{l}]$, thus we are done.
If $\overline{l}$ is special,
then, by Propositions \ref{prop:FN} (4) and \ref{prop:Cd} (2), we are done. If $\overline{l}$
is a bi-secant then by Proposition \ref{prop:Cd1} (4), we are done.

We prove (d).
Let $m$ be a line on $A$ such that
$[m]$ is contained in the support of
$I([l])$. It suffices to prove that for a general $l$,
$[l]$ is contained in the support of
$I([m])$.
For a general $l$, we may assume that 
$m\not =l'$ or $l''$. 
Then it is easy to verify this fact.

Finally we prove (e).
Since $I$ is symmetric and $\deg (\theta+[l])=d-2=g(\sH_1)$,
the divisor is a $(g(\sH_1),g(\sH_1))$-correspondence.    
\end{proof}



\subsection{Duality between $\sH_1$ and $\sH_2$}
\label{subsection:duality}
~

Denote by $\mP^{d-3}$ the projective space dual to $\check{\mP}^{d-3}$.  
The family 
\[
\xymatrix
{
\sD_1\ar[r]\ar[d] & \sH_2\times \sH_1\ar[dl]\\
\sH_1 &
}
\]
induces the morphism 
\begin{eqnarray*}
\sH_1& \to & \mP^{d-3}\\
{[l]} & \mapsto & {[D_l]}.
\end{eqnarray*} 
by the universal property of the Hilbert scheme.
Since $D_l\not =D_{l'}$ for $l\not =l'$, 
$\sH_1\to \mP^{d-3}$ is injective.

Consider
the projection $\sD_{1}\rightarrow\sH_{2}$
and denote by $\widetilde{H}_q$ the fiber over $[q]$.
Since 
$\sD_1$ is a Cartier divisor in a smooth $3$-fold $\sH_1\times \sH_2$ 
then $\sD_{1}$ is Cohen-Macaulay. 
Since no conic on $A$ intersects infinitely many lines on $A$,
$\sD_1\to \sH_2$ is finite. Then $\sD_{1}\rightarrow\sH_{2}$ is flat
since $\sH_{2}$ is smooth.
Note that for a general $q$, $\widetilde{H}_q$ parameterizes
all the lines intersecting $q$.
By considering the morphism
$\pi_{\vert\sH_{1}}\colon
\sH_{1}\rightarrow M\subset\mP^{2}$,
it is easy to see that for a general conic $q$,
$\widetilde{H}_q\in |\pi^*\sO_M(2)-2\delta|$, namely,  
$\widetilde{H}_q\sim 2\theta\sim K_{\sH_1}$.
By the flatness of $\sD_1\to \sH_2$,
it holds $\widetilde{H}_q$ for any $q$.

Recall that we denote by $\{H_q=0\}$ the hyperplane in $\mP^{d-3}$
corresponding to $[q]\in \check{\mP}^{d-3}$.
Note that, for $[l]\in \sH_1$ and $[q]\in \sH_2$, 
$[l]\in \{H_q=0\}$ 
if and only if 
$D_l ([q])=0$.
Thus $\widetilde{H}_q=\{H_q=0\}$. 
Consequently, the injection $\sH_1\to \mP^{d-3}$ is the canonical embedding
$\Phi_{|K_{\sH_1}|}\colon \sH_1\to \mP^{d-3}$
by $\widetilde{H}_q\sim K_{\sH_1}$.

In case $d=5$, a similar construction gives
the duality of the canonical embedding $\sH_1 \subset \mP^2$ and
the birational morphism $\sH_2\to\check{\mP}^2$.

\subsection{Discriminant locus}
 \label{subsection:discr1}
 ~

We follow \cite[7.1.4 p.279]{DK}. 
Let $\Gamma \subset \mP^{g-1}$ be a canonical curve of genus $g$
  and $\theta'$ a non-effective even theta characteristic on $\Gamma$.
By the Riemann-Roch theorem,
it holds that $h^0(\theta'+x)=1$ for a point $x\in \Gamma$.
Let 
$$I:=
\{(x,y)\mid y\ \text{is in the support of the unique member of}\ |\theta'+x|\}
\subset \Gamma\times \Gamma.$$
We call this the {\em{theta-correspondence}},
which is consistent with Definition \ref{mukaiuno}.
We denote by $I(x)$ the fiber of $I\to \Gamma$ over $x$
and call it the {\em{theta-polyhedron}} attached to $x$.
In other words, $I(x)$ is the unique member of
$|\theta'+x|$ as a divisor.

Since the linear hull $\langle I(x)-y\rangle$
  is a hyperplane of $\mP^{g-1}$, then we can define a morphism
  $\pi_{\theta'}\colon I\rightarrow |K_{\Gamma}|=\check{\mP}^{g-1}$ 
as a composition of
  the embedding $I\hookrightarrow \Theta_{\Gamma}$ and the Gauss map
  $\gamma\colon\Theta^{{\rm {ns}}}_{\Gamma}\to\check{\mP}^{g-1}$,

 where
  $\Theta_{\Gamma}\subset J(\Gamma)$ is the theta divisor and 
$\Theta^{{\rm{ns}}}_{\Gamma}$ is the nonsingular locus of $\Theta_{\Gamma}$.

\begin{defn}
 \label{defn:Gamma}  
 The image $\Gamma(\theta')$ of the above morphism 
 $\pi_{\theta'}\colon I\to \check{\mP}^{g-1}$ 
  is called the {\it{discriminant
 locus}} of $(\Gamma,\theta')$. 

Set-theoretically $\pi_{\theta'}$ is the map $(x,y)\mapsto  
\langle I(x)-y\rangle$.
The hyperplane $\langle I(x)-y\rangle$ is called the {\it{face}} 
  of $I(x)$ opposed to $y$.

\end{defn}

From now on in the section 4,
we assume that $d\geq 6$ for the pair
$(\sH_1,\theta)$ and we consider $\sH_2\subset \check{\mP}^{d-3}$. 

For the pair $(\sH_1,\theta)$,
we can interpret $\Gamma(\theta)$
by the geometry of lines and conics on $A$ 
as follows:

\begin{prop}
 \label{prop:discr}
For the pair $(\mathcal{H}_1,\theta)$, 
the discriminant locus $\Gamma(\theta)$ is contained in $\sH_{2}$, and
the generic point of the curve $\Gamma(\theta)$ 
parameterizes line pairs.  
\end{prop}       

\begin{proof}
  
  Take a general point $([l_1],[l_2])\in I$, equivalently,
  take two general intersecting lines $l_1$ and $l_2$.
  $l_1\cup l_2$ is a line pair and
  the lines corresponding to the points of $I([l_1])-[l_2]$ 
  are lines intersecting $l_1$ except $l_2$.
  Thus by \ref{subsection:duality},
  the point in $\check{\mP}^{d-3}$ corresponding to the hyperplane  
  $\langle I([l_1])-[l_2]\rangle$ is nothing but
  $[l_1\cup l_2]\in \mathcal{H}_2$. This implies the assertion.
\end{proof}

\begin{prop} 
\label{prop:gamma}
The curve $\Gamma(\theta)$ belongs to 
the linear system $|3(d-2)h-4\sum_{i=1}^{s}e_i|$ on $\sH_2$.

\end{prop}

\begin{proof} We can write:
\[
\Gamma(\theta) \sim ah-\sum m_i e_i,
\]
where $a\in \mZ$ and $m_i\in \mZ$.
For a general $b\in C$, $L_b$ intersects $\Gamma(\theta)$ simply.
Thus $a$ is the number of line pairs whose images on $B$ pass through $b$.
There exists three lines $l_1$, $l_2$ and $l_3$ through $b$.
It suffices to count the number of reducible conics on $B$
having one of $l_i$ as a component except $l_1\cup l_2$,
$l_2\cup l_3$ and $l_3\cup l_1$. Thus $a=3(d-2)$.

We will count the number
of line pairs belonging to $e_i$.
Each of such line pairs is 
of the form $l_{ij;k} \cup l_{ij}$,
where $l_{ij;k}$ ($k=1,2$) is 
the strict transform of the line through $p_{ij}$
distinct from $\beta_i$.
Thus the number of such pairs is four and $m_i\geq 4$.

Finally
we will count the number
of line pairs intersecting
a general line $l$. By Corollary \ref{noncontenere}, $D_{l}$ does not
contain any line pair $l\cup l'$. 
Since the number of lines on $A$ intersecting a fix line
on $A$ is $d-2$, 
we see that $D_{l}\cdot\Gamma(\theta)\geq (d-2)(d-3)$. Then
\[
(d-2)(d-3)\leq \Gamma(\theta) \cdot D_l=
(d-3)a-\sum_{i=1}^{s}ее m_i.
\]
where $s=\frac{(d-2)(d-3)}{2}$. This implies that $m_i=4$.

\end{proof}

\begin{cor}
\label{cor:DK}
For $(\sH_1,\theta)$,
it holds that
$\deg \Gamma(\theta)=g(g-1)$ and
$p_a(\Gamma(\theta))=\frac{3}{2}g(g-1)+1$.
\end{cor}
\begin{proof}
The invariants of $\Gamma(\theta)$ are easily calculated.
\end{proof}


\subsection{Definition of the Scorza quartic}
\label{subsection:Scorza}
~

By Definition \ref{defn:Gamma}, we have the following diagram:

\begin{equation}
\label{iniezionescorza}
\xymatrix{ & I\subset \Gamma\times\Gamma \ar[dl]_{{ \pi_{\theta'}}}
\ar[dr]^{p}\\
 \Gamma(\theta')\subset\check\mP^{g-1} &  & \Gamma \subset\mP^{g-1}. }
\end{equation}
\noindent
 We can define:
\[
  \overline{D}_H:=\pi_{\theta' *}p^*(H\cap \Gamma),
\]
\noindent
where $H$ is an hyperplane of $\mP^{g-1}$.
It is easy to see: 
\[
  \deg \overline{D}_H=2g(g-1).
\]      

Let $S^{m}\check{V}$ the space of $m$-th symmetric forms on the vector
space $V$. Note that an element of $S^m\check{V}$ defines 
a hypersurface of degree $m$ in $\mP_*V$.
Let $F\in S^{2k} \check{V}$ be a non-degenerate homogeneous form of degree
$2k$ and $\check{F}\in S^{2k} V$ the dual homogeneous form
to $F$ defined as in \cite[\S 2.3]{doldual}.
Following \cite[4.1]{doldual},
we define the variety of {\it{the conjugate pairs}}
\[
\mathrm{CP}(F):=
\{([H_1],[H_2])\in \mP_*\check V\times \mP_*\check V \mid
\langle H_1^k, P_{H_2^k}(\check{F})\rangle=0\},
\]
where $\langle \ , \ \rangle$ is the polarity pairing.
Let 
\[
 \Delta:=\mathrm{CP}(F)\cap (\text{the diagonal of}\ 
 \mP_*\check V\times \mP_*\check V).
 \]   
 Since the diagonal of
 $\mP_*\check V\times \mP_*\check V$ is isomorphic to 
 $\mP_*\check V$ then $\Delta\simeq \{\check{F}=0\}$.

Set $D'_{H}:=P_{H^k}(\check{F})$
for a hyperplane $H\subset \mP_*V$.
Then we can write:
\[
\mathrm{CP}(F)=
\{([H_1],[H_2])\in \mP_*\check V\times \mP_*\check V \mid
D'_{H_2}([H_1])=0\}.
\]

\begin{defn}
A non-degenerate quartic $\{F'_4=0\}$ is called 
the {\em{Scorza quartic}} for $(\Gamma,\theta')$
if
$\{D'_{H}=0\}\cap \Gamma(\theta)=\overline{D}_{H}$ 
for a hyperplane $\{H=0\}$
such that $\Gamma\cap \{H=0\}$ is reduced,
where $D'_H$ is defined as above for $F'_4$. 
\end{defn}
  

\subsection{Dolgachev-Kanev's conjecture
on the existence of the Scorza quartic}\-
\label{subsection:conj}
We show that the following properties
hold for general pairs of canonical curves $\Gamma$ and
even theta characteristics $\theta'$ 
as Dolgachev and Kanev conjectured.

\begin{enumerate}[({A}1)]
 \item
  The number of theta-polyhedrons having a general face in common
is two.
Equivalently, the degree of the map $I\to\Gamma(\theta')$ is two,
 \item
  $\Gamma(\theta')$ is not contained in a quadric, and
 \item
  $I$ is reduced. 
 \end{enumerate}

By \cite[Theorem 9.3.1]{DK},
these three conditions are sufficient for
the existence of 
the Scorza quartic for the pair $(\Gamma, \theta')$.

First we show that for our trigonal curve $\sH_1$ and 
the even theta characteristic $\theta$ defined by intersecting lines 
the above conditions hold.

\begin{lem}
\label{lem:DK}
$(\sH_1,\theta)$ satisfies $(\mathrm{A}1)$--$(\mathrm{A}3)$.
\end{lem}

\begin{proof}
(A1) This condition means that for general lines $l$ and $l'$ on $A$ 
such that $([l],[l'])\in I$ the face $\langle I([l])-[l']\rangle$ belongs only
to $I([l])$ and to $I([l'])$.

By contradiction assume that there exists
a line $m$ on $A$ such that $m\neq l$, $m\neq l'$ and 
  $\langle I([l])-[l']\rangle$ is a face of $I([m])$.
  Then  some $d-3$ points of $I([m])$ lie on the hyperplane
  $\langle I([l])-[l']\rangle$, equivalently,
  $m$ intersects $d-3$ lines on $A$ 
  corresponding to the points of 
  $I([l])\cup I([l'])$ except $l$ and $l'$.
By $d\geq 6$,
  it holds that, for $l$ or $l'$, say, $l$,  there exist
  two lines intersecting both $l$ and $m$.
  Consider the projection $B\dashrightarrow Q$ from
  $f(l)={\overline l_{}}$: 
\begin{equation*}
\label{eq:linebis}
\xymatrix{ & B_{\overline l}е \ar[dl]_{\pi_{1}}
\ar[dr]^{\pi_{2}}\\
 B &  & Q.}
\end{equation*}
We use the notation of Proposition \ref{prop:proj} (2).
Now notice that 
by generality of $l$, ${\overline l}\not =\overline{m}:=f(m)$ is
equivalent to have $l\not =m$. 
Since there exist two lines intersecting both ${\overline l}$ and 
$\overline{m}$, we have ${\overline l}\cap \overline{m}=\emptyset$.
Thus the strict transform $\overline{m}'$ of $\overline{m}$ on $Q$ is a line. 
Since there exist two lines intersecting both ${\overline l}$ and 
$\overline{m}$, 
$\overline{m}'$ intersects the image $E'_{\overline l}$ of $E_{\overline l}$ at two points.
Since $E'_{\overline l}$ is a hyperplane section on $Q$,
this implies that $\overline{m}'\subset E'_{\overline l}$, a contradiction.
\\
(A2) This condition is satisfied by
Theorem \ref{thm:H_2} (4) and Proposition \ref{prop:gamma}.\\
(A3) We prove this in the proof of Proposition \ref{prop:sopratheta}.
\end{proof}

Let 
\[
\Gamma'(\theta):=I/(\tau),
\]
where $\tau$ is the involution on $I$ induced by
that of $\Gamma\times \Gamma$ permuting the factors. 
Note that $I\to \Gamma(\theta)$ factor through $\Gamma'(\theta)$.

\begin{cor}
\label{cor:DK2}
For $(\sH_1,\theta)$,
it holds
$\Gamma(\theta)'\simeq \Gamma(\theta)$.
\end{cor}

\begin{proof}
By Lemma \ref{lem:DK}, (A1) holds for $(\sH_1,\theta)$.  
Thus, by
\cite[Corollary 7.1.7]{DK},
we have 
$p_a(\Gamma(\theta)')=
\frac{3}{2}g(g-1)+1$.
Thus, by Corollary \ref{cor:DK},
$p_a(\Gamma(\theta)')= p_a(\Gamma(\theta))$.
By (A1) again,
the natural morphism
$\Gamma(\theta)'\to \Gamma(\theta)$ is birational.
Therefore
it holds
$\Gamma(\theta)'\simeq \Gamma(\theta)$.
\end{proof}

By a moduli argument we prove the conjecture for a general pair 
$(\Gamma,{\theta}')$.
\begin{thm}
\label{thm:DK}
A general spin curve   
satisfies the conditions $(\textup{A}1)$--$(\textup{A}3)$.
In particular, the Scorza quartic exists for a general spin curve.
\end{thm}

\begin{proof}
Let  
$\sM$ be 
the moduli space of couples $(\Gamma,{\theta}')$, where
$\Gamma$ is a curve of genus $g$ and $\theta'$ is a theta
characteristic 
such that $h^{0}(\Gamma,\theta')=0$. 
Classically, $\sM$ is known to be irreducible
(see \cite{cornalba}). 
Let $U$ be a suitable finite cover of an open neighborhood of a couple
$(\sH_1,\theta)$
such that
there exists the family
$\sG\rightarrow U$ of 
pairs of canonical curves 
and non-effective theta characteristics.
Denote by $(\Gamma_u,\theta_u)$ the fiber of $\sG\to U$ over $u\in U$. 
By Lemma \ref{lem:DK}, 
$(\sH_1,\theta)$ satisfies (A1)--(A3).
Since the conditions (A1) and (A3) are open conditions,
these are true on $U$. 
Thus we have only 
to prove that the condition (A2) is still true on $U$.
Let $\sJ\rightarrow U$ be the family of Jacobians and
$\Theta\rightarrow U$ the corresponding family 
of theta divisors. By \cite[p.279-282]{DK}, the family
$\sI$ of theta-correspondences embeds into $\Theta$, and 
by the family of Gauss maps
$\Theta\rightarrow \check{\mP}^{g-1}\times U$, 
we can construct the family
$\widetilde{\sG}\rightarrow U$ whose fiber
$\widetilde{\sG}_{u}\subset \check{\mP}^{g-1}$ 
is the discriminant $\Gamma(\theta_u)$.
By Corollary \ref{cor:DK2},
it holds $\Gamma(\theta)'\simeq\Gamma(\theta)$ for $(\sH_1,\theta)$.
Thus we have also $\Gamma(\theta_u)'\simeq \Gamma(\theta_u)$ for $u\in U$.
By \cite[Corollary 7.1.7]{DK},
we see that
$p_a(\Gamma(\theta_u))$ and $\deg \Gamma(\theta_u)$ are constant for $u\in U$. 
Thus $\widetilde{\sG}\rightarrow U$ is a flat family
since the Hilbert polynomials are constant.
Since no quadric contains $\Gamma(\theta)$ for $(\sH_1,\theta)$,
neither does $\Gamma(\theta_u)$ for $u\in U$ 
by the upper semi-continuity theorem.
\end{proof} 

\begin{rem}
Let $(\Gamma,\theta')$ be a general pair of 
a canonical curve $\Gamma$ 
and a non-effective theta characteristic $\theta'$.
In the proof of Theorem \ref{thm:DK},
we prove that $\Gamma(\theta')'\simeq \Gamma(\theta')$.
\end{rem}

\subsection{$F_4$ is the Scorza quartic for $(\sH_1,\theta)$}
\label{subsection:coincide}
~

Note that, for $F_4$, it holds $\sD_2=\CP(F_4)_{|\sH_2\times \sH_2}$ and 
$\widetilde{D}_q=P_{H_q^2}(\check{F}_4)$.

Let $\{F'_4=0\}$ be the Scorza quartic for $(\sH_1,\theta)$.
Let $D'_H$ be defined as in \ref{subsection:Scorza} for $F'_4$.
We simply denote $D'_{H_q}$ by $D'_q$.







\begin{prop}
\label{prop:quadrics} 
$F_4$ is the Scorza quartic for $(\sH_1,\theta)$.
\end{prop}

\begin{proof} 
The assertion is equivalent to the following equality: 
\[
\mathrm{CP}(F'_4)_{\vert
\sH_{2}\times\sH_{2}}=\mathrm{CP}(F_4)_{\vert
\sH_{2}\times\sH_{2}}.
\]
Since both sides 
have the structures of quadric section bundles over $\sH_2$,
it suffices to prove $\{D'_q=0\}=\{\widetilde{D}_q=0\}$ for a generic $q$.
Since there is no quadric containing $\Gamma(\theta)$ by Lemma
\ref{lem:DK}, 
it is sufficient to show that 
$\{D'_q=0\} \cap \Gamma(\theta)=\{\widetilde{D}_q=0\} \cap \Gamma(\theta)$.
The set $\{\widetilde{D}_q=0\} \cap \Gamma(\theta)$ consists of points 
corresponding to the line pairs intersecting $q$. 
On the other hand, $\{D'_q=0\}\cap \Gamma(\theta)=\overline{D}_{H_q}$
by the definition of the Scorza quartic 
since $H_q$ is reduced for a general $q$.
By the definition of $\overline{D}_{H_q}$,
we can easily show 
the set $\overline{D}_{H_q}$ also consists of points 
corresponding to the line pairs
intersecting $q$. 
Thus $\{D'_q=0\} \cap \Gamma(\theta)=\{\widetilde{D}_q=0\} \cap \Gamma(\theta)$
as desired.
\end{proof}

\subsection{Moduli space of trigonal spin curves}
\label{subsection:moduli}
~

As in Mukai's case we can reconstruct the threefold 
$\widetilde{A}$, that is the couple $(B,C)$, via the curve $\sH_1$ and
a non-effective theta characteristic $\theta$ on it.
\begin{prop}
\label{prop:recovery}
$\widetilde{A}$ is recovered from $(\sH_1,\theta)$.
\end{prop}

\begin{proof}
From $(\sH_1,\theta)$, we can define
$\Gamma(\theta)$ as in Definition \ref{defn:Gamma}
and $F_4$ by Proposition \ref{prop:quadrics}. 
By Theorem \ref{thm:H_2} and Proposition \ref{prop:gamma}, 
$\sH_2$ is recovered from $\Gamma(\theta)$ as
the intersection of cubics containing $\Gamma(\theta)$. 
By Theorem \ref{diretto} and Lemma \ref{lem:main}, 
$\tA$ is recovered from $F_4$ and $\sH_2$.
\end{proof}

For the next result, we denote $\tA$ by $\tA_d$.

Recall that 
we denote by $\sH_d^{B}$ the union of components
of the Hilbert scheme of $B$ whose general points
parameterize smooth rational curves of degree $d$
obtained inductively as in Proposition \ref{prop:Cd0}.
By the remark after the proof of Proposition \ref{prop:Cdbis},
$\sH_d^B$ is irreducible.

We identify $\sH_d^B$ with the moduli space of $\tA_d$, 
which we denote by $\sM_d$.
Let $\sM'_g$ and $\widetilde{\sM}'_g$ be
the moduli space of trigonal curves of genus $g$ and
the moduli space of pairs of trigonal curves of genus $g$ and 
even theta characteristics, respectively.
We can define the rational map 
$\pi_{\sM}\colon \sM_d\dashrightarrow \widetilde{\sM}'_{d-2}$ by
setting $\tA_d\mapsto (\sH_1,\theta)$.
\begin{cor}
\label{cor:gen}
$\pi_{\sM}$ is birational.
Moreover, $\Ima \pi_{\sM}$ is an irreducible component of 
$\widetilde{\sM}'_{d-2}$ dominating $\sM'_{d-2}$.
In particular a general $\sH_1$ is a general trigonal curve of genus $d-2$.
\end{cor}

\begin{proof}
The first assertion follows from Proposition \ref{prop:recovery}.

Since $\dim \sH_d^B=2d$ and   
$\dim \Aut (B,C_d)\leq 3$,
we see that $\dim \sM_d\geq 2d-3$.
On the other hand, $\dim \sM'_{d-2}=2d-3$
and a smooth curve has only a finite number of theta-characteristics.
Thus the latter part follows from the first. 
\end{proof}

Combining Theorem \ref{diretto}, Proposition \ref{prop:quadrics} and
Corollary \ref{cor:gen}, we obtain:

\begin{cor}
Let $F_4$ be the Scorza quartic for a general trigonal spin curve of genus 
$d-2\, (d\geq 6)$ 
and $\sH_2$ the intersection of cubics containing the discriminant locus
of the trigonal spin curve. Set $n:=\frac{(d-1)(d-2)}{2}$.
Then
the normalization of the main component of $\VSP(F_4, n;\sH_2)$ is
isomorphic to the blow-up of quintic del Pezzo threefold $B$
along a general smooth rational curve of degree $d$
and then the strict transforms of its bi-secant lines on $B$.
\end{cor}

By our study we see some other problems which are of a certain interest:

\begin{prob}
\label{prob:irred}
\begin{enumerate}[(1)]
\item
Is the Hilbert
scheme of curves of $B$ whose generic point corresponds to a
smooth rational curve of degree $d$ irreducible, 
namely, is $\sH_d^B$ the unique irreducible component ?
\item
Is $\widetilde{\sM}'_g$ irreducible ?
\end{enumerate}
\end{prob}

If $d=5$, (2) is true by \cite[Lemma 7.7.1]{DK}.
We show that (1) is true also for $d\leq 6$.
(Probably if $d\leq 5$, then it is known. Our contribution is 
for $d=6$). 


\begin{prop}
\label{prop:Cdbis} 
If $d\leq 6$, then the answer to Problem $\ref{prob:irred}$ $(1)$ is 
affirmative.
\end{prop}

\begin{proof} 
For a smooth projective variety $X$ in some projective space,
let $\sC^{0}_{d}(X)$ be the components of the Hilbert scheme of $X$ 
whose general 
points parameterize smooth rational curves of degree $d$.
By \cite{Per}, $\sC^{0}_{d}(G(a,b))$ is irreducible,
where $G(a,b)$ is the Grassmannian parameterizing
$a$-dimensional sub-vector spaces in a fixed
$b$-dimensional vector space. 
The claim is that $\sC^{0}_{d}(\mP^{6}\cap G(2,5))$ is
irreducible, where $\mP^{6}\subset \mP^9$ is transversal to $G(2,5)$.
The claim is true for $d=1$ since $\sH^B_1\simeq \mP^2$.

Let $\sB$ be the irreducible family of del Pezzo 3-folds
$B=G(2,5)\cap\mP^{6}$, where 
$\mP^{6}\subset \mP^9$ is transversal to $G(2,5)$. 
Let 
$$
J=\{([C^{0}_{d}],[B])\in \sC^0_d(G(2,5)) \times\sB \mid 
C^{0}_{d}\subset B\}.
$$
\noindent
The claim is equivalent to show that 
a general fiber $J\to \sB$ is irreducible. 
Since $d\leq 6$, a smooth rational curve of degree $d$ is contained in
at least a six-dimensional projective space.
Thus a general fiber of $J\to \sC^0_d(G(2,5))$ is non-empty and
irreducible.
Since $\sC^0_d(G(2,5))$ is irreducible,
it holds $J$ is irreducible.
By the argument of \cite[Proof of Theorem 3.1 p.17]{MT},
we have only to show that there is 
one particular component $\sC^{0\star}_{d}(B)$
of a general fiber $J\to \sB$ invariant under monodromy.
 
By induction let us
assume that $\sC^{0}_{d-1}(B)$ is irreducible. Let 
$[C^{0}_{d-1}]\in \sC^{0}_{d-1}(B)$ be a generic element.
The family of lines $[l]\in\sH_{1}^{B}$ which intersect 
a general element of $\sC^{0}_{d-1}(B)$ is irreducible 
by Proposition \ref{prop:Cd} (3). This implies that the family 
$\sC^{0}_{d-1,1}(B)$ of reducible curves $C^{0}_{d}=C^{0}_{d-1}\cup l$
such that 
$[C^{0}_{d-1}]\in \sC^{0}_{d-1}(B)$, $[l]\in\sH_{1}^{B}$ and
${\rm{length}}\, C^{0}_{d-1}\cap l=1$ is irreducible. Similarly to
the proof of 
Proposition \ref{prop:Cd0}, we see that
the locus containing the points corresponding to
the smoothings of curves from $\sC^{0}_{d-1,1}(B)$ is 
an irreducible component of $J$.
\end{proof}

\begin{rem}
The proof of the proposition shows that
$\sH^B_d$ is irreducible for any $d$.
\end{rem}


\begin{thebibliography}{Muk04}

\bibitem[AH95]{AH}
J.~Alexander and A.~Hirschowitz, \emph{Polynomial interpolation in several
  variables}, J. Algebraic Geom. \textbf{4} (1995), no.~2, 201--222.

\bibitem[Cor]{cornalba}
M.~Cornalba, \emph{Moduli of curves and theta-characteristics}, Lectures on
  Riemann surfaces (Trieste, 1987), World Sci. Publ., Teaneck, NJ, 1989,
  pp.~560--589.

\bibitem[DG88]{DG}
E.~Davis and A.~Geramita, \emph{Birational morphisms to $\mathbb{P}^2:$ an
  ideal-theoretic perspective}, Math. Ann. \textbf{279} (1988), 435--448.

\bibitem[DK93]{DK}
I.~Dolgachev and V.~Kanev, \emph{Polar covariants of plane cubics and
  quartics}, Adv. Math. \textbf{98} (1993), no.~2, 216--301.

\bibitem[Dol04]{doldual}
I.~Dolgachev, \emph{Dual homogeneous forms and varieties of power sums}, Milan
  J. of Math. \textbf{72} (2004), no.~1, 163--187.

\bibitem[FN89a]{FuNa}
M.~Furushima and N.~Nakayama, \emph{The family of lines on the {F}ano threefold
  ${V}\sb 5$}, Nagoya Math. J. \textbf{116} (1989), 111--122.

\bibitem[FN89b]{FuNa2}
\bysame, \emph{A new construction of a compactification of $\mathbb{C}^3$},
  Tohoku Math. J. \textbf{41} (1989), no.~4, 543--560.

\bibitem[Fuj81]{Fu2}
T.~Fujita, \emph{On the structure of polarized manifolds with total deficiency
  one, part {II}}, J. Math. Soc. of Japan \textbf{33} (1981), 415--434.

\bibitem[Gim89]{Gimi}
A.~Gimigliano, \emph{On {V}eronesean surfaces}, Nederl. Akad. Wetensch. Indag.
  Math. \textbf{51} (1989), no.~1, 71--85.

\bibitem[HH85]{HH}
R.~Hartshorne and A.~Hirschowitz, \emph{Smoothing algebraic space curves},
  Algebraic geometry, Sitges (Barcelona), 1983, Lecture Notes in Math., vol.
  1124, Springer-Verlag, Berlin-New York, 1985, pp.~98--131.

\bibitem[IK99]{IK}
A.~Iarrobino and V.~Kanev, \emph{Power sums, {G}orenstein algebra and
  determinantal loci}, Lecture Notes in Math., vol. 1721, Springer-Verlag,
  Berlin-New York, 1999.

\bibitem[Ili94]{IlievB5}
A.~Iliev, \emph{The {F}ano surface of the {G}ushel threefold}, Comp. Math.
  \textbf{94} (1994), no.~1, 81--107.

\bibitem[IR01a]{VSP3}
A.~Iliev and K.~Ranestad, \emph{Canonical curves and varieties of sums of
  powers of cubic polynomials}, J. Algebra \textbf{246} (2001), no.~1,
  385--393.

\bibitem[IR01b]{VSP2}
\bysame, \emph{${K}3$ surfaces of genus 8 and varieties of sums of powers of
  cubic fourfolds}, Trans. Amer. Math. Soc. \textbf{353} (2001), no.~4,
  1455--1468.

\bibitem[Isk77]{I1}
V.~A. Iskovskih, \emph{Fano $3$-folds $1$ $(${R}ussian$)$}, Izv. Akad. Nauk
  SSSR Ser. Mat \textbf{41} (1977), English transl. in Math. USSR Izv. 11
  (1977), 485--527.

\bibitem[Isk78]{I2}
\bysame, \emph{Fano $3$-folds $2$ $(${R}ussian$)$}, Izv. Akad. Nauk SSSR Ser.
  Mat \textbf{42} (1978), 506--549, English transl. in Math. USSR Izv. 12
  (1978), 469--506.

\bibitem[Man01]{Man}
M.~Mancini, \emph{Rational projectively {C}ohen-{M}acaulay surfaces of maximum
  degree}, Collect. Math. \textbf{52} (2001), no.~2, 117--126.

\bibitem[Mel06]{Waring}
M.~Mella, \emph{Singularities of linear systems and the {W}aring problem},
  Trans. Amer. Math. Soc. \textbf{358} (2006), no.~12, 5523--5538.

\bibitem[MM81]{MM1}
S.~Mori and S.~Mukai, \emph{Classification of {F}ano $3$-folds with $b_2 \geq
  2$}, Manuscripta Math. \textbf{36} (1981), 147--162.

\bibitem[MM85]{MM3}
\bysame, \emph{Classification of {F}ano $3$-folds with $b_2 \geq 2$, {I}},
  Algebraic and Topological Theories (Kinosaki, 1984), Kinokuniya, Tokyo, 1985,
  to the memory of Dr. Takehiko MIYATA, pp.~496--545.

\bibitem[MR05]{Abelian2}
F.~Melliez and K.~Ranestad, \emph{Degenerations of $(1,7)$-polarized abelian
  surfaces}, Math. Scand. \textbf{97} (2005), no.~2, 161--187.

\bibitem[MS01]{Abelian1}
N.~Manolache and F.-O. Schreyer, \emph{Moduli of $(1,7)$-polarized abelian
  surfaces via syzygies}, Math. Nachr. \textbf{226} (2001), 177--203.

\bibitem[MT01]{MT}
D.~Markushevich and A.~S. Tikhomirov, \emph{The {A}bel-{J}acobi map of a moduli
  component of vector bundles on the cubic threefold}, J. Algebraic Geom.
  \textbf{10} (2001), no.~1, 37--62.

\bibitem[MU83]{MU}
S.~Mukai and H.~Umemura, \emph{Minimal rational threefolds}, Algebraic geometry
  (Tokyo/Kyoto, 1982), Lecture Notes in Math., vol. 1016, Springer, Berlin,
  1983, pp.~490--518.

\bibitem[Muk92]{Mu2}
S.~Mukai, \emph{Fano $3$-folds}, London Math. Soc. Lecture Notes, vol. 179,
  Cambridge Univ. Press, 1992, pp.~255--263.

\bibitem[Muk04]{Mukai12}
\bysame, \emph{Plane quartics and {F}ano threefolds of genus twelve}, The Fano
  Conference, Univ. Torino, Turin, 2004, pp.~563--572.

\bibitem[Per02]{Per}
N.~Perrin, \emph{Courbes rationnelles sur les vari\'et\'es homog\`enes}, Ann.
  Inst. Fourier (Grenoble) \textbf{52} (2002), no.~1, 105--132.

\bibitem[Rei83]{Pagoda}
M.~Reid, \emph{Minimal models of canonical $3$-folds}, Algebraic varieties and
  analytic varieties (Tokyo, 1981), Adv. Stud. Pure Math., vol.~1,
  North-Holland, Amsterdam, 1983, pp.~131--180.

\bibitem[RS00]{VSP}
K.~Ranestad and F.~Schreyer, \emph{Varieties of sums of powers}, J. Reine
  Angew. Math. \textbf{525} (2000), 147--181.

\bibitem[Sch01]{Schr}
F-O. Schreyer, \emph{Geometry and algebra of prime {F}ano $3$-folds of genus
  $12$}, Compositio Math. \textbf{127} (2001), no.~3, 297--319.

\bibitem[Whi24]{White}
M.~P. White, \emph{On certain nets of plane curves}, Proc. Cambridge Phil. Soc.
  \textbf{22} (1924), 1--11.

\end{thebibliography}
\end{document}